\documentclass{article}
\usepackage[utf8]{inputenc}
\usepackage[a4paper]{geometry}
\usepackage{amsmath,amsthm,epsfig,latexsym,graphicx,amssymb}
\usepackage[english]{babel}

\usepackage{thmtools} 
\usepackage{dsfont} 
\usepackage[dvipsnames]{xcolor} 
\usepackage{mathtools}
\usepackage{stmaryrd} 

\usepackage[
    colorlinks = true,  
    linkcolor = PineGreen,
    citecolor = PineGreen,
    ]{hyperref}

\usepackage{tikz-cd}

\usepackage{array}

\usepackage{sectsty} 
\sectionfont{\color{Salmon!70!red}}
\subsectionfont{\color{Salmon!80!red}}

\newcolumntype{C}[1]{>{\centering\arraybackslash}m{#1}}

\usepackage{subfiles} 
\usepackage{graphpap} 
\DeclareTextFontCommand{\emph}{\color{ForestGreen}\em} 

\newcommand{\vvline}{\vspace*{\baselineskip}}

\newcommand{\RR}{\mathbb{R}} 
\newcommand{\NN}{\mathbb{N}}

\newcommand{\ZZ}{\mathbb{Z}}
\newcommand{\QQ}{\mathbb{Q}}
\newcommand{\Class}{\mathcal{C}}
\newcommand{\CoH}{\Class^{1+H}}

\newcommand{\Ca}{\Class^{\alpha}}
\newcommand{\Coa}{\Class^{1+\alpha}}
\newcommand{\Cxa}[1]{\Class^{#1+\alpha}}

\newcommand{\Cinfty}{\Class^\infty}
\newcommand{\Atlas}{\mathcal{A}}
\newcommand{\FF}{\mathcal{F}}
\newcommand{\LL}{\mathcal{L}}
\newcommand{\OO}{\mathcal{O}}

\newcommand{\vectD}[1]{\frac{\partial #1_t}{\partial t}}
\newcommand{\ivectD}[1]{\iota\left(\vectD{#1}\right)}

\newcommand{\MM}{\mathcal{M}}

\newcommand{\PP}{\mathcal{P}}
\newcommand{\PR}{{\mathcal{R}}}
\newcommand{\PC}{{\mathcal{C}}}
\newcommand{\PW}{\Sigma_\PP}
\newcommand{\PWR}{\Sigma_{\PR}}
\newcommand{\inj}{\xhookrightarrow{}}
\newcommand{\surj}{\twoheadrightarrow{}}

\newcommand{\paragraphc}[1]{\paragraph{\color{Salmon!90!red}#1}}

\DeclareMathOperator{\primR}{MinLoop_\PR}
\DeclareMathOperator{\prim}{MinLoop(\PP)}
\newcommand{\SR}{\Sigma_\PR}
\newcommand{\Su}{\Sigma_u}
\newcommand{\caract}[1]{\mathds{1}_{#1}}
\newcommand{\Rec}{\mathfrak{R}}
\newcommand{\LinkForm}{\mathfrak{L}}
\newcommand{\LinkFormi}[1]{\LinkForm_{\PR_{{#1}}}}
\newcommand{\LinkFormR}{\LinkForm_\PR}
\newcommand{\normalize}{\mathcal{N}}
\newcommand{\normR}{\normalize_\PR}

\newcommand{\norm}[1]{\normalize_{\PR_{{#1}}}}
\newcommand{\normB}[1]{\normalize_{{#1}}}

\newcommand{\sus}{M_\PP}
\DeclareMathOperator{\Sup}{Supp}

\newcommand{\TT}{\mathbb{T}}
\newcommand{\algcap}{\mathrlap{\hspace{2.3pt}\cdot}{\cap}}

\newcommand{\OS}{{\mathcal{P}_\phi}}

\newcommand{\intint}[2]{{\llbracket #1,#2 \rrbracket}}
\newcommand{\intoint}[1]{{\rrbracket -\infty,#1 \rrbracket}}
\newcommand{\intinto}[1]{{\llbracket #1,+\infty \llbracket}}

\newcommand{\wb}{\overline}
\newcommand{\wt}{\widetilde}
\newcommand{\wh}{\widehat}

\newcommand{\acts}{\curvearrowright}

\newcommand{\pot}{\mathfrak{p}}         
\newcommand{\qot}{\mathfrak{q}}
\newcommand{\Pot}{\mathfrak{P}}

\newcommand{\cohom}[1]{{[{#1}]_\phi}}     
\newcommand{\cohomG}[1]{{[{#1}]_{\phi^\Delta}}}

\newcommand{\PM}{{\mathbb{P}_\phi M}}  
\DeclareMathOperator*\lowlim{\liminf}
\DeclareMathOperator*\uplim{\limsup}

\DeclareMathOperator{\Leaf}{Leaf}

\DeclareMathOperator{\link}{link}

\DeclareMathOperator{\Int}{Int}
\DeclareMathOperator{\im}{im}
\DeclareMathOperator{\id}{id}

\DeclareMathOperator{\sign}{sign}
\DeclareMathOperator{\Leb}{Leb}

\DeclareMathOperator{\length}{length}
\DeclareMathOperator{\support}{support}

\newtheoremstyle{colorplain}%
{\topsep}   
{\topsep}   
{\itshape}  
{0pt}       
{} 
{.}         
{5pt plus 1pt minus 1pt} 
{\textbf{\textcolor{RoyalBlue}{\textbf{\thmname{#1} \thmnumber{#2}}}}\thmnote{ (#3)}}
{}

\newtheoremstyle{colorremark}%
{\topsep}   
{\topsep}   
{}  
{0pt}       
{\itshape} 
{.}         
{5pt plus 1pt minus 1pt} 
{\textcolor{Purple}{\thmname{#1} \thmnumber{#2}}\thmnote{ (#3)}}
{}

\newtheoremstyle{colordefinition}%
{\topsep}   
{\topsep}   
{}  
{0pt}       
{} 
{.}         
{5pt plus 1pt minus 1pt} 
{\textcolor{Purple}{\textbf{\thmname{#1} \thmnumber{#2}}}\thmnote{ (#3)}}
{}

\theoremstyle{colorplain}

\newtheorem{theorem}{Theorem}
\newtheorem{maintheorem}{Theorem}

\newtheorem{maincorollary}[maintheorem]{Corollary}

\newtheorem{lemma}[theorem]{Lemma}
\newtheorem{proposition}[theorem]{Proposition}
\newtheorem{corollary}[theorem]{Corollary}

\theoremstyle{colorremark}

\newtheorem{remark}[theorem]{Remark}

\newtheorem{question}[theorem]{Question}
\newtheorem{claim}{Claim}

\makeatletter 
\newenvironment{proofabstract}[1][\proofname]{
    \par
    \pushQED{\qed}%
    \normalfont \topsep6\p@\@plus6\p@\relax
    \trivlist
    \item\relax
    {\itshape\color{RoyalPurple}
    #1\@addpunct{.}}\hspace\labelsep\ignorespaces
}{%
    \popQED\endtrivlist\@endpefalse
}

\renewenvironment{proof}[1][\proofname]{
    \setcounter{claim}{0}    
    \setcounter{claimproof}{0}    
    \par
    \pushQED{\qed}%
    \normalfont \topsep6\p@\@plus6\p@\relax
    \trivlist
    \item\relax
    {\itshape\color{RoyalPurple}
    #1\@addpunct{.}}\hspace\labelsep\ignorespaces
}{%
    \popQED\endtrivlist\@endpefalse
}
\makeatother

\newcounter{claimproof} 
\newenvironment{claimproof}
  {\addtocounter{claimproof}{1} \begin{proofabstract}[Proof of Claim \the\value{claimproof}.]}
  {\end{proofabstract}}

  {\begin{proofabstract}[Sketch of proof]}
  {\end{proofabstract}}

\theoremstyle{colordefinition}
\newtheorem{definition}[theorem]{Definition}
\newtheorem{notations}[theorem]{Notations}

\title{Skewed Anosov flows are orbit equivalent to Reeb-Anosov
flows in dimension 3}
\author{Marty Théo}
\date{}

\begin{document}

\maketitle
 
\begin{abstract}
  We prove that in dimension 3, Anosov flows which are $\RR$-covered and skewed are orbit equivalent to Reeb-Anosov flows. 
  A linking number between invariant signed measures is used to characterize the existence of an invariant contact form or of a Birkhoff section with a given boundary. 
  We also prove the existence of open book decompositions with one boundary component for Reeb-Anosov flows. 
\end{abstract}

\section*{Introduction}

Geodesic flows on hyperbolic surfaces have been studied systematically using two approaches: they are classical examples of hyperbolic flows and are the Reeb flows of the natural contact structures on the corresponding manifolds. Anosov generalized hyperbolic geodesic flows in what is now known as Anosov flows. In dimension~3, other families of Anosov flows have been studied: suspension of linear Anosov diffeomorphisms of the torus, Anosov flows built from Dehn surgeries along periodic orbits (also called Fried-Goodman surgeries 
\cite{Fried83,Goodman83}), 
and Anosov flows built from gluing hyperbolic plugs 
\cite{Beguin17}. 

Except for geodesic flows, very few Anosov flows were known to be of Reeb type. Foulon and Hasselblatt 
\cite{Foulon13} 
produced Reeb-Anosov flows from surgeries on hyperbolic geodesic flows. Their construction extends a previous construction of surgery by Handel and Thurston 
\cite{Handel80}. 
Recently Salmoiraghi 
\cite{Salmoiraghi22} 
produced a larger family of Reeb-Anosov flows using bi-contact structures. His construction extends the surgery operation introduced by Goodman. 

Fenley and Barbot 
\cite{Fenley94,Barbot95b} 
independently introduced the key notion of orbit space to understand topological properties of Anosov flows. One remarkable family of Anosov flows is characterized by its orbit space: the skewed~$\RR$-covered Anosov flows. Barbot 
\cite[Theorem~A]{Barbot01}
proved that Reeb-Anosov flows are skewed~$\RR$-covered flows. The main result of the present paper is the converse implication:

\begin{maintheorem}\label{theoremSkewedAreReeb}
    Let~$\phi$ be an Anosov flow on an oriented, closed,~3-dimensional manifold. Suppose that~$\phi$ is~$\RR$-covered and positively skewed. Then~$\phi$ is orbit equivalent to a smooth Reeb-Anosov flow. Additionally if~$\alpha$ is the corresponding contact form, then~$\alpha\wedge d\alpha$ is positive.
\end{maintheorem}

An orbit equivalence between two flows is an homeomorphism exchanging the (oriented) orbits of the flows (see Section \ref{section-AnosovFlow}). 
The theorem gives an answer to Barbot/Barthelmé's conjecture~\cite{Barthelme17}. 
By combining the theorem and previous works, we obtain a series of equivalences, represented in 
Table~\ref{tableTetrachotomy}. 
Birkhoff sections and linking numbers between invariant signed measures are introduced later in the introduction. We denote by~$\MM^0_p(\phi)$ the set of null-homologous~$\phi$-invariant probability measures.

\begin{table}[h]
    \renewcommand{\arraystretch}{1.7}
    \setlength\tabcolsep{6pt}
    \hspace*{-20px}
    \begin{tabular}{c|C{3.9cm}|C{3.5cm}|C{3.9cm}|C{1.6cm}|}
        \multicolumn{1}{c}{} & \multicolumn{1}{c}{Positively twisted} & \multicolumn{1}{c}{Flat flows} & \multicolumn{1}{c}{Negatively twisted} & \multicolumn{1}{c}{Others} \\ 
        \cline{2-5}
       1 & positively skewed $\RR$-covered & $\RR$-covered with trivial bi-foliated plane & negatively skewed $\RR$-covered & non $\RR$-covered \\
        \cline{2-5}
       2 & $\exists$ Birkhoff section  with positive boundary & $\exists$ Birkhoff section with no boundary & $\exists$ Birkhoff section with negative boundary & other \\ 
        \cline{2-5}
       3 & $\exists \mu\in\MM_p^0(\phi)$ with positive linking numbers with all $\MM_p^0(\phi)$ & $\MM_p^0(\phi)=\emptyset$ & $\exists \mu\in\MM_p^0(\phi)$ with negative linking numbers with all $\MM_p^0(\phi)$ & other \\ 
        \cline{2-5}        
       4 & $\exists\psi\simeq\phi$, $\exists \alpha$ contact, $\psi$-invariant, $\alpha\wedge d\alpha >0$ & $\exists\psi\simeq\phi$, $\exists \alpha\neq 0$ $\psi$-invariant, $d\alpha=0$ &  $\exists\psi\simeq\phi$, $\exists \alpha$ contact $\psi$-invariant, $\alpha\wedge d\alpha <0$  & other \\ 
        \cline{2-5}
    \end{tabular}
    \caption{Tetrachotomy for Anosov flows on closed 3-manifolds. Each column corresponds to equivalent statements. Each line corresponds to a different point of view.}
    \label{tableTetrachotomy}
\end{table}

\newpage

\begin{maincorollary}\label{mainCorollaryEquivalence}
    In a given column in Table~\ref{tableTetrachotomy}, two cells are equivalent to each other. The columns correspond to exclusive natures for Anosov flows. See the paragraph Tetrachotomy for some explanations.
\end{maincorollary}

\paragraphc{Tetrachotomy.}
In the Table \ref{tableTetrachotomy}, we present a tetrachotomy for Anosov flows on oriented closed 3-manifolds: positively and negatively twisted flows, flat flows and others. Each column contains four equivalent properties. Line~1 corresponds to Fenley-Barbot classification of the orbit spaces, upgraded to take into account the orientation of the 3-manifold in the skewed case. Line~2 corresponds to the existence of a particular Birkhoff section. Line~3 corresponds to the existence of a particular null-homologous invariant probability measure. Line~4 corresponds to the existence of an orbit equivalent Anosov flow which admits a particular invariant differential~1-form. 

The new implications are proven at the end of 
Appendix~\ref{sectionBSexistence},
together with an alternative proof to one of Barbot's theorem (see below).

In the columns~1 and~3, the ambient manifold is oriented. Reversing the orientation exchanges the role of these two columns.
We now discuss the already existing equivalences. The tetrachotomy was first discovered by 
Barbot~\cite{Barbot95b} and Fenley~\cite{Fenley94} 
simultaneously. The equivalence between Lines~1 and~2 is proven 
in~\cite[Theorem~A]{ABM22}.

For the first column, 
Barbot~\cite[Theorem A]{Barbot01} 
proved that the~$4^{th}$ cell implies the~$1^{st}$ cell. Barbot does not specify a sign for~$\alpha\wedge d\alpha$, but it can be recovered. We give a second proof of this implication in 
Appendix~\ref{sectionBSexistence}, 
with the sign. The main theorem of this article is the converse implication: the~$1^{st}$ cell implies the~$4^{th}$ cell. We additionally introduce the~$3^{rd}$ line which is interesting on its own. The equivalence between the~$2^{nd}$ and~$3^{rd}$ cells is more elementary and in the spirit of the work of 
Ghys~\cite{Ghys09}. 

For the second column, Solodov\footnote{V.V. Solodov did not publish the proof. The main argument is that actions on~$\RR$, for which every non-trivial element acts with at most one fixed point, are conjugated to affine actions. Apply this to the stable and unstable leaf spaces in the orbit space, defined in Section \ref{section-orbitspace}} proved the equivalence between cell~1 and cell~2. The equivalence between the~$2^{nd}$ and~$3^{rd}$ cells follows from Schwartzman-Sullivan's Theorem on cross-sections 
\cite{Schwartzman57,Sullivan76}. 
The equivalence between the $2^{nd}$ and $4^{th}$ cells is well known and more elementary.

\newpage

One important corollary of the main theorem was communicated by Barthelmé and Bowden. They prove together with Mann
\cite[Theorem 1.11]{Barthelme21} 
that on a given 3-manifold, there exists at most finitely many Reeb-Anosov flows, up to orbit equivalence. 

\begin{maincorollary}\label{mainCorollaryFiniteness}
   On a given closed 3-manifold, there exists at most finitely many $\RR$-covered Anosov flows up to orbit equivalence.
\end{maincorollary}

\vline

To prove Theorem \ref{theoremSkewedAreReeb}, we roughly follow these steps. Starting from a $\RR$-covered Anosov flow, we need to find an orbit equivalent Anosov flow which preserves a smooth volume form~$V$, where~$V$ would correspond to~$\alpha\wedge d\alpha$ for a contact form~$\alpha$. Given a transitive codimension one Anosov flow (which is the case of $\RR$-covered Anosov flows in dimension 3), Asaoka constructed an orbit equivalent Anosov flow, which is additionally volume preserving~\cite{Asaoka07}. 
Asaoka uses critically the notion of Gibbs measures, which are an important family of invariant probability measures that satisfy some dynamical properties similar to the invariant volume forms.
Given a Gibbs measure, Asaoka builds a new~$\Class^{1+H}$ differential structure on the ambient manifold, for which the Gibbs measure is induced by a Hölder continuous volume form. Then he approximates the flow by smooth Anosov flows, smooth for the new differential structure, each of them preserving a smooth volume form.

For an Anosov flow in general, preserving a smooth volume form is not enough to be a reparametrization of a Reeb-Anosov flow. We need two additional conditions, which we express in terms of homology and linking number. To a $\phi$-invariant signed measure $\mu$ corresponds a homology class $\cohom{\mu}\in H_1(M,\RR)$ (called Schwartzman asymptotic cycle) which we define in 
Section \ref{sectionHomology}. 
Denote by~$\MM_s(\phi)$ the set of~$\phi$-invariant signed measures and by~$\MM^0_s(\phi)$ the set of null-homologous~$\phi$-invariant signed measures.

\begin{maintheorem}[Made precise in Theorem~\ref{theoremMeasureLinkingNumber}]\label{theoremTempMeasureLK}
    Let~$\phi$ be a transitive Anosov flow in  dimension 3. There exists a unique continuous bilinear map $\link_\phi\colon\MM^0_s(\phi)\times\MM^0_s(\phi)\to\RR$ extending the linking number for two disjoint, null-homologous knots.
\end{maintheorem}

Ghys \cite{Ghys09} 
proposed a strategy to define the linking number between invariant measures for any flow on oriented 3-manifold, but his strategy is not explicitly developed\footnote{It was privately communicated by several experts of the subject that Ghys's strategy is insufficient. Dehornoy and Boulanger are working on filling the missing parts in Ghys's strategy.}. The difficult part is to prove a uniform continuity statement, which enables to extend the linking number from a dense subset to all invariant measures. We state the required properties of the linking number in 
Section \ref{sectionLKintro}. 
The proofs (for Anosov flows only) are technical and postponed to 
Section \ref{sectionLK}. 
The use of the linking number is motivated by McDuff's criterion\footnote{
    McDuff's theorem is written in terms of structural boundaries. 
    Prasad recently interpreted McDuff's theorem in terms of linking number, for right-handed flows on rational homology spheres. Prasad defined the linking number with an invariant volume form with the equation in 
    Lemma \ref{lemmaLinkSmoothMeasure}} 
which characterize Reeb flows.

\begin{maintheorem}[McDuff~{\cite[Theorem 5.2]{McDuff87}} - Prasad~\cite{Prasad22}]\label{theoremReebLikeCondition}
    Let~$\phi$ be a smooth flow on a closed, oriented 3-manifold, and denote $X=\vectD \phi$. Then~$\phi$ is a smooth reparametrization of a Reeb flow for a contact form $\alpha$ with $\alpha\wedge d\alpha>0$, if and only if it preserves a smooth positive volume form $V$ such that $\iota_XV$ is exact and such that $V$ has positive linking number with all of $\MM^0_p(\phi)$. 
\end{maintheorem}

In the previous theorem, the linking number between $V$ (seen as an invariant measure) and $\mu\in\MM^0_s(\phi)$ is equal to $\link_\phi(V,\mu)=\int_M \iota_{X}\beta d\mu$, for any 1-one form $\beta$ satisfying $d\beta=\iota_{X}V$. We give a proof of the reformulation for our context in
Section~\ref{sectionLinkReeb}. 
Prasad's reformulation of McDuff criterion motivates the following definition. A signed measure $\mu\in\MM_s(\phi)$ is said to be Reeb-like if it satisfies $\cohom{\mu}=0$ and 
$\min_{\nu\in\MM_p^0(\phi)}\link(\mu,\nu)> 0$.
The condition in the previous theorem can be reformulated as follows: $\phi$ is a smooth reparametrization of a Reeb flow if and only if it preserves a smooth Reeb-like invariant measure (with full support). 

Our main goal is to build an orbit equivalent flow preserving a smooth invariant measure with the Reeb-like property. We first build a Reeb-like measure $\Leb_\gamma$ supported on a single periodic orbit $\gamma$, as we explain later. 
Then we approximate $\Leb_\gamma$ with a Reeb-like Gibbs measure. Finally we adapt Asaoka's strategy to regularize the Gibbs measure. We obtain a smooth Anosov flow preserving a smooth volume form, for which McDuff's criterion is satisfied.

It remains to explain the existence of the measure~$\Leb_\gamma$.
Fried~\cite[Theorem 2]{Fried83} 
proved the existence of Birkhoff sections for transitive Anosov flow. Birkhoff sections are roughly speaking transverse surfaces, bounded by periodic orbits, which intersect efficiently the orbits of the flow. The Birkhoff sections constructed by Fried bring little information on the topology of the ambient manifold and of the flow. In a previous work 
\cite[Theorem A]{ABM22}, 
we characterized positively skewed~$\RR$-covered Anosov flows using the existence of Birkhoff sections with positive boundaries. In 
Section~\ref{sectionBS}, 
we improve this result.

\begin{maintheorem}\label{theoremStrongOpenBook}
  Let~$\phi$ be a $\RR$-covered and positively skewed Anosov flow on a closed, oriented~3-dimensional manifold. Then~$\phi$ admits an embedded positive Birkhoff section with only one boundary component, whose corresponding periodic orbit has orientable stable and unstable leaves. In particular, it gives an open book decomposition of~$M$, adapted to~$\phi$.
\end{maintheorem}

Take a periodic orbit~$\gamma$ bounding a Birkhoff section as in the theorem. The invariant Lebesgue measure~$\Leb_\gamma$ on~$\gamma$ is a Reeb-like measure. Hence one can apply the above discussion on a rescaling of~$\Leb_\gamma$.

\vline 

In our strategy, we change the differential structure on the ambient manifold, and the parametrization of the flow comes as a consequence of the technique used. It comes as a surprise since the parametrization of the flow seems to play an equally important role than the smooth structure. 
Barbot introduced a notion of topologically contact Anosov
flow~\cite{Barbot01}. 
Being topologically contact is invariant under continuous conjugation, so it does not depend on the smooth structure of the ambient manifold.

\begin{question}
    Can we characterize all topologically contact parametrizations of a given Reeb-Anosov flow?
\end{question}

\vline

We end the introduction by discussing two connected theorems, which we obtain with the same techniques. Our smoothing strategy also works for Anosov flows said to by homologically full, that is every element in~$H_1(M,\ZZ)$ is the homology class of a periodic orbit. 
Sharp~\cite{Sharp93}
proved that an Anosov flow is homologically full if and only it admits a null-homologous Gibbs measure. We can smooth the Gibbs measure with Asaoka's technique and obtain. 

\begin{maintheorem}\label{theoremNullCohomVolumeFrom}
    Let~$\phi$ be a homologically full Anosov flow on a closed~3-dimensional manifold~$M$. There exists a smooth Anosov flow~$\psi$ and a smooth volume form~$V$ on~$M$ which satisfy the following:~$\psi$ is orbit equivalent to~$\phi$, it preserves~$V$, and the interior product~$\iota_{\vectD\psi}V$ is null-cohomologous.
\end{maintheorem}

In~\cite{Ghys09,Colin22}, 
a positive linking number condition is used to construct a Birkhoff section. For Anosov flows, their idea and the linking number map give a criterion for the existence of Birkhoff section bounding a specific boundary. Note that the result is not specific to transitive Anosov flows, and should hold for every flow $\phi$ admitting a honest linking function on $\MM^0_s(\phi)$.

\begin{maintheorem}\label{theoremExistenceBSwithLink}
    Let~$\phi$ be a transitive Anosov flow on a closed, oriented,~3-dimensional manifold~$M$. Let~$\Gamma$ be a collection of periodic orbits, with multiplicities in~$\ZZ$, supposed to be null-homologous in~$H_1(M,\QQ)$. Then the following statements are equivalent:
    \begin{enumerate}
        \item there exist~$n>0$ and a Birkhoff section bounded~$n\Gamma$,
        \item the invariant Lebesgue measure~$\Leb_\Gamma$ supported on~$\Gamma$ is Reeb-like.
    \end{enumerate}
\end{maintheorem} 

Notice that the theorem holds even when there exists no null-homologous probability measure. In this case both conditions are satisfied independently on~$\Gamma$. We do not know if one can always take~$n=1$ in the first assertion, assuming that~$\Gamma$ is null-homologous in~$H_1(M,\ZZ)$. We prove the theorem in 
Appendix~\ref{sectionBSexistence}.

\paragraphc{Acknowledgments.} I am grateful to M Asaoka for the discussion and help regarding the use of Gibbs measures. I thank the Max Plank Institute in Bonn for its financial support. I thank P Dehornoy and A Boulanger for the discussions on the linking number.
I also thank A Rechtman, T Barbot, M Postic and V Colin for the various discussions on the subject.

\setcounter{tocdepth}{3}  
{\hypersetup{linkcolor=black}
\tableofcontents\label{ToC}}

\section{Preliminary}

We introduce the notions appearing in the main theorem: Anosov flows, orbits space, positively skewed $\RR$-covered Anosov flows and Reeb flows.

\subsection{Anosov flows}\label{section-AnosovFlow}

Let~$M$ be a connected closed three-dimensional manifold and~$\phi$ a~$\Class^1$ flow on~$M$. The flow~$\phi$ is said \emph{Anosov} if there exists a~$\phi$-invariant splitting of~$TM$ into three line bundles~$TM=E^s\oplus X\oplus E^u$ and two real numbers~$A,B>0$ such that for one/any Riemannian norm~$ \lVert. \rVert$ on~$TM$, we have:
\begin{itemize}
    \item~$X$ is tangent to the flow,
    \item for all~$t\geq 0$,~$ \lVert d_{E^s}\phi_t \rVert\leq A\exp^{-Bt}$,
    \item for all~$t\leq 0$,~$ \lVert d_{E^u}\phi_t \rVert\leq A\exp^{-B|t|}$. 
\end{itemize}

Here~$d_{E^s}\phi_t$ and~$d_{E^u}\phi_t$ correspond to the differential~$d\phi_t$ restricted to the line bundles~$E^s$ and~$E^u$.
The flow is said \emph{transitive} if there exists an orbit of~$\phi$ dense inside~$M$. We compare three properties for Anosov flows, each of them implies that the flow is transitive: $\RR$-covered Anosov flows 
\cite{Barbot95b}, Reeb-Anosov flows~\cite{Barbot01} 
and Anosov flows admitting a Birkhoff section~\cite{Fried83}. In the rest of the article,~$M$ is supposed oriented and~$\phi$ is supposed transitive.

The bundles~$E^s$ and~$E^u$ are integrable into two 1-foliations, which we denote respectively by~$\FF^{ss}$ and~$\FF^{uu}$. They are called the strong stable and unstable foliations. The plane bundles~$E^s\oplus X$ and~$E^u\oplus X$ are also integrable in two 2-foliations, invariant by~$\phi$. They are denoted by~$\FF^s$ and~$\FF^u$ and called weak stable and unstable foliations. 

Every leaf of the foliations~$\FF^{ss}$ and~$\FF^{uu}$ is homeomorphic to~$\RR$ and each leaf of the foliations~$\FF^s$ and~$\FF^u$ is homeomorphic to either~$\RR/\ZZ \times \RR$ or an open Möbius strip if it contains a periodic orbit, or~$\RR^2$ otherwise.
Through the whole article, we use the notation~$\FF(x)$ for the leaf of the foliation~$\FF$ containing the point~$x$. We also denote by~$\Leaf(\FF)$ the set of leaves of the foliation~$\FF$, equipped with the quotient topology.

Let~$\phi$ and~$\psi$ be two flows on respectively~$M$ and~$N$. An \emph{orbit equivalence}~$h\colon (M,\phi)\to(N,\psi)$ is a homeomorphism~$h\colon M\to N$ such that for each orbit arc~$\gamma$ of~$\phi$,~$h(\gamma)$ is an orbit arc of~$\psi$ and~$h$ preserves the orientation by the flow. When~$M$ and~$N$ are oriented, we additionally suppose that~$h$ preserves the orientation. It is well known that given an Anosov flow $\phi$, any flow close enough to $\phi$, for the $\Class^1$ topology, is orbit equivalent to $\phi$, with an orbit equivalence close to the identity (see \cite[Proposition~3]{Katok91}).

Important notions in this article are invariant under orbit equivalence: the homeomorphism class of the bi-foliated plane, existence of a Birkhoff section with positive boundary and linking number between invariant signed measures.

\vvline

Let~$M'$ be the manifold given as the set of tuples~$(x,o^s_x,o^u_x)$ with~$x\in M$, and~$o^s_x,o^u_x$ are local orientations of respectively the stable and unstable foliations at~$x$. The projection~$M'\to M$ on the first coordinate is a degree 4 covering map.
Denote by~$\wh M$ one connected component of~$M'$. The projection~$\wh\pi\colon\wh M\to M$ is a covering map, which is of degree one if and only~$\phi$ has orientable stable and unstable foliations. Recall that~$M$ is supposed orientable, so~$\wh M\to M$ is of degree one or two. We call~$\wh M$ the \emph{orientations-bundle covering} of~$\phi$. We also denote by~$\wh \phi$ the lift of~$\phi$ to~$\wh M$. 
Notice that~$\wh\phi$ is Anosov, its stable and unstable foliations are the lifts of the corresponding foliations for~$\phi$, and they are orientable.

Notice that~$\wh\pi_*(\pi_1(\wh M))$ is a subgroup of~$\pi_1(\wh M)$ of index one or two, so it is normal.
We denote by~$\sign\colon\pi_1(M)\to \{-1,1\}$ the group morphism obtained as~$\pi_1(M)\surj\pi_1(\wh M)/\wh\pi_*(\pi_1(\wh M))\inj\{-1,1\}$. That is given a simple loop~$\gamma$ based on~$x$,~$\sign([\gamma])=1$ if and only if the stable and unstable foliations of~$M$ are orientable in a neighborhood of~$\gamma$. Notice that~$\sign$ is a group morphism. The map~$\sign$ is defined so that for a closed and injective loop~$\gamma$, the stable and unstable foliations of~$\phi$ are orientable on a small neighborhood of~$\gamma$ if and only if~$\sign(\gamma)=1$.

\subsection{Orbit space}\label{section-orbitspace}

Take an Anosov flow~$\phi$ on a closed 3-manifold~$M$. We denote by~$\pi_M\colon\wt M \to M$ the universal covering map of~$M$. The flow and the foliations lift to a flow~$\wt\phi$ and foliations~$\wt\FF^s$,~$\wt\FF^u$,~$\wt\FF^{ss}$ and~$\wt\FF^{uu}$ in~$\wt M$.
We call \emph{orbit space} of~$\phi$, the quotient~$\OS$ of~$\wt M$ by orbits of the flow~$\wt\phi$. Barbot \cite[Theorem 3.2]{Barbot95b} and Fenley simultaneously proved 
\cite[Proposition 2.1]{Fenley94}
that~$\OS$ is homeomorphic to a plane. We denote by~$\pi_\OS\colon\wt M\to\OS$ the projection. The foliations~$\wt\FF^s$,~$\wt\FF^u$ are saturated by~$\wt\phi$, so they project in~$\OS$ to two transverse foliations by curves, which we denote by~$\LL^s$ and~$\LL^u$. If the flow~$\phi$ is smooth, the set~$\OS$ comes with a~$\Cinfty$ structure, for which the foliations~$\LL^s$ and~$\LL^u$ are of class~$\Class^1$ 
(see~\cite[Corollary 9.4.11]{Fisher19}). 
We call bi-foliated plane of~$\phi$ the set~$(\OS,\LL^s,\LL^u)$. 

Suppose that~$M$ is oriented. We lift the orientation on~$M$ to an orientation on~$\wt M$. We equip any surface~$S\subset\wt M$ transverse to the flow~$\wt\phi$, with the orientation which satisfies that~$S$ is positively transverse to the flow. The projection~$\pi_\OS$ restricts to an immersion. We fix on~$\OS$ the orientation for which for any positively transverse surface~$S\subset\wt M$, the immersion~$S\to\OS$ is orientation preserving. 

A particular family of Anosov flow are the ones for which the space of leaves~$\Leaf(\FF^s)$ (or equivalently~$\Leaf(\LL^s)$) is homeomorphic to~$\RR$. This condition is equivalent to having~$\Leaf(\FF^u)$ and~$\Leaf(\LL^u)$  homeomorphic to~$\RR$. These flows are said to be \emph{$\RR$-covered}.

\begin{figure}
    \begin{center}
        \begin{picture}(100,32)(0,0)
        \color{green!50!black}
        \put(28.5,16.5){$\LL^s$}
        \color{red!60!black}
        \put(4,29){$\LL^u$}
        \put(0,0){\includegraphics[width=100mm]{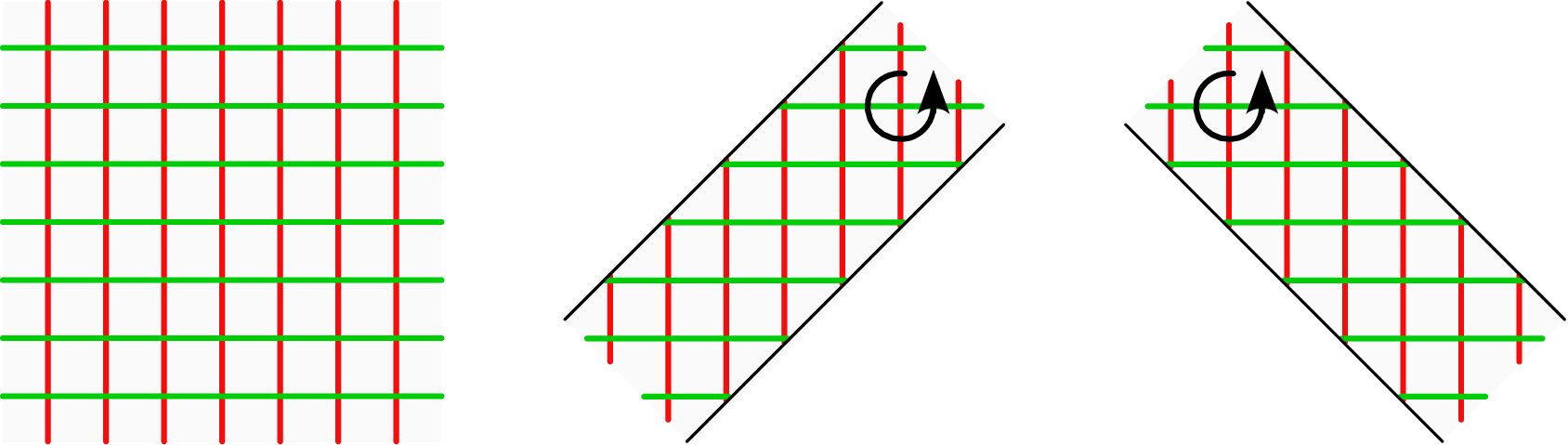}}
        \end{picture}
    \end{center}
    \caption{Classification of the bi-foliated plane of $\RR$-covered Anosov flows: suspension flows, positively skewed and negatively skewed.}
    \label{figureOrbitSpace}
\end{figure}

\begin{theorem}[Barbot \cite{Barbot95b}-Fenley {\cite[Theorem 3.4]{Fenley94}}]
    Let~$\phi$ be a $\RR$-covered Anosov flow on an oriented 3-manifold. Then there is a homeomorphism, preserving the orientation and the foliations, from the bi-foliated plane of~$\phi$ to one of the three following cases (illustrated in Figure~\ref{figureOrbitSpace}):
    \begin{enumerate}
        \item~$(\RR^2,H,V)$ where~$H$ and~$V$ are the foliations made of horizontal and vertical lines,
        \item~$(D^+,H,V)$ where~$D^+$ is the diagonal strip~$\{(x,y)\in\RR^2, |x-y|<1\}$,
        \item~$(D^-,H,V)$ where~$D^-$ is the anti-diagonal strip~$\{(x,y)\in\RR^2, |x+y|<1\}$.
    \end{enumerate}
\end{theorem}
In the theorem~$\RR^2,D^+,D^-$ are equipped with the trigonometric orientations. Fenley's and Barbot's original theorems are stated without the orientation condition: the bi-foliated plane of a $\RR$-covered Anosov flow is either trivial (case 1) or skewed (cases 2 and 3). Barbot noticed that in the skewed case, the manifold $M$ is orientable \cite[Theorem C]{Barbot95b} and so we can distinguish two sub-cases.

\begin{theorem}[Barbot {\cite[Theorem 2.5]{Barbot95b}}]\label{theoremRcoveredTransitive}
    Any $\RR$-covered Anosov flow is transitive. 
\end{theorem}

In the first case, an argument from Solodov implies that~$\phi$ is orbit equivalent to a suspension of a linear Anosov diffeomorphism on the torus. In the second and third case,~$M$ is additionally orientable, and we say that~$\phi$ is respectively positively skewed and negatively skewed (these are sometimes called twisted instead of skewed). Notice that these two cases are mutually exclusive because of the choice of the orientation. In particular reversing the orientation on~$M$ switches the roles of the cases 2 and 3.

The fundamental group~$\pi_1(M)$ acts on~$\wt M$. That action preserves the foliations of~$\wt\phi$, so it induces an action~$\pi_1(M)\acts\OS$ which preserves~$\LL^s$ and~$\LL^u$. A point~$\xi$ in~$\OS$ corresponds to the orbit~$\pi_M(\pi_\PP^{-1}(\{\xi\}))$ of~$\phi$. Additionally a point~$\xi$ invariant by a non-trivial element in~$\pi_1(M)$ corresponds to a periodic orbit of~$\phi$. More precisely, let~$\wt x\in\wt M$,~$x=\pi_M(\wt x)\in M$ and $\xi=\pi_\OS(\wt x)\in\OS$. Then~$x$ is~$\phi$-periodic if and only if there exists~$g\in\pi_1(M)\setminus\{0\}$ such that~$g\cdot \xi=\xi$. Suppose~$x$ to be periodic, then~$g$ can be taken as the homotopy class of the orbit of~$x$ in~$\pi_1(M,x)$ 
(see~\cite{Verjovsky73} 
for the non-triviality of that element). Denote by~$g\in\pi_1(M,x)$ the homotopy class of the orbit of~$x$. Then~$g$ preserves the leaves~$\LL^s(\xi)$ and~$\LL^u(\xi)$. If the stable leaf~$\FF^s(x)$ is orientable, then~$g$ preserves the orientation on~$\LL^s(\xi)$ and~$\LL^u(\xi)$. Additionally~$g$ contracts~$\LL^u(\xi)$ and expends~$\LL^s(\xi)$
(see \cite[Proposition~5]{ABM22}). 
When~$\FF^s(x)$ is not orientable, then~$g^2$ satisfies that property.

\begin{lemma}[Fenley~{\cite[Theorem 4.1]{Fenley94}}]\label{lemmaAntiHomotopicOrbit}
    Let~$\phi$ be a skewed $\RR$-covered Anosov flow. For any periodic orbit~$\gamma$ of~$\phi$, there exists a periodic orbit~$\delta$ and~$n,m\geq 1$ for which~$\gamma^n$ is anti-homotopic to~$\delta^m$.
\end{lemma}

In the lemma,~$\gamma^n$ stands for the orbit~$\gamma$ traveled~$n$ times. 

\begin{lemma}\label{lemmaHomologicallyFull}
    A skewed $\RR$-covered Anosov flow is homologically full.
\end{lemma}

\begin{proof}
    This is a consequence of two facts: the homology class of periodic orbits span~$H_1(M,\ZZ)$ 
    (see the discussion following Proposition 9 in~\cite{Parry86}), 
    and every periodic orbit has a positive multiple which is anti-homotopic to a positive multiple of another periodic orbit 
    (see Lemma~\ref{lemmaAntiHomotopicOrbit}). 
    It follows that there exists a finite family of periodic orbits~$(\gamma_i)_i$, for which the convex hull of the homology class~$[\gamma_i]\in H_1(M,\RR)$ contains~$0$ in its interior. So the flow is homologically full according to the Sharp's characterization~\cite[Theorem 1]{Sharp93}.
\end{proof}

\subsection{Reeb flows}

A \emph{contact form} on a 3-dimensional manifold~$M$ is a~$\Class^1$ differential form~$\alpha$ such that~$\alpha\wedge d\alpha$ is never zero. In particular~$\alpha\wedge d\alpha$ is a volume form on~$M$, so~$M$ is orientable. Given a contact form~$\alpha$, the plane distribution given by~$\xi=\ker(\alpha)$ is called a contact structure. 
Given a contact form~$\alpha$, there exists a unique vector field~$Y_ \alpha$ on~$M$ satisfying~$\iota_{Y_\alpha}\alpha=1$ and~$\iota_{Y_\alpha}d\alpha=0$. This vector field is called the Reeb vector field of~$\alpha$. When~$\alpha$ is of class~$\Cinfty$,~$Y_\alpha$ is also of class~$\Cinfty$, so it generates a~$\Cinfty$ flow~$\psi^\alpha$, called \emph{Reeb flow} of~$\alpha$. Flows which are both Anosov and Reeb flows are called Reeb-Anosov flows.

\begin{theorem}[Barbot~{\cite[Theorem A and Proposition 5.1]{Barbot01}}]\label{theoremBarbotTheorem}
    In dimension 3, any Reeb-Anosov flow is $\RR$-covered and skewed, and in particular transitive.
\end{theorem}

Barbot does not relate the nature positively/negatively skewed of the flow with the orientation given by the contact form. We give a new proof of his theorem in Appendix~\ref{sectionBSexistence}, which additionally gives the sign correspondence.

The contact structure of the Reeb flow of a smooth contact form~$\alpha$, that is the plane distribution given by~$\ker(\alpha)\subset TM$, is smooth, transverse to the flow and invariant by the flow. When the flow is Anosov, only one plane distribution satisfies this property, the sum~$E^s\oplus E^u$ of the strong stable and strong unstable directions of the flow. In particular this sum is smooth for Reeb-Anosov flows. 
It is known that when~$E^s\oplus E^u$ is of class~$\Class^1$, then either the flow is a constant time suspension of an Anosov diffeomorphism of the torus, or~$E^s\oplus E^u$ is the kernel of an invariant~$\Class^1$ contact form (see \cite[Theorem 3.2]{KatokS90}). In general, this plane field is only Hölder, but it can still satisfy some topological version of being of the contact type. Barbot call these flows topologically contact, and proved 
\cite{Barbot01} 
that these flows are also $\RR$-covered and skewed. 
In Section~\ref{sectionLinkReeb}, we give a condition for an Anosov flow to be of Reeb type, in terms of linking number.

section{Invariant measures}

We give our notation for invariant measures and current in this section. We denote by~$\MM_p(\phi)$ and~$\MM_s(\phi)$ the sets of~$\phi$-invariant probability measures and~$\phi$-invariant signed measures on~$M$. 
We equip these sets with the weak topology: a sequence of signed measures~$(\mu_n)_n$ converges toward~$\mu_\infty$ if for any continuous function~$f\colon M\to\RR$, the integral~$\int_Mf\mu_n$ converges toward~$\int_Mf\mu_\infty$. The set~$\MM_p(\phi)$ is compact for the weak topology.

A particular type of invariant signed measures are the one supported by finitely many periodic orbits. We call \emph{algebraic multi-orbit} an algebraic sum of the form $\Gamma=\sum_{i=1}^n a_i\gamma_i$ where $\gamma_i$ is a periodic orbit and with $a_i\in\RR$. The collection of the $\gamma_i$ is called the support of $\Gamma$.  
Given a periodic orbit $\gamma$ parametrized by $f\colon\RR/T\ZZ\to M$, $\gamma(t)=\phi_t(f(0))$, we define a Lebesgue measure $\Leb_\gamma$ on $\gamma$ by the pushforward of the Lebesgue measure on $[0,T]$ by $f$. Its total mass is the length of $\gamma$ for the flow. 
Given an algebraic multi-orbit~$\Gamma=\sum_i a_i\gamma_i$, we define the invariant signed measure~$\Leb_\Gamma=\sum_i a_i\Leb_{\gamma_i}$. These signed measures are used later in some density arguments, and to build a first Reeb-like measure (see Section \ref{sectionLinkReeb} for the definition and Lemma \ref{lemmaBStoReebLike}). 

A \emph{transverse (signed) measure}~$\mu$ is the data for any compact topologically embedded surface~$S$, topologically transverse to~$\phi$, of a signed measure~$\mu_S$ on~$S$, such that for any compact embedded surfaces~$S_1,S_2$ transverse to~$\phi$, and any holonomy map~$h\colon U_1\subset S_1\to U_2\subset S_2$ along the flow,~$h^*\mu_{S_1}$ and~$\mu_{S_2}$ coincide on~$h(U_1)$. 
Given a signed measure~$\mu\in\MM_s(\phi)$, we denote by~$\mu^\perp$ the corresponding transverse measure. It is defined by~$\mu^\perp(A)=\frac{1}{\epsilon}\mu(\phi_{[0,\epsilon]}(A))$ for any measurable set~$A\subset S$ transverse to~$\phi$, and~$\epsilon>0$ small enough (which depends on~$S$).
Denote by~$\Leb_\phi$ the unique measure on orbits of~$\phi$ in~$M$, satisfying that for every~$x\in M$ and~$\epsilon>0$ small enough,~$\phi_{[0,\epsilon]}(x)=\epsilon$. Given a transverse measure~$\mu^\perp$, we denote by~$\mu^\perp\otimes\Leb_\phi$ the invariant signed measure obtained locally as a product of~$\mu^\perp$ and~$\Leb_\phi$ along the orbits of~$\phi$. The maps~$\mu\mapsto\mu^\perp$ and~$\mu^\perp\mapsto\mu^\perp\otimes\Leb_\phi$ are continuous and inverse one to each other.

Any orbit equivalence~$h\colon(M,\phi)\to(N,\psi)$ induces a natural homeomorphism~$\mu^\perp\mapsto h^*\mu^\perp$ between the sets of transverse measures, by pushing forward a measure $\mu^\perp_S$ on a transverse surface $S\subset M$ to a measure $h^*\mu^\perp_S$ on the transverse surface $h(S)\subset N$. 
There is a natural homeomorphism~$\Theta_h\colon \MM_s(\phi)\to\MM_s(\psi)$ between the set of invariant signed measures of the two flows. For a signed measure~$\mu\in\MM_s(\phi)$,~$\Theta_h(\mu)$ is the signed measure~$h^*(\mu^\perp)\otimes_\psi\Leb_\psi$.
We introduce later the homology class and the linking number of signed measures. These notions are invariant under orbit equivalence.

\begin{remark}
    Since they do not depend on the parametrization of the flow nor on the smooth structure on the manifold, transverse measures are the natural notion to define the homology class and the linking number. For convenience, we still use invariant signed measures instead, but we keep in mind the underlining transverse signed measures.
\end{remark}

\subsection{Homology of invariant measures} \label{sectionHomology}

Here $M$ is only a compact 3-manifold and $\phi$ a flow for which $\partial M$ is invariant. 
Denote by~$n=3$ the dimension of~$M$. 
Given a $k$-form $\alpha$, a $(n-k)$-form $\beta$ and a $k$-chain $b$, we write $\alpha\cdot\beta$ and $\alpha\cdot b$ the two integrals $\int_M\alpha\wedge\beta$ and $\int_b\alpha$.

A~$k$-current is a linear form on the set~$\Omega^k(M)$ of smooth~$k$-differential form on~$M$, which is continuous for the~$\Cinfty$ topology on~$\Omega^k(M)$. See 
\cite{DeRham84}
for more background on currents.
Given a~$k$-current~$C$ and a smooth vector field~$Y$ on~$M$, we define the currents~$dC$ and~$\iota_YC$ to be the~$k-1$ and~$k+1$ currents given respectively by:
\begin{itemize}
    \item~$dC(\alpha)=C(d\alpha)$,
    \item~$\iota_YC(\alpha)=C(\iota_Y\alpha)$.
\end{itemize}

Given an invariant signed measure~$\mu$, the map~$\Omega^0(M)=\Cinfty(M,\RR)\to\RR$ given by~$f\mapsto \int_Mf(x)\mu(x)$ is a 0-current, which we denote by~$C_\mu$. In particular if~$X$ is the vector field generating the flow~$\phi$,~$\iota_XC_\mu$ is a 1-current, which is canonically associated to the transverse measure~$\mu^\perp$. 

For any closed $k$-current $C$, that is satisfying $dC=0$, there is a natural element $[C]$ in the dual of $H^1(M,\RR)$ given by $[C]\cdot [\alpha]=C(\alpha)$. The dual of $H^1(M,\RR)$ is $H_1(M,\RR)$, which is naturally isomorphic to $H^2(M,\partial M,\RR)$ (or simply $H^2(M,\RR)$ when $M$ is closed) by Lefschetz-Poincaré duality. We denote by $[C]\in H_1(M,\RR)\simeq H^2(M,\partial M,\RR)$ the homology of a closed current $C$. Depending on the
context, one thinks about $[C]$ as a homology element or as a cohomology element.

\begin{lemma}\label{lemmaInvarianceImpliesClosed}
    For any invariant signed measure~$\mu\in\MM_s(\phi)$, we have~$d(\iota_XC_\mu)=0$. 
\end{lemma}

\begin{proof}
    Let~$\mu\in\MM_s(\phi)$ be an invariant signed measure. For any smooth function~$f\colon M\to\RR$, we have~$C_\mu(\phi_t^*f)=C_\mu(f)$. It follows that~$\LL_XC_\mu(f)=C_\mu(\LL_Xf)=0$, where~$\LL_X$ is the Lie derivative along~$X=\vectD{\phi}$. Therefore we have~$d(\iota_XC_\mu)=0$.
\end{proof}

\begin{definition}
We denote by~$\cohom{\mu}\in H_1(M,\RR)$ the homology class of~$\iota_XC_\mu$, and by~$\MM^0_s(\phi)$ the set of signed measures~$\nu\in\MM_s(\phi)$ which are null-homologous.
\end{definition}

\begin{lemma}\label{lemmaCohomologyContinuity}
    The map~$\mu\in\MM_s(\phi)\to \cohom{\mu}\in H_1(M,\RR)$ is continuous.
\end{lemma}

\begin{proof}
    The homology class of a signed measure~$\mu\in\MM_s(\phi)$ is determined by the integrals~$\int_M\iota_X\alpha_kd\mu$ for a finite family~$(\alpha_k)_k$ of closed 1-form for which~$([\alpha_k])_k$ spans~$H^1(M,\RR)$. Since the functions~$\iota_X\alpha_k$ are continuous, the maps~$\mu\mapsto \int_M\iota_X\alpha_kd\mu$ are also continuous. Hence the map~$\mu\mapsto\cohom{\mu}$ is continuous.
\end{proof}

\begin{lemma}\label{lemmaCohomologyMeasureVolume}
    Assume $M$ closed. Let~$V$ be a $\Class^1$ invariant volume form and~$\mu_V$ be the invariant measure induced by~$V$. 
    Then~$\iota_XV$ is closed and~$\cohom{\mu_V}$ and $[\iota_XV]$ are Poincaré dual.
\end{lemma}

\begin{proof}
    For any closed 1-form~$\alpha$, we have~$(\iota_X\alpha) V = \alpha\wedge(\iota_XV)$. Hence the following holds
    $$\cohom{\mu}\cdot[\alpha] 
            = \int_M(\iota_X\alpha) V  
            = \int_M\alpha\wedge(\iota_XV) 
            = [\iota_XV]\cdot[\alpha]$$
    The lemma follows from Poincaré duality.
\end{proof}

An algebraic multi-orbit $\Gamma$ can be seen as a singular cycle for the singular homology, and thus induces a homology class $[\Gamma]\in H_1(M,\RR)$.

\begin{lemma}\label{lemmaMeasureCohomology}
    For any algebraic multi-orbit $\Gamma$, we have $\cohom{\Leb_\Gamma}=[\Gamma]$.
\end{lemma}

\begin{proof}
    For any closed 1-form~$\beta$ and any periodic orbit~$\gamma$ of period~$T$, we have:$$\cohom{\Leb_\gamma}\cdot[\beta]=\int_M\iota_X\beta\Leb_\gamma=\int_0^T\iota_X\beta(\phi_t(x))dt=\int_\gamma\beta=[\beta]\cdot[\gamma]$$
\end{proof}

Denote by $h_*\colon H_1(M,\RR)\to H_1(N,\RR)$ the push forward map.

\begin{proposition}\label{propFluxInvariance}
    Suppose $\phi$ Anosov and transitive. Let~$h\colon (M,\phi)\to (N,\psi)$ be an orbit equivalence. Then for any~$\mu$ in~$\MM_s(\phi)$, we have~$[\Theta_h(\mu)]_\psi=h_{*}\cohom{\mu}$.
\end{proposition}

In other words, the homology class is a topological data which depends only on the transverse measure corresponding to~$\mu$. 
We give a proof adapted to Anosov flow, but the result remains true for all flows.

\begin{proof}
    When $\phi$ is Anosov, it is well known that the set of signed measures supported on finitely many orbits is dense in $\MM_s(\phi)$. We give a stronger result in the next lemma.
    Let $\mu\in\MM_s(\phi)$ be a signed measure supported on finitely many periodic orbits.     
    Lemma \ref{lemmaMeasureCohomology} implies that $[\Theta_h(\mu)]_\psi=h_{*}\cohom{\mu}$.    
    So by continuity and density, $[\Theta_h(\mu)]_\psi=h_{*}\cohom{\mu}$ holds for all signed measures $\mu$ in~$\MM_s(\phi)$. 
\end{proof}

Finally we give a statement whose proof is delayed to Section \ref{sectionLK}. 

\begin{lemma}\label{lemmaFinitBoundedMeasureDense}
    Given a signed measure $\mu\in\MM_s(\phi)$ and a finite union of periodic orbit $\Delta\subset M$, there exists a sequence of algebraic multi-orbits $\Gamma_n$ with rational coefficients, such that $\Leb_{\Gamma_n}$ converges toward $\mu$, and so that for all $n$:
    \begin{itemize}
        \item the support of $\Gamma_n$ is disjoint from $\Delta$,
        \item $\cohom{\Leb_{\Gamma_n}}=\cohom{\mu}$.
    \end{itemize}
    Additionally if $\Delta$ is empty and $\mu$ is a probability measure, $\Gamma_n$ can be taken so that $\Leb_{\Gamma_n}$ is a probability measure.
\end{lemma}

The lemma follows from the 
Lemmas \ref{lemmaDensityFiniteOrbitPW} and \ref{lemmaDensityFiniteMeasureDisjoint},
which are similar but for a symbolic dynamic system.

\subsection{Linking number of invariant measures} \label{sectionLKintro}

The proof of Theorem \ref{theoremSkewedAreReeb} relies on some results regarding linking numbers. The proofs are postponed to
Section \ref{sectionLK}.

Fix one transitive Anosov flow~$\phi$ on a closed oriented 3-dimensional manifold~$M$. Consider~$\Gamma_1$ and~$\Gamma_2$ two disjoint knots (not necessarily connected), supposed to be null-homologous (with coefficients in $\ZZ$). The linking number between $\Gamma_1$ and $\Gamma_2$ is classically defined as $$\link(\Gamma_1,\Gamma_2)=S_1\algcap \Gamma_2$$ where $S_1$ is a 2-chain bounded by $\Gamma_1$. The linking number does not depend on the choice of $S_1$ and is equal to $\Gamma_1\algcap S_2$ for any 2-chain $S_2$ bounded by $\Gamma_2$. 

\begin{theorem}\label{theoremMeasureLinkingNumber}
    Let~$\phi$ be a transitive Anosov flow on an oriented 3-manifold. There exists a unique continuous bilinear map~$\link_\phi\colon\MM^0_s(\phi)\times\MM^0_s(\phi)\to\RR$ which satisfies~$\link_\phi(\Leb_{\Gamma_1},\Leb_{\Gamma_2})=\link(\Gamma_1,\Gamma_2)$ for all pairs of disjoint and null-homologous knots.
    Additionally the map~$\link_\phi$ is symmetric.
\end{theorem}

\begin{remark}
    For two exact 2-forms $\alpha_1, \alpha_2=d\beta_2$, a similar linking number can be defined by $\link(\alpha_1,\alpha_2)=\int_M\alpha_1\wedge\beta_2$. For an exact 2-form $\alpha=d\beta$ and a null-homologous knot $L=\partial S$, similarly their linking number is defined by $\link(\alpha,L)=\int_S\alpha=\int_L\beta$. Exact 2-forms and null-homologous knots which are $\phi$-invariant corresponds to specific null-homologous invariant measures. The linking number in Theorem \ref{theoremMeasureLinkingNumber} is equal to the linking numbers given above in both cases (see Lemma \ref{lemmaLinkSmoothMeasure}).
\end{remark}

Denote by~$\MM^0_p(\phi)$ the set of null-homologous invariant probability measures. As a consequence of the continuity and the compactness of~$\MM^0_p(\phi)$ (see \cite[Corollary 13.30]{Klenke2008} for the compactness), we obtain the following.

\begin{corollary}\label{corollaryContinuityMinLK}
    Suppose that~$\MM_p^0(\phi)$ is not empty. The map~$\MM^0_p(\phi)\to\RR$ defined as follows is well-defined and continuous $$\mu\mapsto\min_{\nu\in\MM_p^0(\phi)}\link_\phi(\mu,\nu)$$
\end{corollary}

From the density of finitely supported signed measures (Lemma \ref{lemmaFinitBoundedMeasureDense}) and the continuity of the linking number, we deduce the following. 

\begin{proposition}\label{lemmaInvarianceLKbyTopEquiv}
    For any orbit equivalence $h\colon (M,\phi)\to(N,\psi)$ between two Anosov flows on oriented closed 3-dimensional manifolds, we have~$\link_\phi=\link_\psi\circ(\Theta_h,\Theta_h)$.
\end{proposition}

We need a last lemma to compute linking number with volume forms. Given a $\phi$-invariant volume form $V$, denote by~$\nu_V\in\MM^0_s(\phi)$ the $\phi$-invariant measure induced by $V$. 

\begin{lemma}\label{lemmaLinkSmoothMeasure}
    Let $V$ be a $\phi$-invariant volume form, and assume that $\iota_X V=d\alpha$ for some smooth 1-form $\alpha$. Then for any invariant signed measure $\mu\in\MM^0_s(\phi)$, we have $$\link_\phi(\nu_V,\mu)=\int_M\iota_X\alpha d\mu$$
\end{lemma}

The proof is delayed to the very end of Section \ref{sectionLK}.

\subsection{Reeb-like measures}\label{sectionLinkReeb}

We prove Theorem~\ref{theoremReebLikeCondition}: 
a smooth flow is a reparametrization of a Reeb flow if and only if it preserves a smooth volume form which satisfies the Reeb-like condition.

\begin{definition}
    A signed measure $\mu\in\MM_s(\phi)$ is said to be \emph{Reeb-like} if it satisfies $\cohom{\mu}=0$ and 
    $$\min_{\nu\in\MM_p^0(\phi)}\link(\mu,\nu)> 0$$
\end{definition}

We also say that an invariant volume form is Reeb-like if the corresponding invariant measure is Reeb-like. The Reeb-like condition is invariant under orbit equivalence of $\phi$. That is if $h\colon (M,\phi)\to(N,\psi)$ is an orbit equivalence between two transitive Anosov flows, then a signed measure $\mu\in\MM_s(\phi)$ is Reeb-like for $\phi$ if and only if $\Theta_h(\mu)$ is Reeb-like condition for $\psi$ 
(see Propositions~\ref{propFluxInvariance} and~\ref{lemmaInvarianceLKbyTopEquiv}).

\begin{remark}[Flux]
    One can interpret~$\cohom{\mu}$ as the flux of~$\phi$ for the measure~$\mu$. That is~$\cohom{\mu}\cdot[S]$ corresponds to the quantity of flow going through a closed surface~$S$, relatively to~$\mu$.
    In~\cite{Colin08}, 
    the authors define a notion of flux for area-preserving diffeomorphisms. They establish that an area-preserving diffeomorphism (on a surface with boundary) is the first-return map of a Reeb flow if and only if it has zero flux. It is an exercise to verify that the zero flux condition in their article is equivalent to the null-homological condition in the Reeb-like condition. Our linking number condition corresponds to a control of the contact structure on the boundary of the 3-manifold on which the Reeb flow lives.
\end{remark}

\begin{lemma}\label{lemmaNullHomologousIntersectionIsEnough}
    Let $E$ be a finite dimensional vector space and $F\subset E\times\RR$ be a compact convex set such that $F\cap (0\times(-\infty,0])=\emptyset$. Then there exist $\eta>0$ and a linear map $f\colon E\to\RR$ satisfying $$F\subset \big\{(x,t)\in E\times\RR,t\geq f(x)+\eta\big\}$$
\end{lemma}

\begin{proof}
    Since $F$ is convex and compact, there exists a hyperplane $H\subset E\times\RR$ separating~$F$ and~$(0,0)$. Additionally $H$ can be chosen so that $H\cap(0\times\RR)=\{(0,c)\}$ for some $c>0$. Then there exists a linear map $f\colon E\to\RR$ satisfying that $(x,y)\in E\times\RR$ lies in $H$ if and only if $y=f(x)+c$, and $f(0,1)>0$. Then for all $(x,y)\in F$, we have $y>f(x)+c$.
\end{proof}

Recall the statement Theorem~\ref{theoremReebLikeCondition}: a smooth flow is a reparametrization of a Reeb flow if and only if it preserves a smooth volume form with the Reeb-like property. McDuff wrote her theorem using a different vocabulary. Prasad interpreted the theorem for flows on homology-spheres only. His proof of the reformulation also works on non-homology spheres, as long as we add the null-homologous conditions for measures. For completeness, we sketch a proof in the general case.

\begin{proof}[Sketch proof of Theorem~\ref{theoremReebLikeCondition}]
    Suppose that a smooth reparametrization $\psi$ of $\phi$ is the Reeb flow of a contact form $\alpha$ with $\alpha\wedge d\alpha>0$. 
    Denote by $X=\vectD{\phi}$, by $Y=\vectD{\psi} = e^fX$ for some smooth function $f\colon M\to\RR$, and by $V=e^f\alpha\wedge d\alpha$. Notice that~$\iota_XV=d\alpha$ is exact. So $\LL_XV=0$
    and $V$ is a $\phi$-invariant positive volume form.
    According to 
    Lemma~\ref{lemmaLinkSmoothMeasure}, 
    the linking number between $V$ and any probability measure $\mu\in\MM^0_p(\phi)$ is equal to
    $$\link_\phi(V,\mu)=\int_M \iota_X\alpha d\mu=\int_Me^{-f}d\mu>0$$

    Therefore the volume form $V$ is Reeb-like for the flow $\phi$. 

    Suppose now that $\phi$ admits an invariant volume form $V$, supposed smooth, positive and Reeb-like. By assumption, there exists a smooth differential 1-form~$\beta$ on $M$ such that $d\beta =\iota_XV$. We transform $\beta$ in a contact form in two steps. 
    By the assumption on the linking number, the convex $$F=\Big\{\Big(\cohom{\mu},\int_M\iota_X\beta d\mu\Big)\in H_1(M,\RR)\times \RR,\mu\in\MM_p(\phi)\Big\}$$
    satisfies $F\cap (0\times(-\infty,0])=\emptyset$. Therefore
    Lemma~\ref{lemmaNullHomologousIntersectionIsEnough} applied to the $F$ yields a closed 1-form $\delta$ and $c>0$ satisfying that for all $\mu\in\MM_p(\phi)$, we have 
    $$\int_M\iota_X(\beta+\delta)d\mu=\int_M\iota_X\beta d\mu +\cohom{\mu}\cdot[\delta] \geq c$$
    Then by averaging $\beta+\delta$ along the flow for a time $T>0$ large enough, one obtains a differential 1-form 
    $$\alpha=\frac{1}{T}\int_0^T\phi_t^*(\beta+\delta)dt$$
    which satisfies $\iota_X\alpha>0$ and $d\alpha=\iota_XV$. It follows that $\alpha\wedge d\alpha=(\iota_X\alpha)V>0$ and that the flow generated by $\frac{1}{\iota_X\alpha}X$ is a Reeb flow.
\end{proof}

\section{Markov partition} \label{sectionMarkovPartition}

Let $\phi$ be a transitive Anosov flow. We define Markov partitions and prove technical lemmas used 
in Sections \ref{sectionBS} and \ref{sectionLK}. 
Markov partitions are well studied, some standard references can be found in \cite{Bowen08} (for diffeomorphisms) and \cite{Ratner69}. We call \emph{Markov cuboid} for $\phi$ a compact subset~$\PC\subset M$ for which there exists a homeomorphism $f\colon[0,1]^3\to\PC$, such that for all $x,y$ in $[0,1]$, we have: 
\begin{itemize}
    \item $f(x\times[0,1]^2)$ is included in a weak unstable leaf,
    \item $f([0,1]\times y\times[0,1])$ is included in a weak stable leaf.
\end{itemize}
Similarly a \emph{Markov rectangle} $\PR$ is the image of a topological embedding $f\colon[0,1]^2\to \PR\subset M$, topologically transverse to $\phi$, which sends the horizontal and vertical foliations in $[0,1]^2$ to the weak stable and unstable foliations in $M$. In this section, we want to differentiate the cuboids from the rectangles. In the following sections, only the rectangles are used.

Take a Markov cuboid $\PC$ and a map $f\colon[0,1]^3\to\PC$ given as above, and suppose that for every $x,y\in[0,1]$, we have $t\in[0,1]\mapsto f(x,y,t)$ is orientation preserving, for the orientation given by the flow. We define the sets:
\begin{itemize}
    \item $\partial^{s}\PC=f([0,1]\times \{0,1\}\times[0,1])$ the stable boundary of $\PC$,
    \item $\partial^{u}\PC=f(\{0,1\}\times[0,1]\times[0,1])$ the unstable boundary of $\PC$,
    \item $\partial^{su}\PC=\partial^s\PC\cup\partial^u\PC$,
    \item $\partial^-\PC=f([0,1]^2\times 0)$, $\partial^+\PC=f([0,1]^2\times 1)$, which are Markov rectangles. 
    \item we denote by $\Int\partial^\epsilon\PC$ for $\epsilon\in\{s,u,su,-,+\}$, the image by $f$ of the same sets, where $[0,1]$ is replaced by $(0,1)$. That is $\Int\partial^\epsilon\PC$ is the interior of the rectangles corresponding to $\partial^\epsilon\PC$.
\end{itemize}

The weak stable and unstable foliations on $M$ induce stable and unstable foliations on any Markov rectangle $\PR$ and cuboid $\PC$. Denote by $W_\PR^s,W_\PR^u,W_\PC^s,W_\PC^u$ the stable and unstable foliations on the rectangle $\PR$ and the cuboid $\PC$. Given a Markov rectangle $\PR$, a subset $U\subset\PR$ is said to be a \emph{vertical} (resp. \emph{horizontal}) sub-rectangle of $\PR$ if it is connected and is the union of unstable (resp. stable) leaves on $\PR$. 
A \emph{Markov partition}~$\PP$ is a finite family of Markov rectangles $\{\PR_1,\hdots,\PR_n\}$ for which there exist some Markov cuboid $\{\PC_1,\hdots,\PC_n\}$ that satisfies:
\begin{enumerate}
    \item $\PR_i=\partial^-\PC_i$,
    \item $\cup_i\PC_i=M$,
    \item for any distinct $i,j$, the interior of $\PC_i$ and $\PC_j$ are disjoint,
    \item for any $i,j$, either the intersection $(\Int\partial^+\PC_i)\cap\Int\PR_j$ is empty, or $\PC_i\cap\PC_j$ is included in $\partial^+\PC_i\cap\PR_j$, and it is a horizontal sub-rectangle of $\partial^+\PC_i$ and a vertical sub-rectangle of $\PR_j$. 
\end{enumerate}

The Markov cuboids in the definition are unique. If $\PR$ is a Markov rectangle of a Markov partition, we denote by $\PC_\PR$ the corresponding cuboid.

Any transitive Anosov flow admits a Markov partition (see \cite{Ratner69} or \cite[Theorem~2.1]{Ratner73}).
Our definition of Markov partition is more restrictive than the usual definition, in the sense that a Markov cuboid $\PC$ has a connected boundary $\partial^+\PC$. By taking Markov sub-partitions, any Markov partition in the usual sense can be transformed into a Markov partition in our sense. We define $\partial^s\PP,\partial^u\PP,\partial^{su}\PP$ as the union of the corresponding boundary of the Markov cuboids of the Markov partition $\PP$. It follows from the definition that for all $t\geq 0$, we have~$\phi_t(\partial^s\PP)\subset\partial^s\PP$ and $\phi_{-t}(\partial^u\PP)\subset\partial^u\PP$.

\paragraphc{Symbolic dynamics.}
Here we fix a Markov partition $\PP=\{\PR_1,\hdots,\PR_n\}$. Denote by $A\in M_n(\ZZ)$ the transition matrix of $\PP$, that is the (invertible) matrix which contains a 1 in position $(i,j)$ if $\Int(\partial^+\PC_i)\cap\Int\PR_j$ is non-empty, and a zero otherwise.
We call \emph{$\PP$-word} any sequence $u\colon I\subset\ZZ\to \PP$ 
satisfying $A^{j-i}_{u_i,u_j}>0$ whenever $i$ and $j$ are in $I$. 
The set of bi-infinite $\PP$-words (that is sequences $u\colon\ZZ\to\PP$) is denoted by $\PW$. For two $\PP$-words $u,v\colon I\subset\ZZ\to\PP$, we denote by $d_\PP$ the distance between $u$ and $v$, defined by:
$$d_\PP(u,v)=\sum_{\substack{i\in I\\ u_i\neq v_i}}2^{-|i|}$$

It induces on $\PW$ the topology for which cylinders form a basis of the topology. Here cylinders are the subsets of the form $\{u\in\PW, u_i=v_i\text{ for all }i\in I\}$, where $I$ is a finite subset of $\ZZ$ and $v\colon I\to\PP$ is a $\PP$-word.
We denote by $\sigma$ the shift map on $\PW$, given by $\sigma(u)_i=u_{i+1}$. Notice that~$\sigma$ is 2-Lipschitz. 

Like Anosov flows, $\sigma$ admits a stable foliation and an unstable foliation, denoted by $W^s$ and $W^u$. Here $W^s(v)$ (resp. $W^u(v)$) is the set of sequence $w\in\PW$ with $w_i=v_i$ for all $i$ large enough  (resp. for all $i\leq i_0$ for some $i_0\in\ZZ$). The shift map divides the distances by 2 in each stable leaf, and multiplies the distances by 2 on each unstable leaf. We additionally call the local stable and unstable leaves of $v$ the sets defined by 
$$W^s_l(v)=\{w\in\PW,w_i=v_i\text{ for all }i\geq 0\}$$
$$W^u_l(v)=\{w\in\PW,w_i=v_i\text{ for all }i\leq 0\}$$

A \emph{cyclic $\PP$-word} is a map $u\colon \ZZ/n\ZZ\to\PP$ such that we have $A_{u_i,u_{i+1}}>0$ for all $i\in\ZZ/n\ZZ$. We call $n$ the \emph{length of $u$}. We say that $u$ is \emph{primitive} if it is not $k$-periodic for any $k<n$. To a cyclic $\PP$-word $u$ corresponds a $\PP$-word $\wb u\colon\ZZ\to\PP$ given by $\wb u_i=u_{(i\bmod n)}$. 
For two cyclic $\PP$-words $u,v$ of length $n,m$ such that $u_0=v_0$, we denote by $uv$ or $u\circ v$ the concatenation of $u$ and $v$, which is the cyclic $\PP$-word of length $n+m$ given for $i\in\intint{0}{n-1}$ by $(u v)_i=u_i$, and for $i\in \intint{n}{m+n-1}$ by $(u v)_i=v_{i-n}$. 

\paragraphc{Encoding the dynamics.} Given a point $x\in M$, a $\PP$-word $u\in\PW$ and an increasing sequence~$(t_n)_n$ in $\RR$ with $t_0\leq0$ and $t_1>0$, we say that~$x$ admits an \emph{itinerary} $(u_n,t_n)_{n\in\ZZ}$ if for any $n\in\ZZ$, the point~$\phi_{t_n}(x)$ belongs to the rectangle $u_n$, and $\phi_{[t_n,t_{n+1}]}(x)$ lies in the corresponding cuboid (in particular $x$ belongs to the cuboid $\PC_{u_0}$). We also say that~$x$ admits an itinerary along~$u$. It implies that $t_n$ converges to $-\infty$ and $+\infty$ when $n$ goes to~$-\infty$ and respectively~$+\infty$. Every point $x\in M$ admits an itinerary, but there is no uniqueness when the orbit of $x$ intersects the stable/unstable boundary of the Markov partition. Additionally every $u$ in $\PW$ corresponds to an itinerary of a point in $M$. Given a set $I\subset\ZZ$ and a $\PP$-word $u\colon I\to\PP$, we say that $x\in M$ admits a \emph{short itinerary} along $u$ if $x$ admits an itinerary of the form $(v_n,t_n)_n$ where $v$ is a bi-infinite $\PP$-word extending $u$. When $u$ a a cyclic $\PP$-word, we similarly say that $x\in M$ admits a short itinerary along $u$ if it admits a short itinerary (for the previous definition) along the (not-cyclic) $\PP$-word $u_{|\intint{0}{|u|}}$ (which starts and ends by the rectangle $u_0=u_{|u|}$).

Let $\PR$ be a Markov rectangle of a Markov partition. Take $n\geq 1$, a finite $\PP$-word $u\colon\intint{0}{n}\to\PP$ satisfying $u_0=\PR$, and denote by $\PR_u$ the set of points in $\PR\subset M$ admitting a short itinerary along~$u$. 
We denote by $T^u\colon \PR_u\to\RR^*_+$ the map obtained the following way. For $x\in\PR_u$, the point $x$ admits an itinerary of the form $(w_i,t_i)_{i\in\ZZ}$ such that $w_{|\intint{0}{n}}=u$ and $t_0=0$. Set $T^u(x)=t_{n}$, so that $\phi_{T^u(x)}(x)$ lies in the Markov rectangle $u_{n}\cap\partial^+\PC_{u_{n-1}}$. The element $T^u(x)$ does not depend on the specific choice of $w$. We also denote by $\sigma^u\colon \PR\to u_{n}\cap\partial^+\PC_{u_{n-1}}$ the map given by $\sigma^{u}(x)=\phi_{T^u(x)}(x)$.

\begin{lemma}\label{lemmaShortItinerary}
    The following holds:  
    \begin{itemize}
        \item $\PR_u$ is a closed horizontal sub-rectangle of $\PR$. 
        \item $\sigma^u(\PR_u)$ is a closed vertical sub-rectangle of $u_n$.
    \end{itemize}
\end{lemma}

\begin{proof}
    When $n=1$, it is a consequence of the definition of Markov partitions. When $n\geq 2$, it can be verified by induction on $n$.
\end{proof}

We denote by $\pi_\PP\colon \PW\to M$ the \emph{encoding map} given by $\pi_\PP(u)=x$, where $x$ is the unique point in $M$ which admits an itinerary of the form $(u_n,t_n)_n$ for some sequence $(t_n)_n$ satisfying $t_0=0$. In particular $\pi_\PP(u)$ lies in the rectangle $\partial^-u_0$. The map $\pi_\PP$ is well-defined, continuous, and sends the foliations $W^s_l,W^u_l$ to the foliations $W_\PR^s,W_\PR^u$ in each rectangle $\PR\in\PP$.

\subsection{Suspension of a Markov partition}\label{sectionMarkovPartitionSuspension}

Fix a Markov partition $\PP$ of $\phi$. Given $u$ in $\PW$, denote by $T_\PP(u)\in\RR^+$ the first-return time from $\pi_\PP(u)\in \partial^-u_0$ to $\pi_\PP\circ\sigma(u)\in\partial^+u_0\cap\partial^-u_1$ for $\phi$. The map $T_\PP:\PW\to\RR^+$ is continuous since~$\pi_\PP$ is continuous. We denote by $(\sus,\phi^\PP)$ the \emph{suspension flow} for the map $\sigma$ and the time function~$T_\PP$. That is $\sus$ is the compact quotient $\PW\times\RR/\simeq$ where $\simeq$ is the relation equivalence induced by $(u,s)\simeq(\sigma(u),s-T_\PP(u))$. The flow $\phi^\PP$ is defined by  $\phi^\PP_t(u,s)=(u,t+s)$. 

\begin{lemma}\label{lemmaSymboliqueSuspension}
    The map $h_\PP\colon \sus\to M$ given by $h_\PP(u,s)=\phi_s(\pi_\PP(u))$ is well-defined and is a semi-conjugacy between $\phi^\PP$ and $\phi$. That is $h_\PP\circ \phi^\PP_t=\phi_t\circ h_\PP$.
    Additionally it is continuous, and is injective on the set of point $x$ satisfying that the orbit of $h_\PP(x)$ is disjoint from $\partial^{su}\PP$.
\end{lemma}

See \cite[Theorem 6.6.5]{Fisher19} for a proof of Lemma \ref{lemmaSymboliqueSuspension}.
Denote by $\MM_s(\sigma)$ and $\MM_s(\phi^\PP)$ the set of signed measures on $\PW$ and $\sus$ which are invariant by $\sigma$ and $\phi^\PP$. For any signed measure $\mu\in\MM_s(\sigma)$, we denote by $\mu\otimes\Leb_\PP\in\MM_s(\phi^\PP)$ the signed measure which coincide locally with the product of $\mu$ and the Lebesgue measure on $\RR$. We define $\Theta_{\PP}\colon\MM_s(\sigma)\to \MM_s(\phi)$, the surjective map given by $\Theta_{\PP}(\mu)=h_\PP^*(\mu\otimes\Leb_\PP)$. The map $\Theta_\PP$ is used in 
Section \ref{sectionLK} 
to relate a linking number between $\phi$-invariant signed measures with a linking number between $\sigma$-invariant signed measures.

\begin{figure}
    \begin{center}
        \begin{picture}(147,70)(0,0)
          \color{green!50!black} 
          \put(28.5,16.5){$\LL^s$} 
          \color{red!60!black} 
          \put(4,29){$\LL^u$}
          \put(0,0){\includegraphics[height=60mm]{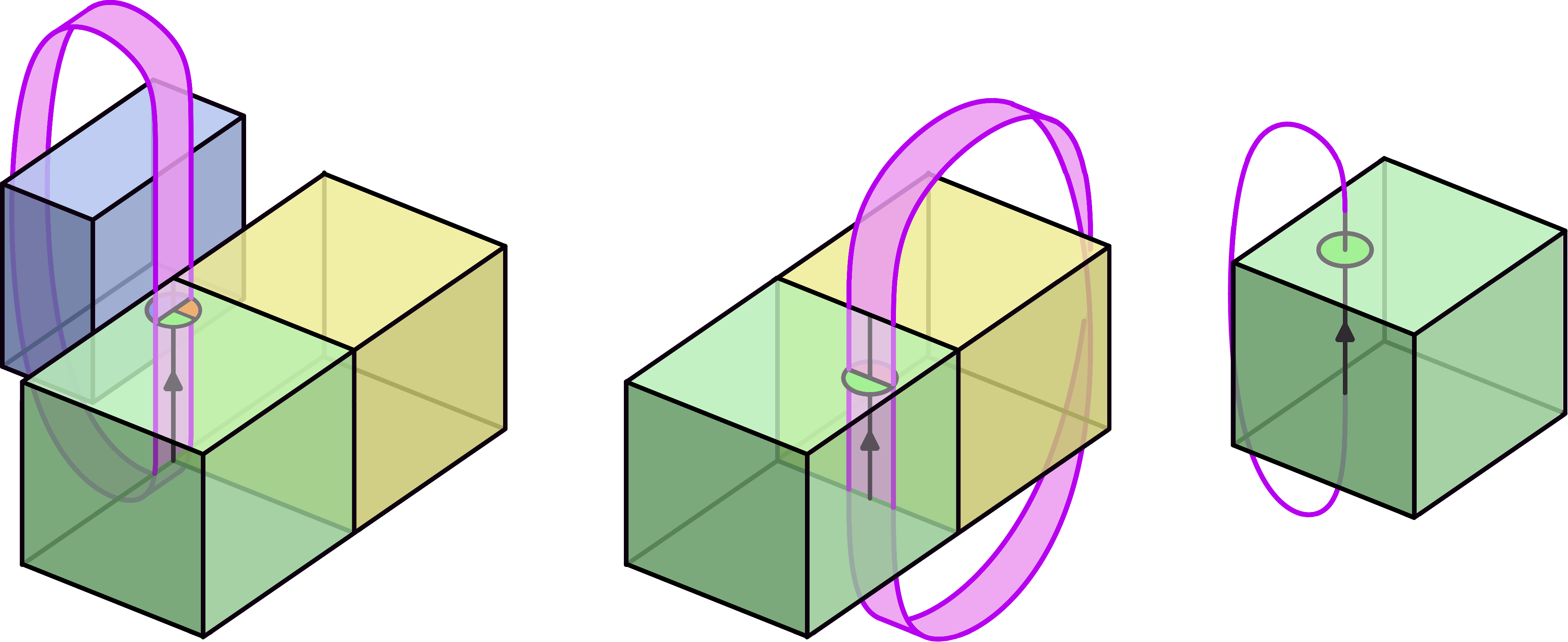}}
        \end{picture}
    \end{center}
    \caption{Germs of $\PP$-quadrant (in dark blue, bright green and orange) on periodic orbits and their orbit by the flow (in purple). Three periodic orbits of $\phi$ are represented, with respectively 3, 2 and 1 germs of $\PP$-quadrants. The orbit in the middle has non-oriented stable/unstable leaves, and the two germs of $\PP$-quadrant are exchanged by the flow.}
    \label{figureMarkovPartition}   
\end{figure}

Given a cyclic $\PP$-word $u$, denote by $\wt\gamma_u$ the periodic orbit of $(\wb u,0)$ for $\phi^\PP$, and by $\gamma_u$ the periodic orbit of $\pi_\PP(\wb u)$ for $\phi$. We call the \emph{realization} of $u$ the map obtained as a restriction of~$\pi_\PP$ to $\wt\gamma_u\to\gamma_u$. We later need two technical lemmas on the realization of periodic words. They should be known by experts but we do not know if a proof is already present in the literature.

\begin{lemma}\label{lemmaRealizationDegree}
    Let $u$ be a cyclic $\PP$-word. When $\pi_\PP(\wb u)$ is on the boundary $\partial^{su}\PP$ and have a non-orientable stable or unstable leaf, the realization of $u$ is a degree two covering map. Otherwise, it is a homeomorphism.
\end{lemma}

\begin{lemma}\label{lemmaDistinctRealization}
    Let $u,v$ be two distinct periodic $\PP$-words with $u_0=v_0=\PR$. Then $\pi_\PP(u)$ and $\pi_\PP(v)$ are on distinct stable and unstable leaves of $W_\PR^{s}$ and $W_\PR^{u}$.
\end{lemma}

\begin{proof}[Proof of Lemmas \ref{lemmaRealizationDegree} and \ref{lemmaDistinctRealization}]
    The realization map is a cover since $h_\PP$ is a semi-conjugacy between two periodic orbits. 

    Let $p$ in $M$ be a periodic point. Given an embedding $f\colon[-1,1]^2\to\FF^s(p)$ which sends $(0,0)$ to $p$ and $\{0\}\times[-1,1]$ inside the orbit of $p$, we call local stable half leaf of $p$ the images by $f$ of the sets $[-1,0]\times[-1,1]$ and $[0,-1]\times[-1,1]$. Two local stable half leaves of $p$ (for two different functions $f$) are said equivalent if they contain a common local stable half leaf of $p$. Define a germ of stable half leaf of $p$ to be an equivalence class of local stable half leaves of $p$. The flow $\phi$ induces a natural flow on the germs of stable half leaves.
    
    A germ of stable half leaf $g$ is said to be inside the stable boundary $\partial^s\PP$ of the Markov partition if one (or any small enough) local stable half leaf in $g$ lies inside $\partial^s\PP$. We define similarly the germs of unstable half leaves of periodic points. 
    
    We claim that, the set of germs of stable half leaves (based on periodic points) which lies inside $\partial^s\PP$ is invariant by the flow $\phi$, in the future and in the past. Indeed by definition of the Markov partition, it is invariant by pushing along the flow in the future. And by periodicity, any germ of stable half leaf obtained by pushing along the flow for a time $-t$ for $t\geq0$ can be obtained by pushing along the for a time $2nT-t\geq 0$ for any large enough $n\in\NN$, and $T$ is the period of the orbit it is based on. One consequence is that when $p$ lies on a stable boundary of Markov rectangle of $\PP$, any small enough local stable half leaf is either included inside $\partial^s\PP$, or it intersects $\partial^s\PP$ only along the orbits of $p$ (this is not necessarily true for non-periodic points).

    We define in a similar manner a local $\PP$-quadrant at $p$ (illustrated in Figure \ref{figureMarkovPartition}). Given a neighborhood $U$ of $x$, the connected components of $U\setminus\partial^{su}\PP$ which contains $x$ in their closures are called $\PP$-quadrants. Two $\PP$-quadrants are said equivalent if their intersection contains a common $\PP$-quadrant. When restricted to small enough $\PP$-quadrants (for $U$ included in a small enough ball around $x$), its induces a equivalence class. Call germ of $\PP$-quadrant an equivalence class of small enough local~$\PP$-quadrant. Since the flow preserves the set of germs of stable and unstable half leaves inside~$\partial^{su}\PP$ (based at periodic points), it induces a flow on the set of germs of $\PP$-quadrant (based at periodic points).

    Given a periodic point $q\in M_\PP$, there exists a unique germ $\OO(q)$ of $\PP$-quadrant based at $h_\PP(q)$ which is included in the same Markov cuboid than $h_\PP(q)$. More precisely, if $q=(v,t)$ for some $0\leq t<T_\PP(v)$, then the germ $\OO(q)$ intersects the interior of the cuboid $\PC_{v_0}$ (and even is included in the cuboid when $t>0$). 
    
    Notice that the map $\OO$ conjugates the flows $\phi$ and $\phi^\PP$ (restricted to periodic points). We claim that $\OO$ is a bijection (on periodic points). Write $q=(v,t)$ for some $0\leq t<T_\PP(v)$, and $p=h_\PP(q)$ the point on which $\OO(q)$ is based. The sequence $v$ is exactly the bi-infinite sequence of rectangles whose interiors are intersected by the orbit of the germ of $\PP$-quadrant $\OO(q)$. The time $t$ is equal to zero if and only if $p$ lies one the Markov rectangle $\partial^-v_0$, and otherwise $t$ is the distance (along the flow) from $p$ to $\partial^-v_0$ (with a positive sign). Therefore the map $\OO(q)\to(v,t)$ described above induces a reciprocal map to $q\to\OO(q)$.
        
    If $T$ denote the period of $p$, the map $\phi_T$ permutes the germ of quadrant at $p$. This permutation is the identity either when $p$ is in the interior of a Markov cuboid (because there is only one $\PP$-quadrant at $p$), or if its stable and unstable leaves are both orientable. In the remaining case, $\phi_T$ induces an involution (with no fixed points) on the set of germs of $\PP$-quadrant at $p$. So each germ of quadrant are $2T$ periodic.

    The injectivity of the map $q\to\OO(q)$ also proves the second lemma. Indeed if $\pi_\PP(u)$ and $\pi_\PP(v)$ are on the same stable leaves, then they are equal (there is at most one periodic point per stable leaf), and by injectivity $u=v$. 
\end{proof}

\subsection{Vertical and horizontal orders}

Here we suppose that the foliations $\FF^{ss}$ and $\FF^{uu}$ are oriented. 
The set of leaves $\Leaf(W_\PR^u)\simeq\Leaf(W_\PC^u)$ are homeomorphic to a segment, and are endowed with a natural orientation coming from the orientation on $\FF^{ss}$. 

Identify the leaf space $\Leaf(W_\PR^s)$ with any unstable leaf in $\Leaf(W_\PR^u)$, and view it as a vertical segment, oriented upward.
Denote by $\leq_v$ the order on $\Leaf(W_\PR^s)$, coming from its orientation. We define the (large and strict) orders $\leq_v$ and $<_v$ on $\PR$ by pushing back the corresponding order on $\Leaf(W_\PR^s)$ using the map $x\in\PR\mapsto W^s_\PR(x)\in\Leaf(W_\PR^s)$. That is $x\leq_v y$ if $W_\PR^s(x)\leq_v W_\PR^s(y)$, and similarly for $<_v$. Note that $<_v$ is not the strict order on $\PR$ coming from $\leq_V$ since two different points on the same stable leaf are comparable for $\leq_v$ but not for $<_v$. We similarly define the orders $\leq_h$ on $\Leaf(W_\PR^u)$ and on $\PR$, were we view $\Leaf(W_\PR^u)$ as a horizontal segment oriented from left the right. The orders $\leq_v$ and $\leq_h$ are called the \emph{horizontal and vertical orders}. Here we picture the stable foliation to be made of horizontal leaves, oriented to the right, and the unstable foliation to be made of vertical leaves oriented upward.

\begin{figure}
    \begin{center}
        \begin{picture}(140,40)(0,0)
          \put(88,5){$W_\PR^s$}
          \put(57,5){$W_\PR^u$}
          \put(20,5){$\PC_0$}
          \put(61,19){$\PC_1$}
          \put(15,35){$\PC_2$}
          \put(91,27){$\partial^+\PC_0$}
          \put(135,21){$\PR_1$}
          \put(83,19){$\PR_2$}
          \put(0,0){\includegraphics[width=140mm]{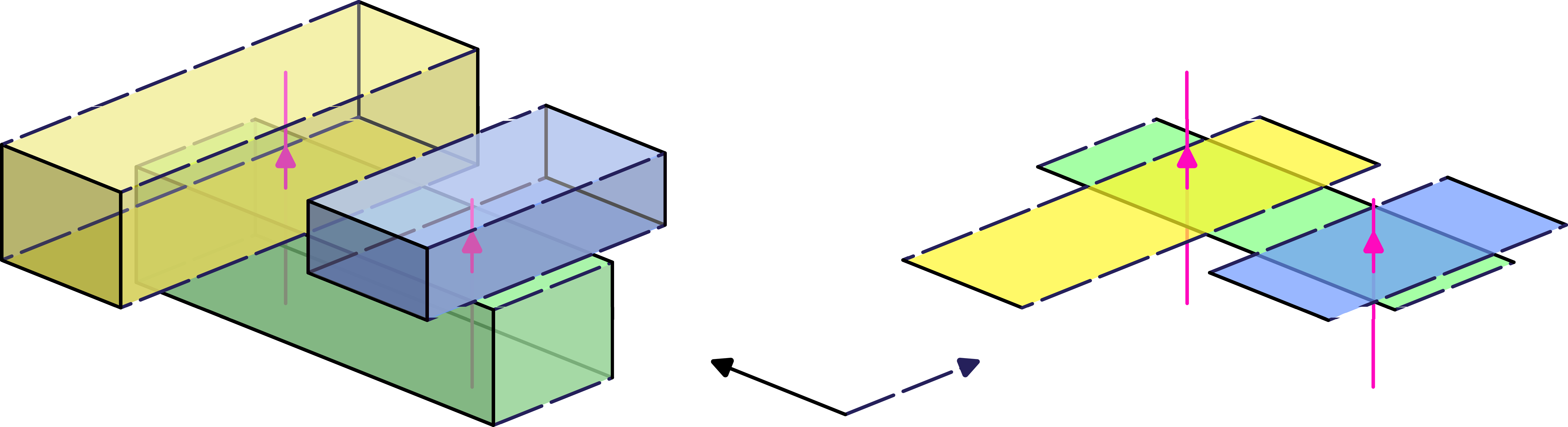}}
        \end{picture}  
    \end{center}
    \caption{Lexicographic orders illustrated on three Markov rectangles $\PR_0,\PR_1,\PR_2$ and on the corresponding cuboids $\PC_0,\PC_1,\PC_2$. The rectangle $\PR_1$ is below $\PR_2$ relatively to $\PR_0$. Whenever $\PR_0\PR_1u_1$ and $\PR_0\PR_2u_2$ are $\PP$-words in $\PW^+$, we have  $\PR_0\PR_1u_1<_v\PR_0\PR_2u_2$.}
    \label{figureOrder}
\end{figure}

\paragraphc{Lexicographic order in $\PW$.} We define in the same manner horizontal and vertical orders on $\PW$. We do not simply lift the orders on $\PR$, because the map $\PW\to M$ is not injective. 
Given some rectangles $\PR_0$, $\PR_1$, $\PR_2$ in $\PP$ and $\PC_0,\PC_1,\PC_2$ the corresponding cuboid, we say that $\PR_1$ is below $\PR_2$ relatively to $\PR$ if $\Int \partial^+\PC_0$ intersects both $\Int\PR_1$ and $\Int\PR_2$, and the intersection $\PC_0\cap\PR_1$ is inferior to $\partial^+\PC_0\cap\PR_2$ for the order $\leq_v$ on $\Leaf(W_{\PC_0}^u)$. The notion is illustrated in Figure~\ref{figureOrder}.

Denote by $\PW^+$ and $\PW^-$ the sets of $\PP$-word of the form  $u\colon\intinto{0}\to\PP$ and respectively $u\colon\intoint{0}\to\PP$. They correspond to the sets of leaves of the foliations $W_l^s$ and $W_l^u$. They come with the distances induced by $d_{\PP}$. 
We define on $\PW^+$ the order $\leq_v$ defined by $u\leq_v v$ if and only if either $u=v$, or there exists $n\geq 0$ satisfying that $u_i=v_i$ for all $0\leq i\leq n$, and $u_{n+1}$ is below $v_{n+1}$ relatively to $v_n=u_n$. In particular, $u$ and $v$ are comparable for $\leq_h$ if and only if $u_0=v_0$. We define similarly an order $\leq_h$ on $\PW^-$. These orders are used in Section \ref{sectionLK} to compute explicitly a linking number.

Similarly to above, we define the orders $\leq_v,<_v,\leq_h,w_h$ on $\PW$ by pushing back the orders $\leq_v,<_v$ by the projection $\PW\to\PW^+$ and the orders $\leq_h,w_h$ by the projection $\PW\to\PW^-$. One can verify that for $u\in\PW$, the leaf $W^s_\PR(\pi_\PP(u))\subset u_0$ depends only on $p^+(u)$.
One can verify that for $u\in\PW$, the leaves $W^s_\PR(\pi_\PP(u))$ and $W^u_\PR(\pi_\PP(u))$ of the rectangle $\PR=u_0$ depend only the images of $u$ in $\PW^+$ and respectively $\PW^-$. Additionally the maps $u\mapsto\pi_\PP(u)$ is not-decreasing for the orders $\leq_v$ and $\leq_h$.

\subsection{Word combinatorics}

We prove a lemma which is essential to build primitive periodic orbits of an Anosov flow. Let $A$ be a finite set and $A^*$ the free monoid generated by $A$. That is elements in $A^*$ are finite words (possibly empty) whose letters are inside $A$. For $u\in A^*$, we denote by $|u|$ the length of $u$. We also denote by $u\leq v$ if $u$ is a prefix of $w$, that is $v\in u\cdot A^*$.

\begin{lemma}\label{lemmaPrimitiveWord}
    Take two elements $x,y\in A^*$ which are not powers of a common element in $A^*$. Then for all $n,m\geq 1$, the cyclic $\PP$-word $x^{2n}y^{2m}$ is not a power of another element.
\end{lemma}

To prove the lemma, we use the following fact on word combinatorial
\cite[Corollary~6.2.5]{Lothaire97}.
 If two words $x,y\in A^*$ satisfy a non-trivial relation, then they are power of a common word. By satisfy a non-trivial relation, one should understand that the unique monoid morphism $f\colon \{a,b\}^*\to A$, satisfying $f(a)=x$ and $f(b)=y$, is not injective. In particular the hypothesis of the lemma implies that $x^k$ and $y^l$ are distinct for all $k,l\geq 1$.

\begin{proof}
    Assume that the lemma is true for $n=m=1$ and prove it for $m,n$ arbitrary. Take $u=x^n$ and $v=y^m$. If $u$ and $v$ are powers on a common a word, then $u^{k}=v^{l}$ for some integers $k,l\geq 1$, then~$x$ and $y$ satisfy a non-trivial relation, so they are powers of a common word. It is not possible by assumption. So one can apply the case $n=m=1$ to $u$ and $v$, which proves the lemma for any other $n,m$.

    We prove now the case $n=m=1$. Suppose that $x^{2}y^{2}=z^k$ for some $z\in A^*$ and $k\geq 2$. Without lost of generality, we can suppose $|x|\geq|y|$. We consider several cases for the lengths of~$x$ and $z$, and prove that they are all impossible. 

    \textbf{Case 1.} Suppose that $|x|=|z|$. Since $z^k$ and $x^2y^2$ are equal and start with $z$ and $x$, one has~$z=x$. It follows that some power of $x$ is equal to some power of $y$, which is impossible.

    \textbf{Case 2.}  Suppose that $|x|>|z|$. Since $x^2\leq z^k$, there exist $u,v$ prefixes of $z$ and integers $1\leq a,b< k$ such that $x=z^au$ and $x^2=z^bv$. Notice that $b\geq 2a\geq a+1$. We claim that~$xz=zx$. Indeed $z\leq x$, so $xz\leq x^2\leq z^k$. Similarly we have $zx=z^{a+1}u$, so $zx\leq z^bu\leq z^k$ when $b=a+1$, and $zx\leq z^{a+2}\leq z^b\leq z^k$ when $b>a+1$. Hence $zx$ and $xz$ have the same length and are prefixes of the same word. So $xz=zx$. It follows that $x=t^c$ and $z=t^d$ for some $t\in A^*$ and~$c,d\geq 1$. Then we have $y^2=t^{kd-2c}$, which contradicts that $x$ and $y$ are not power of a common element.

    \vline

    We suppose now that $|z|>|x|$ Sine $4|x|\geq 2|x|+2|y|=k|z|$, it follows that $k$ is equal to~2 or~3.
    
    \textbf{Case 3.}  Suppose that $|z|>|x|$ and $k=2$. One has $x\leq z\leq x^2\leq z^2$, so there exists~$u,v\in A^*$ such that $x^2=zv$ and $z=xu$. It follows that $x=uv$. 
    We have $z^2=uvuuvu$ and $x^2y^2=uvuvy^2$, so $vy^2=uvu$. Thus $y$ and $u$ have the same length, and they are both suffix of the same word, so $y=u$. Then $uvu=vy^2=vuu$, so $uv=vu$. Hence $u$ and $v$ are power of a common word, which contradicts the assumption on $x,y$.

    \textbf{Case 4.} Suppose that $|z|>|x|$ and $k=3$. 
    Similarly we can write $x=uv$ and $z=uvu$. Notice that we have $2|y|=3|z|-2|x|=4|u|+|v|\leq 2|x|=2|u|+2|v|$, so $|v|\geq 2|u|$.
    On has $x^2=uvuv\leq z^2=uvuuvu$, so $v\leq uv$. Denote by $w$ the suffix of $v$ satisfying $uv=vw=x$. Since le length of $w$ is equal to $|u|\leq |v|/2$, there exists an integer $p\geq 2$ and $t\in A^*$ with $|t|<|w|$ such that $v=tw^p$. We can write
    \begin{align*}
        z^3 &= uvuuvuuvu \\
            &= uvuvwuvwu \\
            &= x^2wvw^2u
    \end{align*}

    Hence $y^2=wvw^2u=wtw^{p+2}u$. Since we have $|w|=|u|$ and $|t|< |w|<(p+2)|w|$, the middle index in $wtw^{p+2}u$ is inside the term $w^{p+2}$ . So there exist two words $a,b\in A^*$ and $q,r\in\NN$ such that $w=ab$, $q+r+1=p+2$ and $y=wtw^qa=bw^ru$. We cannot have $r=0$, otherwise $|wtw^qa|\geq (q+1)|w|\geq 4|w|>2|w|\geq |bu|$. Hence the words $wtw^qa$ $bw^ru$ are equal and start with $ab$ and $ba$. So $ab=ba$, and they are power of a common element $s\in A^*$. It follows the sequence of equalities
    \begin{align*}
        u &= ba\in \{s\}^* \\
        wtw^qa &= bw^ru\in\{s\}^* \text{ so } t\in\{s\}^* \\
        v &=tw^p\in \{s\}^* \\
        uv &= vw \text{ so } u\in\{s\}^*
    \end{align*}
    Hence $x,y$ are powers of $s$, which is not possible. Therefore the last case is also impossible.
\end{proof}

\begin{lemma}\label{lemmaProducWordInRectangle}
    Let $x,y$ be two cyclic $\PP$-words with $x_0=\PR=y_0$, such that  $\wb x\neq\wb y$. Then the point~$\wb{xy}$ lies horizontally and vertically strictly between the points $\wb x$ and $\wb y$.
\end{lemma}

\begin{proof}
    We mimic a classical argument used to find invariant points in the orbit space. We prove the result for the vertical order only. Write $z=xy$ the product cyclic $\PP$-word, and denote by $x^+,y^+,z^+$ the images of $\wb{x},\wb{y},\wb{z}$ in $\PW^+$. We suppose that~$x^+<_v y^+$ and prove $x^+ <_v z^+ <_v y^+$, the other case is similar.

    Denote by~$\SR^+$ the set of $u\in\PW^+$ satisfying $u_0=\PR$.
    Consider the two maps $f_x,f_y\colon\SR^+\to\SR^+$ satisfying that $f_x(w)$ is the element in $\SR^+$ obtained by concatenating $x$ to the left of $w$, and similarly for $f_y$. That is $f_x$ starts with the same first~$|x|$ elements than $x$, and $\sigma^{|x|}\circ f_x=\id$.
    Denote by $d_\PP^+$ the distance on $\SR^+$ given for~$a\neq b$ by $$d_\PP^+(a,b)=2^{-\inf\{i\geq0,a_i\neq b_i\}}$$

    The map $f_x$ is increasing and $2^{-|x|}$-Lipschitz for the distance $d_\PP^+$.
    Hence the map $f_x\circ f_y=f_{xy}$ is $2^{-|xy|}$-Lipschitz on a compact set, so it has a unique fixed point, which is $z^+$. The goal is to prove that the fixed point is between $x^+$ and $y^+$. 
    Given two points $a,c\in\SR^+$, we denote by~$[a,c]_v$ the set of points $b\in\SR^+$ satisfying $a\leq_v b\leq_v c$. For any point $b\in[a,c]_v\subset\SR^+$, we have 
    \begin{equation}\label{eq-distance}
        d_\PP^+(a,b)\leq d_\PP^+(a,c)
    \end{equation}
    Indeed write $n=\inf\{i\geq 0,a_i\neq b_i\}$. By definition of the order $\leq_v$, $b$ have the same $(n-1)$ coefficients than $a$ and $c$. It follows that $d_\PP^+(a,b)\leq 2^{-n}=d_\PP^+(a,c)$. Applying this for $b=y^+$ and for the interval $[x^+,w]_v$ (for any $w>y^+$ in $\SR^+$) yields $d_\PP^+(x^+,y^+)\leq d_\PP^+(x^+,w)$.
    
    It follows that $f_x([x^+,y^+]_v)\subset[x^+,y^+]_v$. Indeed if $w$ belongs to $f_x([x^+,y^+]_v)\setminus[x^+,y^+]_v$, it is neither smaller than $x$ since $f_x$ is increasing, not greater than $y$ as explained above. So it is in~$[x^+,y^+]_v$ since the order $\leq_v$ is a total order.
    
    So $[x^+,y^+]_v$ is invariant by $f_x$. Similarly it is invariant by $f_y$ and by $f_x\circ f_y$. The set $[x^+,y^+]_v$ is compact and invariant by the contracting map $f_x\circ f_y$, so $f_x\circ f_y$ admits a fixed point~$w$ in~$[x^+,y^+]_v$. By uniqueness of the fixed point, we have $w=z^+$ and thus $x^+<_v z^+<_v y^+$.
\end{proof}

\section{Birkhoff Section}\label{sectionBS}

We give a brief overview on partial sections and Birkhoff sections. Here the 3-manifold~$M$ is supposed oriented. We give a variation of the usual definition of partial sections using the orientation on~$M$.

A \emph{partial section} for~$\phi$ is a pair~$(S^*,f)$ where~$S^*$ is a compact~$\Cinfty$ oriented surface with (possibly empty) boundary and~$f$ is a smooth immersion of $S^*$ in $M$, which satisfies that~$f_{|\Int S^*}$ is positively transverse the flow, and~$f(\partial S^*)$ is a finite union of periodic orbits. It implies that each boundary component of~$S^*$ is sent to a periodic orbit~$\gamma$ of~$\phi$, and~$f\colon \gamma^*\to\gamma$ is a covering map. We usually denote the image $f(S^*)$ by~$S$, and identify the immersed surface~$S$ with~$(S^*,f)$.
A \emph{Birkhoff section} is a partial section~$(S^*,f)$ embedded in its interior, such that there exists~$T>0$ for which any orbit arc of~$\phi$ of length~$T$ intersects~$S$. 

\begin{proposition}[Fried-desingularisation, introduced in~\cite{Fried83}]\label{propFriedDesingularization}
    Let $\phi$ be a smooth flow on a 3-dimensional closed manifold $M$. 
    Let~$S_1$ be a Birkhoff section and~$S_2$ be a partial section for $\phi$. Then there exists a Birkhoff section for $\phi$, relatively homologous to~$[S_1]+[S_2]$ in~$H_2(M,\partial S_1\cup\partial S_2,\RR)$.
\end{proposition}

Fried introduced a desingularisation and used it (together with the Fried sections we define below) to prove that any transitive Anosov flow admits a Birkhoff section. Fried description is not precise and works well only in a specific context. We write a precise proof of the proposition in the 
Appendix~\ref{appendixFS}, 
which additionally works for non-Anosov flows. Note that other proofs of similar results have already be given in close setups \cite{Hryniewicz2020,Marty20}.

\paragraphc{Sign of boundary components.}

Fix a partial section~$S=(S^*,f)$. The orientation on the surface~$S^*$ induces an orientation on its boundary, so that the normal direction going outside~$S^*$ plus the positive direction tangent to the boundary, induce the orientation on $S^*$. We orient the orbits of the flow with the orientation given by the flow.
A boundary component~$\gamma^*$ is said positive if~$f\colon \gamma^*\to f(\gamma)$ is orientation preserving. A partial section is called \emph{positive} if all its boundary components are positive.

\begin{theorem}[Asaoka-Bonatti-Marty {\cite[Theorem A]{ABM22}}]\label{theoremABM}
    An Anosov flow is~$\RR$-covered and positively skewed if and only if it admits a positive Birkhoff section.
\end{theorem}

The proof rely on a construction of Birkhoff sections given some data in the orbit space. 
In Section \ref{sectionGibbsMeasures}, we need more than what is prescribed by  Theorem \ref{theoremABM}, namely we need a Birkhoff section with only one boundary component (Theorem \ref{theoremStrongOpenBook}).
To prove the stronger version, we use the construction of Birkhoff sections developed in \cite{ABM22} (see the proof in Section \ref{section-BS} for a sketch of proof of Theorem \ref{theoremABM}). 

In Section \ref{sectionGibbsMeasures}, we need more than what is prescribed by the theorem, namely we need a Birkhoff section with only one boundary component.
Given a partial section~$(S^*,f)$, we view~$\partial S$ as an algebraic multi-orbit. It is given by $$\partial S=\sum_{\gamma^*\subset\partial S^*}\deg(\gamma^*\to f(\gamma^*))f(\gamma^*)$$

Here~$\deg(\gamma^*\to f(\gamma^*))$ is the algebraic degree of the map~$\gamma^*\xrightarrow[]{f} f(\gamma^*)$. Note that $S$ is positive if and only if $\partial S$ has positive coefficients.

Recall that a $\phi$-invariant signed measure $\mu$ is said to be Reeb-like if it is null-homologous and has a positive linking number with all null-homologous $\phi$-invariant probability measure.

\begin{lemma}\label{lemmaBStoReebLike}
    For any positive Birkhoff section~$S$, the invariant measure $\Leb_{\partial S}$ is Reeb-like.
\end{lemma}

\begin{proof}
    According to Lemma \ref{lemmaMeasureCohomology}, 
    $\Leb_{\partial S}$ is null-homologous. 
    Take $T>0$ as in the definition of Birkhoff section, that is every orbit arc of length $T$ intersect $S$. For any algebraic multi-orbit $\Gamma$ with coefficients in $\NN$, and with support disjoint from $\partial S$, we have $$\link_\phi(\Leb_{\partial S},\Leb_\Gamma)=S\algcap\Gamma \geq \frac{1}{T} \length(\Gamma)=\frac{1}{T}\Leb_\Gamma(M)$$ 
    since $\Gamma$ intersects $\Int S$ only positively.
    It follows from the continuity of the linking number and the density 
    Lemma \ref{lemmaFinitBoundedMeasureDense} 
    that $\link_\phi(\Leb_{\partial S},\mu)\geq \frac{1}{T}$ holds true for any probability measure $\mu\in\MM^0_p(\phi)$.
\end{proof}

\subsection{Fried section}\label{section-FriedSection}

For Anosov flows, Fried constructed~\cite{Fried83}
 a type of partial section which we call a \emph{Fried section}. We give a more precise definition to additionally control the boundary orbits and their multiplicities.

Let us discuss Fried's construction before stating the technical lemma. Suppose that the stable and unstable foliations are orientable. Take two periodic points $x_1$ and $x_2$, close to each other, and denote by $\gamma_i$ the periodic orbit of $x_i$. Consider the almost orbit which consists of traveling along $\gamma_1$ $k_1$-times, then jumping from $x_1$ to $x_2$, traveling along $\gamma_2$ $k_2$-times and jumping back to $x_1$. Assuming 
$x_1$ and $x_2$ close enough, the Shadowing Lemma (see \cite[Theorem 5.3.3]{Fisher19}) on Anosov flows produces a periodic orbit $\delta$ which remains close to this almost orbit. 

\begin{figure}
    \begin{center}
        \begin{picture}(116,32)(0,0)
        \put(0,0){\includegraphics[height=33mm]{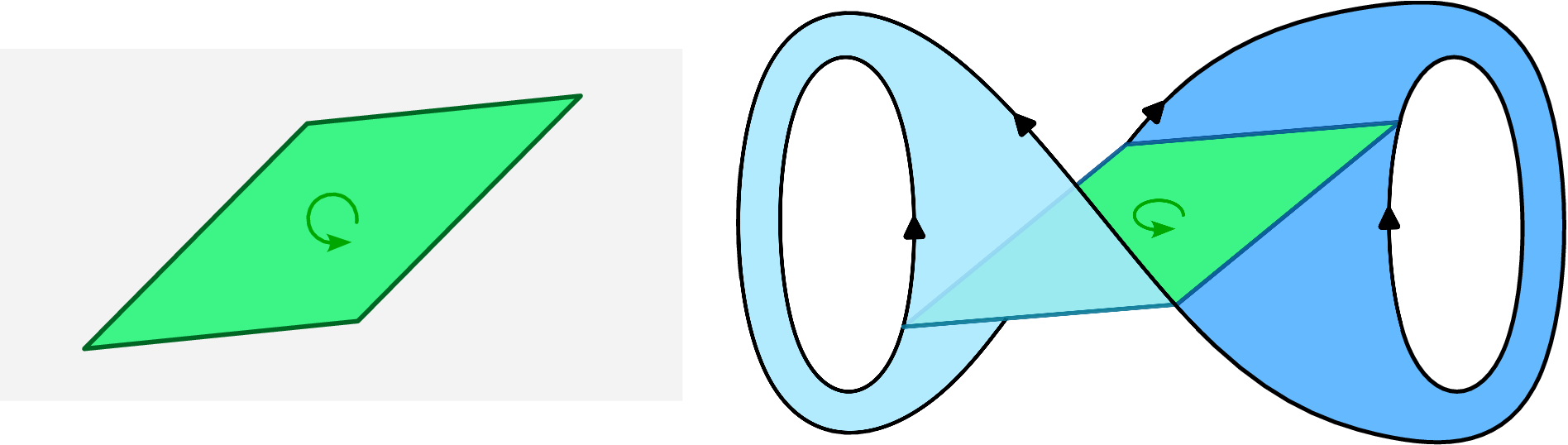}}
        \put(2,4){$\wb{u_1}$}
        \put(41,27){$\wb{u_2}$}
        \put(22,4){$\wb{u_1^{k_1}u_2^{k_2}}$}
        \put(14,26){$\wb{u_2^{k_2}u_1^{k_1}}$}
        \put(64,14){$\gamma_1$}
        \put(104,14){$\gamma_2$}
        \put(89.5,3){$\gamma_w$}
        \end{picture}
    \end{center}
    \caption{Illustration of a 2-chain modeling a Fried Section. In green is represented a transverse polygon (it is~$f(P)\subset\PR$ in the proof of Lemma~\ref{lemmaFriedSection}). In blue are represented two tangent rectangles. Their union can be smoothed to obtain a Fried section.}
    \label{figureFriedSection}
\end{figure}
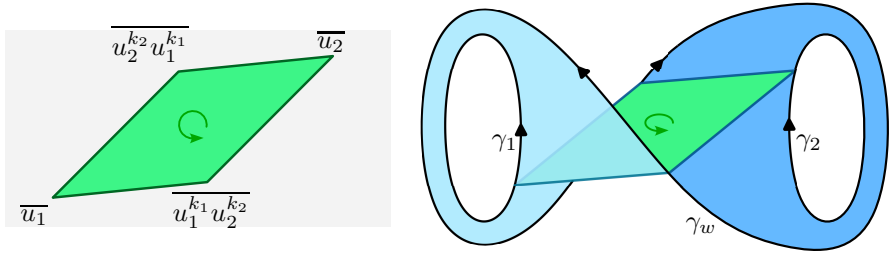

Fried constructed a partial section by drawing a small 4-gon $P$ transverse to the flow, so that $x_1$ and $x_2$ are two opposite corners, the two other corners are on $\delta$, and the two edges adjacent to $x_1$ (resp. $x_2$) are isotopic along the flow. Take one edge adjacent to $x_1$ and push it along the flow on the other edge. The image of the isotopy is a rectangle $R_1$ tangent to the flow. The same procedure on $x_2$ yield a rectangle $R_2$. The union $P\cup R_1\cup R_2$, illustrated in Figure~\ref{figureFriedSection}, is a surface partially tangent and partially transverse to the flow. It can be smoothed in a partial section bounded by $\gamma_1$, $\gamma_2$ and $\delta$, without any consideration for the multiplicities. 

Two problems arise here. If the stable and unstable leaves of $x_i$ are not orientable, we need to travel along $x_i$ an even amount of times to make the construction work. Additionally the multiplicity of the partial section at $\delta$ can be complicated to determine in general. Using a Markov partition, we describe precisely the partial section and determine its multiplicities.

\vvline

Let~$\PP$ be a Markov partition,~$u_1,u_2$ be two distinct, primitive cyclic~$\PP$-words, which start with the same letter~$\PR$. Denote by~$\gamma_i$ the periodic orbit of~$\phi$ corresponding to~$u_i$. We denote by~$N_i\in\{1,2\}$ the degree of the realization of~$u_i$, given in 
Lemma~\ref{lemmaRealizationDegree}. 
Fix two integers~$k_1,k_2\in\NN_{\geq 1}$. We write~$w=u_1^{k_1}u_2^{k_2}$ and denote by~$\gamma_w\subset M$ the periodic orbit corresponding to~$w$. If~$\gamma_i$ has non orientable leaves, we suppose that~$k_iN_i$ is even.
    
\begin{lemma}\label{lemmaFriedSection}
    Assume the stable and unstable foliations oriented. Suppose that~$\wb{u_1}$ is smaller than~$\wb{u_2}$ for both the vertical and horizontal orders. Then there exist~$m\geq 1$ and a partial section~$S$ satisfying~$\partial S =(k_1N_1\gamma_1+k_2N_2\gamma_2)-m\gamma_w$.
    
    Suppose now that~$\wb{u_1}<_v\wb{u_2}$ and~$\wb{u_2}<_h\wb{u_1}$. Then there exist~$m\geq 1$ and a partial section~$S$ satisfying~$\partial S =-(k_1N_1\gamma_1+k_2N_2\gamma_2)+m\gamma_w$.

    If additionally~$k_1$ and~$k_2$ are both even, then~$m=1$ holds in the two previous statements.
\end{lemma}

We call \emph{Fried section} the partial section~$S$ given in the lemma.

\begin{proof}
    We prove the first case, the second case is similar. 
    We use the notation~$\PR_u$,~$\sigma^u$ and~$T^u$ from 
    Section~\ref{sectionMarkovPartition}. 
    That is~$\PR_u$ is the set of elements~$x\in\PR$ having a short itinerary along~$u$,~$\sigma^u$ is the map~$\PR_u\to\PR$ corresponding to pushing along the flow along the itinerary~$u$, and~$T^u$ is the time we push along~$\phi$ to go from~$x$ to~$\sigma^u(x)$.
    
    We write~$v_i=u_i^{k_i}$. Denote by~$P$ the 4-gon given by~$[0,1]^2$, equipped with the anti-clockwise orientation. By assumption, either~$\gamma_i$ has orientable leaves, or~$k_i$ is even, or~$N_i=2$ (and~$\pi_\PP(\wb{u_i})$ lies in the boundary of the Markov partition~$\PP$). So for~$j\in\{1,2\}$,~$\pi_\PP(\wb{v_1v_2})=\sigma^{v_2}(\pi_\PP(\wb{v_2v_1}))$ and~$\pi_\PP(\wb{v_2v_1})=\sigma^{v_1}(\pi_\PP(\wb{v_1v_2}))$ are on the same side relatively to~$\pi_\PP(\wb{u_j})$, for the horizontal and the vertical orders. In particular~$\pi_\PP(\wb{v_1v_2})$ and~$\pi_\PP(\wb{v_2v_1})$ are in the Markov sub-rectangle of~$\PR$ which admits~$\pi_\PP(\wb{u_1})$ and~$\pi_\PP(\wb{u_2})$ as opposite corners. Therefore, we can find an immersion~$f\colon P\to \PR\subset M$ such that:
    \begin{itemize}
        \item~$(0,0)$,~$(1,0)$,~$(1,1)$ and~$(0,1)$ are respectively sent to~$\pi_\PP(\wb{u_1})$,~$\pi_\PP(\wb{v_1v_2})$,~$\pi_\PP(\wb{u_2})$,~$\pi_\PP(\wb{v_2v_1})$,
        \item the restrictions of~$f$ to each edge of~$P$ is transverse to the stable and unstable foliations,
        \item~$\sigma^{v_1}\circ f(t,0)=f(0,t)$ for every~$t\in[0,1]$,
        \item~$\sigma^{v_2}\circ f(t,1)=f(1,t)$ for every~$t\in[0,1]$.
    \end{itemize}

    In the~$3^{rd}$ point,~$f(t,0)$ lies in the smaller horizontal rectangle containing~$f(0,0)=\pi_\PP(\wb{v_1})$ and~$f(1,0)=\pi_\PP(\wb{v_1v_2})$, so~$f(t,0)$ lies in the horizontal sub-rectangle~$\PR_{v_1}$. Therefore~$\sigma^{|v_1|}\circ f(t,0)$ is well-defined, and similarly for the~$4^{th}$ point. 

    Denote by~$Q$ the blow-up of~$P$ along its four corners, that is with each corners of $P$ replaced by a quarter of a circle, and by $\pi_Q\colon Q\to P$ be the blow-down projection. We still denote by~$(t,0)$ and~$(0,t)$ the points in~$\partial Q$ corresponding to lifts of~$(t,0),(0,t)\in\partial P$, for $t\in(0,1)$. Denote by~$(0^+,0)$ and~$(0,0^+)$ the limits in~$\partial Q$ of~$(t,0)$ and~$(0,t)$ when~$t$ goes to~$0$.
    Take~$\tau\colon Q\to\RR$ a smooth map, and consider the map~$g\colon Q\to M$ given by $$g(x)=\phi_{\tau(x)}(f\circ \pi_Q(x))$$
    
    We choose the map~$\tau$ on the edges corresponding to~$[0,1]\times 0$ and~$0\times [0,1]$ so that~$\tau(t,0)=T^{v_1}\circ f(t,0)+\tau(0,t)$. Then we have~$g(t,0)=g(0,t)$. That is the image of~$g$ is an immersed 8-gon, whose edges corresponding to~$[0,1]\times 0$ and~$0\times [0,1]$ coincide. Similarly we choose~$\tau$ on~$[0,1]\times 1$ and~$1\times [0,1]$, so that the corresponding edges in~$g(Q)$ coincide. We additionally choose~$\tau$ so that the surface obtained by gluing the coinciding edges, is smooth. Denote by~$S$ that surface. 
    
    The surface~$S$ is constructed by pushing~$\im f\subset\PR$ along the flow, so it is transverse to the flow in its interior. It is in fact positively transverse. Indeed we have~$\pi_\PP(\wb{u_1})<_v\pi_\PP(\wb{u_2})$, so~$\pi_\PP(\wb{v_1v_2})<_v\pi_\PP(\wb{v_2v_1})$. Similarly we have~$\pi_\PP(\wb{v_1v_2})>_h\pi_\PP(\wb{v_2v_1})$. It follows that the curve~$f(\partial [0,1]^2)$ is going anti-clockwise. So~$f$ preserves the orientation and~$S$ is positively transverse in its interior. 
    
    By construction,~$S$ has some potential boundaries only on the images by $g$ of the blown-up corners of $Q$, which are sent in the orbits going through the corners of~$f([0,1]^2)$, namely~$\gamma_1$,~$\gamma_2$ and~$\gamma_w$. We need to determine the boundary of~$S$. 
    
    Recall the definition of the realization of a primitive cyclic $\PP$-word $u$. We constructed in 
    Section~\ref{sectionMarkovPartitionSuspension} 
    a suspension flow $\phi^\PP$ for the homeomorphism $\sigma\colon\PW\to\PW$, with the time function corresponding to the first-return time from $\PR$ to the boundary $\partial^+\PC_\PR$ of the corresponding cuboid. Let $u\colon\ZZ/n\ZZ\to\PP$ be a primitive cyclic $\PP$-word (it is not $k$ periodic for $k<n$). To $u$ corresponds a periodic orbit $\wt\gamma_u$ of $\phi^\PP$. The realization is the covering map, from the orbit $\wt\gamma_u$, to the $\phi$-periodic orbit of corresponding to $u$. Here we additionally call the length of the realization of $u$, the length of the orbit $\wt\gamma_u$ (counted with multiplicity if $u$ is a multiple of a primitive cyclic $\PP$-word).

    We have~$\tau(t,0)=T^{v_1}\circ f(t,0)+\tau(0,t)$ for~$t\in(0,1)$, so at the limit we have~$\tau(0^+,0)=T^{v_1}\circ f(0,0)+\tau(0,0^+)$. Therefore~$T^{v_1}\circ f(0,0)$ is the length of the boundary component of~$S$ corresponding to the blow-up of the corner $(0,0)$ (viewed as an orbit traveled multiple times). The time~$T^{v_1}\circ f(0,0)$ is equal to~$k_1$ times the length of the realization of~$u_1$. So the boundary component of~$S$ in~$\gamma_1$ has multiplicity~$k_1N_1\neq 0$ (recall that the realization of~$u_i$ has degree~$N_i$). Therefore we can choose~$\tau$ so that~$S$ is additionally immersed on the boundary component corresponding to the edge~$[(0,0^+),(0^+,0)]\subset\partial Q$. Similarly we can choose~$\tau$ so that~$S$ is everywhere immersed, the boundary component in~$\gamma_2$ has multiplicity~$k_2N_2$, and the boundary component in~$\gamma_w$ has negative multiplicity. 
    
    To be more specific, the length of the boundary of~$S$ in~$\gamma_w$, is given by~$-T^{v_1}\circ f(1,0)-T^{v_2}\circ f(0,1)=-T^{v_1v_2}\circ f(1,0)$. Suppose that~$k_1$ and~$k_2$ are both even. Then~$w=u_1^{k_1}u_2^{k_2}$ is not the power of a smaller cyclic~$\PP$-word (see 
    Lemma~\ref{lemmaPrimitiveWord}). 
    Additionally~$\gamma_w$ is not in the boundary of~$\PP$ (since it intersects the interior of the Markov rectangle~$\PR$ in~$f(1,0)$), so the realization of~$\wb w$ is of degree one. By definition of the realization, the length of the realization of $w$ is~$T^{w}\circ f(1,0)$, so the multiplicity of~$S$ in~$\gamma_w$ is~-1.
\end{proof}

\subsection{Birkhoff section with one boundary component}\label{section-BS}

We prove Theorem~\ref{theoremStrongOpenBook}: 
when an Anosov flow is positively skewed, it admits a Birkhoff section with only one boundary component, with multiplicity one and the corresponding periodic orbit lies in orientable stable and unstable leaves. For that, we follow the proof of Theorem 
\ref{theoremABM} given in~\cite{ABM22}. 

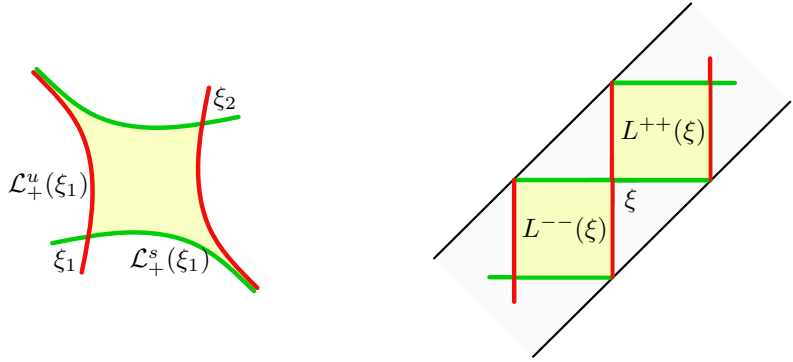
\begin{figure}
    \begin{center}
        \begin{picture}(100,40)(0,0)
        \put(0,0){\includegraphics[width=100mm]{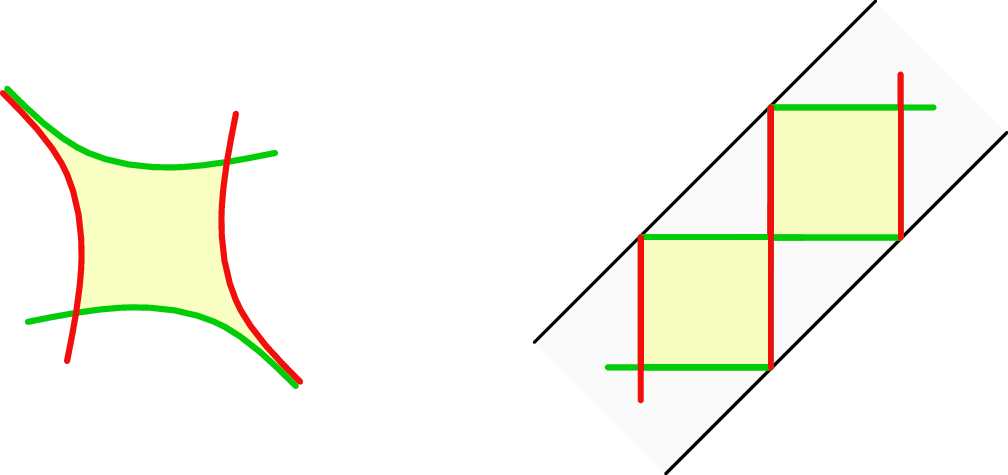}}
        \put(3,12){$\xi_1$}
        \put(24,33.5){$\xi_2$}
        \put(78,20){$\xi$}
        \put(77.5,29){$L^{++}(\xi)$}
        \put(64.5,16){$L^{--}(\xi)$}
        \put(13,12){$\LL^s_+(\xi_1)$}
        \put(-3,22){$\LL^u_+(\xi_1)$}
        \end{picture}
    \end{center}
    \caption{Lozenge with corners $\xi_1$ and $\xi_2$ on the left. The two lozenges with a corner in $\xi$ on the right. On the right, the flow is $\RR$-covered and positively skewed. }
    \label{figureLozenge}
\end{figure}

A key notion to prove the Theorem are lozenges (illustrated in Figure \ref{figureLozenge}). Denote by $(\OS,\LL^s,\LL^u)$ the bi-foliated plane of $\phi$. Fix a pair of orientations on $\LL^s$ and $\LL^u$ so that the orientation on $\LL^s$ plus the orientation on $\LL^u$ gives the orientation on $\OS$. For $x\in\OS$, we denote by $\LL^s_+(x)$ the open half leaves of $x$ (excluding $x$), so that the orientation goes from $x$ to $\LL^s_+(x)$. We define similarly $\LL^s_-(x)=\LL^s(x)\setminus(\LL^s_+(x)\cup\{x\})$, $\LL^u_+(x)$ and $\LL^u_-(x)$.

For $\xi\in\OS$ and $\sigma,\epsilon\in\{-,+\}$, define the quadrant $Q_{\sigma\epsilon}(\xi)$ as the set of points $\zeta\in\OS$ satisfying $\LL^s_\sigma(\xi)\cap\LL^u(\zeta)\neq\emptyset$ and $\LL^u_\epsilon(\xi)\cap\LL^s(\zeta)\neq\emptyset$. 
A (positive) \emph{lozenge} $L\subset\OS$ is a subset of the orbit space which is simultaneously the quadrant $(+,+)$ of a point $\xi_1$ and the quadrant $(-,-)$ of a point $\xi_2$, that is $L=Q_{++}(\xi_1)=Q_{--}(\xi_2)$. In general, one also have to consider a second type of lozenges, defined by $Q_{+-}(\xi_1)=Q_{-+}(\xi_2)$, but only the first type appears in positively skewed orbit spaces. The points $\xi_1$ and $\xi_2$ are respectively called the corners $\xi_{(-,-)}$ and $\xi_{(+,+)}$ of the lozenges $L$. We use lozenges later to construct positive Birkhoff sections.

When the flow is $\RR$-covered and positively skewed, each point $\xi\in\OS$ bounds two positive lozenges. We denote by $L^{++}(\xi)$ the lozenge for which $\xi$ is a corner $\xi_{(-,-)}$, and $L^{--}(\xi)$ the other lozenge. 

\begin{proof}[Proof of Theorem~\ref{theoremStrongOpenBook}]
    We adapted a proof from~\cite{ABM22}. 
    Denote by~$\wh M$ the orientations-bundle covering and~$\wh \phi$ the lifted flow. The lifted flow is also $\RR$-covered and positively twisted, and so transitive (see Theorem \ref{theoremRcoveredTransitive}). Given some $\epsilon>0$, it follows from the Shadowing Lemma (see \cite[Theorem 5.3.3]{Fisher19}) that there exists a periodic orbit~$\wh\gamma$ of~$\wh\phi$ for which any point in~$\wh M$ is a distance at most~$\epsilon$ from~$\wh \gamma$. Lift the orbit $\wh\gamma$ to an orbit $\wt\gamma$ inside the universal covering space $\wt M$, and project it in a point $\xi\in\OS$. Denote by $L$ the lozenge 
    $L^{++}(\xi)$. We take a double covering to ensure that the action $\pi_1(\wh M)$ preserves the orientation of the foliations on $\OS$. Therefore for any $g\in\pi_1(\wh M)$, the lozenge $g\cdot L$ is equal to $L^{++}(g\cdot \xi)$.

    We can chose $\epsilon$ so that every ball of radius $2\epsilon$ is included in a region of $\wh M$ trivially bi-foliated by the weak stable and unstable foliations. It implies that any point in $\OS$ lies in (the interior of) a lozenge of the form $g\cdot L$ for some $g\in\pi_1(\wh M)$. Indeed for $\zeta\in\OS$, take a lift $\wt x$ if $\zeta$ inside $\wt M$. Then there exists a lift $g\cdot\wt\gamma$ of $\gamma$ which goes at distance at most $2\epsilon$ from $\wt x$. The point $\zeta$ and the image $g\cdot \xi$ of $g\cdot\wt\gamma$ in the orbit space lie in a common trivially bi-foliated region of $\OS$. We can additionally that $\zeta$ is above and to the right of $g\cdot \xi$, which implies that $\zeta$ lies in the interior of $g\cdot L$.

    Let~$K$ be a compact subset of the orbit space~$\OS$ that satisfies 
    \begin{equation*}
        \cup_{g\in\pi_1(M)}g\cdot K=\OS    
    \end{equation*}
    For example, take $K$ to be the image in $\OS$ of a big compact ball in~$\wt M$. It follows from the density of periodic points that every point in~$K$ admits an open neighborhood lying inside a periodic lozenge (a lozenge whose corners are periodic). By compactness, there exists a finitely many $g_i\in\pi_1(\wh M)$ which satisfy that 
    \begin{equation*}
        K\subset\cup_{i}g_i\cdot L
    \end{equation*}
    It follows 
    \begin{equation}\label{eq-lozenge-cover-OS}
        \cup_{g\in\pi_1(M)}g\cdot L=\OS
    \end{equation}
    
    Denote by~$\gamma$ the periodic orbits in~$M$ corresponding to the points~$\xi$ (the image in $M$ of the preimage of $\xi$ in $\wt M$).
    Also denote by $\gamma'$ the periodic orbits corresponding to the corner $\xi_{(+,+)}$ of the lozenge $L$. Recall that $g_i\cdot L=L^{++}(g_i\cdot\xi)$, so the corners of the lozenge $g_i\cdot L$ induces in $M$ the same two periodic orbits $\gamma$ and $\gamma'$.
    
    In the proof of the implication~$\Leftarrow$ of~\cite[Theorem A]{ABM22} (end of Section 3.2), we construct a positive Birkhoff section~$S_1$ bounded by~$\gamma$ and by~$\gamma'$. In short, Barbot \cite{Barbot95} constructed an immersed partial section $A$ (it is like a partial section but the interior of the surface is only immersed and not embedded). The surface $A$ is topologically an annulus, and is bounded by two periodic orbits $\gamma$ and $\gamma'$. It satisfies the key property that an orbit in $M$ intersects the interior of $A$ if and only it admits a lift in the orbit space which lies inside $L^{++}(\xi)$. Then $A$ can be desingularised into an embedded partial section, which turns out to be a Birkhoff section because of this property and of the equality \ref{eq-lozenge-cover-OS}.

    \begin{claim}
        If $\gamma$ is $\epsilon$-dense for $\epsilon$ small enough, there exists a Markov partition $\PP$ and a Markov rectangle $\PR$ of $\PP$ which intersects $\gamma$ and $\gamma'$ in two points $p,q$ which satisfy $p<_v q$ and $q<_h p$, that is $p$ is below and on the right of $q$.  
    \end{claim}
     
    Assume the claim and continue the proof. Denote by~$n$ and~$m$ the multiplicity of~$S_1$ along~$\gamma$ and~$\gamma'$. According to Lemma~\ref{lemmaFriedSection}, there exist a partial section~$S_2$ and a periodic orbit~$\delta$ which satisfy~$\partial S_2=\delta - 2n\gamma-2m\gamma'$. According to Proposition~\ref{propFriedDesingularization}, there exists a Birkhoff section~$S_3$ relatively homologous to~$2[S_1]+[S_2]$ in~$H_2(M,\gamma\cup\gamma'\cup\delta,\RR)$. The boundary of~$2[S_1]+[S_2]$ is equal to~$\delta$, with multiplicity one. Hence~$S_3$ has only one boundary component in~$\delta$, with multiplicity one, so it is a positive Birkhoff section.

    Notice that~$\delta$ is homologous to~$2(n\gamma+m\gamma')$. Recall that the map~$sign\colon\pi_1(M)\to\{-1,1\}$ is a group morphism, which measures the orientability of the stable and unstable foliations, on a neighborhood of a simple closed loop. Since~$H_1(M,\ZZ)$ is the abelianization of~$\pi_1(M)$ and~$\{-1,1\}$ is Abelian, it factorizes into a group morphism~$\wt\sign\colon H_1(M,\ZZ)\to\{-1,1\}$. So~$\wt\sign(\delta)=1$. By definition of the function~$\sign$, the stable and unstable leaves of~$\delta$ are orientable.

    \begin{claimproof}
        Ratner \cite{Ratner69}\footnote{see \cite[Theorem 1.3.15]{Iakovoglou2023} for a proof written in english.} proved the following fact for transitive Anosov flows in dimension 3. Given any periodic orbit $\delta$, $\phi$ admits a Markov partition whose stable and unstable boundaries are included in the weak stable and unstable leaves of $\delta$. It follows that for any point~$x$ in $M$, there exists a Markov partition which contains $x$ in the interior of one of its cuboid.

        For all $x\in M$, take an open ball $B_x$ centered at $x$, a Markov partition $\PP_x$ and a rectangle~$\PR_x$ in $ \PP_x$ (and $\PC_x$ the corresponding cuboid) for which the closure of $B_x$ is included in the interior of~$\PC_x$. By compactness, there exists finitely many points $(x_i)_i$ whose corresponding open balls~$B_{x_i}$ cover the manifold. 
        
        Given $i$ fixed, denote by $B^\PR_{x_i}$ the projection of $B_{x_i}\subset\PC_{x_i}$ onto $\PR_{x_i}$ along the backward flow. By assumption, the closure of $B^\PR_{x_i}$ is disjoint from the boundary of $\PR_{x_i}$, so there exists an small ball $b_i\subset\PC_{x_i}$ which satisfies the following. For all point $p'\in b_i$, its projection on $\PR_{x_i}$ is a point~$p$ which satisfies that for all point $q\in B^\PR_{x_i}$, we have $p<_v q$ and $q<_h p$.

        Take $\epsilon$ to be smaller than the radius of all the balls $b_i$. Given a point $q'\in\gamma'$, we can find~$i$ for which $y$ belongs to $B_{x_i}$. Since $\gamma$ is $\epsilon$-dense, it intersects the ball $b_i$ in (at least) a point~$p'$. Then the projections $p,q$ of $p',q'$ onto $\PR_{x_i}$ are in relative positions $p<_v q$ and $q<_h p$.
    \end{claimproof}
\end{proof}

\section{Reeb-like property for Anosov flows}\label{sectionGibbsMeasures}

Gibbs measures are a family of invariant probability measures which satisfy properties similar to invariant smooth volume forms. On the one hand a Gibbs measure corresponds to a volume form in some other differential structure on the manifold. On the other hand Gibbs measure offer some control, which we use to obtain the Reeb-like property.

\subsection{Reeb-like Gibbs measure}\label{sectionRLGm}

We introduce some thermodynamic formalism we need later, more can be found in 
\cite{Bowen08}. 
A Hölder continuous map~$\pot\colon M\to\RR$ is call a \emph{Hölder potential}. Denote by~$\Pot(M)$ the set of Hölder potentials on~$M$, which we equip with the topology given by the sup norm~$\lVert .\rVert_\infty$.
Denote by~$h_\mu(\phi)$ the measured entropy of~$\phi$ relatively to an invariant probability measure~$\mu$. Given a Hölder potential~$\pot\colon M\to\RR$ the pressure of~$\pot$ is the real number $$P(\pot)=\sup_{\mu\in\MM_p(\phi)}\Big(h_\mu(\phi)+\int_M\pot d\mu\Big)$$

Any invariant probability measure~$\mu$ for which the supremum is achieved is call a \emph{Gibbs measure} of~$\pot$ (or the equilibrium state of~$\pot$). It is known that for Anosov flows, Gibbs measures exist and are unique for Hölder potentials 
\cite{Bowen08}. 
Additionally a Gibbs measure has no atom and charges every non-empty open sets.
Given a Hölder potential~$\pot$, we denote by~$\mu_\pot$ its Gibbs measure. 

Recall that an Anosov flow is said homologically full when all homology classes in $H_1(M,\ZZ)$ can be represented by periodic orbits.

\begin{theorem}[Sharp {\cite[Theorem 1]{Sharp93}}]\label{theoremNullCohomGibbsMeasure}
    An Anosov flow is homologically full if and only if it admits a null-homologous Gibbs measure.
\end{theorem}

Similarly we prove

\begin{theorem}\label{theoremReebLikeGibbsMeasure}
    Any~$\RR$-covered and positively skewed Anosov flow admits a Gibbs measure satisfying the Reeb-like condition.
\end{theorem}

In order to prove the theorem, we get some control on a Gibbs measure, and use the continuity of~$\mu_\pot$ in~$\pot$ to obtain implicitly the Gibbs measure in the theorem.

\begin{lemma}\label{lemmaContinuousInPot}
    For an Anosov flow, the maps~$\pot\in\Pot(M)\mapsto P(\pot)\in\RR$ and~$\pot\in\Pot(M)\mapsto \mu_\pot\in\MM_p(\phi)$ are continuous.
\end{lemma}

\begin{proof}
    The pressure is continuous as supremum of 1-Lipschitz maps. Take~$(\pot_n)_n$ a sequence of Hölder potentials converging toward a Hölder potential~$\pot_\infty$. By compactness of~$\MM_p(\phi)$, the sequence~$\mu_{\pot_n}$ accumulates to a probability measure~$\nu\in\MM_p(\phi)$. Then we have:
    
    \begin{align*}
        \bigg|\int_M\pot_n d\mu_{\pot_n} - \int_M\pot_\infty d\nu\bigg|
                &\leq \int_M\big| \pot_n-\pot_\infty\big| d\mu_{\pot_n} + \bigg| \int_M\pot_\infty d\mu_{\pot_n} -\int_M\pot_\infty d\nu  \bigg|\\
                &\leq \lVert \pot_n-\pot_\infty\rVert_\infty + \bigg| \int_M\pot_\infty d\mu_{\pot_n} -\int_M\pot_\infty d\nu  \bigg|\\
                &\xrightarrow[n\to +\infty]{} 0
    \end{align*}

    Since the entropy~$h_\mu(\phi)$ is upper semi-continuous in~$\mu$, the entropy~$h_\nu(\phi)$ is not less than the limit inf of the entropies~$h_{\mu_{\pot_n}}(\phi)$. Hence we have:

    \begin{align*}
        P(\pot_\infty)
                &= \liminf_n P(\pot_n) \\
                &= \liminf_n \Big(h_{\mu_{\pot_n}}(\phi)+\int_M\pot_n d\mu_{\pot_n}\Big) \\
                &\leq h_{\nu}(\phi)+\int_M \pot_\infty d\nu \leq P(\pot_\infty)
    \end{align*}

    Hence~$\nu$ is the Gibbs measure of~$\pot_\infty$, so the map~$\pot\mapsto\mu_\pot$ is continuous.
\end{proof}

\paragraphc{Null-homologous Gibbs measures.}
We first reprove Sharp's Theorem \ref{theoremNullCohomGibbsMeasure}, by constructing a null-homologous Gibbs measure. Our construction allows us to have more control on the linking numbers.
The proof goes as follows. 
Lemma~\ref{lemmaEquiStateControle} 
gives a way to control the Gibbs measure. Given a Hölder potential which mostly charges one periodic orbit, its Gibbs measure charge mostly a neighborhood of that orbit. We apply this lemma on many periodic orbits, to obtain a better control on the Gibbs measure's homology class. 

\vline

Here we suppose that the flow~$\phi$ is homologically full.
Denote by~$n$ the dimension of~$H_1(M,\RR)$. We construct a map~$f:\RR^n\xrightarrow[]{}\Pot(M)\xrightarrow[]{\cohom{\mu_{\cdot}}}H_1(M,\RR)\simeq\RR^n$, such that for~$r$ large enough,~$f$ has no zero on the sphere of radius~$r$, and induces a map~$\partial B(0,r)\to S^{n-1}$ of degree one. 

Let us fix some notations used in the next sections.

\begin{notations}\label{notation-orbit}
    Fix a family of periodic orbits~$\gamma_1^+,\cdots, \gamma_p^+$, such that~$([\gamma_1^+],\cdots,[\gamma_p^+])$ is a basis of~$H_1(M,\RR)$ 
    (see~\cite{Parry86} for the existence). 
    According to Lemma~\ref{lemmaAntiHomotopicOrbit},
    for each index~$i$ there exists a periodic orbit homologous to~$-[\gamma_i^+]$, which we denote by~$\gamma_i^-$.

    Take two families of closed neighborhoods~$U_i^+$ and~$U_i^-$ of~$\gamma_i^+$ and~$\gamma_i^-$, which we take pairwise disjoint. We assume that $U_i^\pm$ is small, in the sense that~$\cap_{t\in\RR}\phi_t(U_i^\pm)=\gamma_i^\pm$.

    Then fix Hölder potentials~$g_i^+$ and~$g_i^-$ such that~$g_i^\epsilon\equiv 1$ on~$\gamma_i^\epsilon$ and~$g_i^\epsilon\equiv 0$ outside~$U_i^\epsilon$, for~$\epsilon\in\{+,-\}$. Define the map~$\pot\colon \RR^n\to \Pot(M)$ given by
    $$\pot(x)=\sum_{\substack{1\leq i\leq n\\ x_i>0}}x_i\cdot g_i^+ + \sum_{\substack{1\leq i\leq n\\ x_i<0}}|x_i|\cdot g_i^-$$
\end{notations}

For each~$x\in\RR^n$, the map~$\pot(x)$ is a Hölder potential, so it admits a Gibbs measure $\mu_{\pot(x)}$.
Take a basis~$(\omega_i)_i$ of~$H^1(M,\RR)$ which is Poincaré dual to the basis~$([\gamma_i^+])_i$. More precisely~$\omega_i\cdot[\gamma_j^+]$ equals 1 when~$i=j$, and it equals 0 otherwise. We define the map 
$f_\pot\colon\RR^n\to\RR^n$ defined by 
$$f_\pot(x)=\big(\cohom{\mu_{\pot(x)}}\cdot \omega_1,\cdots, \cohom{\mu_{\pot(x)}}\cdot \omega_n\big)$$

\begin{proposition}\label{propDegreOneFluxMap}
    There exists~$r>0$ such that~$f_\pot^{-1}(0)\subset [-r,r]^n$, and such that the map $\big(x\in\partial[-r,r]^n\big)\mapsto \Big(\frac{f_\pot(x)}{\lVert f_\pot(x)\rVert_2}\in S^{n-1}\Big)$ is of degree 1.
    Additionally there exists~$x\in(-r,r)^n$ satisfying~$f_\pot(x)=0$.
\end{proposition}

Theorem~\ref{theoremNullCohomGibbsMeasure} is a consequence of the proposition, since $f_\pot(x)=0$ if and only if $\cohom{\mu_{\pot(x)}}=0$.
To prove the proposition, we first need some control on the Gibbs measure. Denote by~$h_{top}$ the topological entropy and recall the variational principle:
\begin{equation}\label{eq-variation-principle}
    h_{top}(\phi)=\sup_{\mu\in\MM_p(\phi)} h_\mu(\phi)
\end{equation}

\begin{lemma}\label{lemmaEquiStateControle}
    Let~$\gamma$ be a periodic orbit of the flow,~$D>C>0$ be two real numbers,~$U\subset M$ be a neighborhood of~$\gamma$, and~$\pot$ be a Hölder potential. Suppose that~$\pot\leq D$ with equality on~$\gamma$, and~$\pot\leq C$ outside~$U$. Then we have $$\mu_\pot(M\setminus U)\leq \frac{h_{top}(\phi)}{D-C}$$ 
\end{lemma}

\begin{proof}
    By definition a Gibbs measure maximizes the pressure, so we have 
    \begin{align*}
        D\mu_\pot(M\setminus U)+D\mu_\pot(U)
            &= D = \int_M\pot d\Leb_\gamma \\
            &\leq \int_M\pot d\Leb_\gamma+h_{\Leb_\gamma}(\phi) \\
            &\leq \int_M\pot d\mu_\pot+h_{\mu_\pot}(\phi)\\
            &\leq C\mu_\pot(M\setminus U)+D\mu_\pot(U)+h_{top}(\phi)
    \end{align*}
    We use the variational principle (Equation \ref{eq-variation-principle}) for the last inequality.
    It follows 
    $$(D-C)\mu_\pot(M\setminus U)\leq h_{top}(\phi)$$
\end{proof}

\begin{lemma}\label{lemmaEquiStateThinControle}
    Let $\gamma_1\hdots\gamma_n$ be periodic orbits and $U_1\hdots U_n$ be neighborhood of these orbits, so that for each $i$, we have $\cap_{t\in\RR}\phi_t(U_i)=\gamma_i$.
    Let $(\pot_n)_n$ be a sequence of potentials and $D_n>C_n>0$ be constants satisfying:
    \begin{itemize}
        \item $\pot_n\leq D_n$ with equality on one of the $\gamma_i$,
        \item $\pot_n\leq C_n$ outside $\cup_i U_i$
        \item $D_n-C_n\xrightarrow[n\to+\infty]{}+\infty$.
    \end{itemize}
    Then all the accumulation points of the sequence $(\mu_{\pot_n})_n$ are probability measures of the form $\sum_i a_i\Leb_{\gamma_i}$, where $a_i\geq 0$.
\end{lemma}

\begin{proof}
    Let $\nu$ be an accumulation point of the sequence $(\nu_{\pot_n})_n$. Lemma \ref{lemmaEquiStateControle} implies that for each $n$, we have 
    $$\mu_{\pot_n}(M\setminus \cup_iU_i)\leq \frac{h_{top}(\phi)}{D_n-C_n}$$
    
    So the mass $\mu_{\pot_n}(M\setminus \cup_iU_i)$ of the open set $M\setminus \cup_iU_i$ converges toward zero. According to the Portmanteau Lemma, we have $\nu(M\setminus \cup_iU_i)=0$. Using the assumption that $U_i$ is small, we obtain $$\nu(M\setminus \cup_i\gamma_i)\leq \cup_{t\in\QQ}\nu(M\setminus\cup_i\phi_t(U_i))=0$$ Therefore $\nu$ is an invariant probability measure supported on the orbits $\gamma_i$.
\end{proof}

\def\sign{sign}
For~$x\in\RR^n$, we denote by ~$x_i$ its~$i^{th}$ coordinates, such that~$x=(x_1,\cdots,x_n)$. Similarly for a map~$f$ with values in~$\RR^n$, we write~$f=(f_1,\cdots,f_n)$.
Denote by $\length(\gamma)$ the length of a periodic orbit $\gamma$ for the flow $\phi$.

\begin{proof}[Proof of Proposition~\ref{propDegreOneFluxMap}]
    We prove that for any $x\in\RR^n$ whose norm is large enough, there exists a linear form $l_x\colon\RR^n\to\RR$ which satisfies that $l_x(x)$ and $l_x(f_\pot(x))$ are both positive. Then for all $t\in[0,1]$, $l_x(tx+(1-t)f_\pot(x))$ is also positive. Therefore outside a large ball in $\RR^n$, $t\mapsto t\id +(1-t)f_\pot$ is a homotopy inside $\RR^n\setminus\{0\}$. The claim follows.

    Take $\epsilon_1\hdots\epsilon_n$ in $\{+,-\}$ and take $x\in\RR^n$ so that for all $i$, either $x_i=0$ or $\sign(x_i)=\epsilon_i$. Consider the accumulation points of the measures of the form $\mu_{\pot(x)}$ for these points $x$. According to Lemma \ref{lemmaEquiStateThinControle}, these measures are of the form $\mu=\sum_ia_i\Leb_{\gamma_i^{\epsilon_i}}$, with $a_i\geq 0$ and $\sum_ia_i\length(\gamma_i^{\epsilon_i})=1$. So the accumulation points of $f_\pot(x)$ for these points $x$ are of the form $(\epsilon_1a_1,\hdots,\epsilon_na_n)$.

    It follows that when $\lVert x\rVert_\infty$ is large enough, $f_\pot(x)$ is close to $(\epsilon_1a_1,\hdots,\epsilon_na_n)$.
    Therefore if one write 
    $$l_x(y)=\sum_i\epsilon_i\length(\gamma_i^{\epsilon_i})y_i$$ 
    the value $l_x(f_\pot(x))$ is close to $\sum_i a_i \length(\gamma_i^{\epsilon_i})=1$, so $l_x(f_\pot(x))>0$ when $\lVert x\rVert_\infty$ is large enough. It is clear that $l_x(x)>0$, so the criterion given above applies. 
\end{proof}

\paragraphc{Gibbs measures with positive linking number.}
We deform the Gibbs measure obtained in 
Proposition~\ref{propDegreOneFluxMap}
to find a Reeb-like Gibbs measure. 
According to Theorem 
\ref{theoremStrongOpenBook}, 
there exists a Birkhoff section bounded by a single orbit~$\delta$. The measure~$\Leb_\delta$ satisfies the Reeb-like property according to 
Lemma \ref{lemmaBStoReebLike}.
We extend the map~$\pot$ given in the previous subsection. 
Take the sets~$U_i^\pm$ (from Notation \ref{notation-orbit}) disjoint from~$\delta$, and take~$U_0$ to be a compact neighborhood of~$\delta$, disjoints from each set~$U_i^\pm$. Also take a map~$g_0\colon M\to[0,1]$ such that~$g_0\equiv 1$ on~$\delta$ and~$g_0\equiv 0$ outside~$U_0$. We additionally suppose that~$g_0<1$ outside~$\delta$ and~$g_i^\pm<1$ outside~$\gamma_i^\pm$.
We define the map~$\qot\colon\RR^n\times\RR\to \Pot(M)$ by
$$\qot(x,t)=\pot(x)+t\cdot g_0(x)$$

\begin{lemma}\label{lemmaDivergingPotSequence}
    There exists a sequence~$(x_n,z_n)_n$ in $\RR^n\times\RR$ such that~$\lVert (x_n,z_n)\rVert_\infty$ goes to~$+\infty$, and such that for all~$n$, we have~$z_n\geq 0$ and~$[\mu_{\qot(x_n,z_n)}]=0$.
\end{lemma}

\begin{lemma}\label{lemmaConvergingEquilibriumMeasure}
    For any sequence~$(x_n,z_n)_n$ given by the previous lemma, the sequence $\mu_{\qot(x_n,z_n)}$ converges toward~$\frac{\Leb_{\delta}}{\length(\delta)}$.
\end{lemma}

\begin{proof}[Proof of Lemma~\ref{lemmaDivergingPotSequence}]
    Suppose that all points~$(x,z)$ with~$x\geq 0$ and with~$[\mu_{\qot(x,z)}]=0$ remain in a bounded region. Take~$r>0$ large, and define the maps~$f_\qot\colon \RR^n\times\RR\to\RR^n$ and~$h\colon \partial[-r,r]^n\times [0,1]\to S^{n-1}$ by
     $$f_\qot(x,t)=(\cohom{\mu_{\qot(x,t)}}\cdot \gamma_i^+)_i$$
     $$h(x,t)=\frac{f_\qot((1-t)x,rt)}{\lVert f_\qot((1-t)x,rt)\rVert_2}$$
     For~$r>0$ large enough,~$h$ is well-defined, continuous and is a homotopy between a degree one map (see Lemma~\ref{propDegreOneFluxMap}) and a constant map. This is not possible, concluding.
\end{proof}

\begin{proof}[Proof of Lemma~\ref{lemmaConvergingEquilibriumMeasure}]
    Consider a sequence~$(x_n,z_n)_n$ as above and $\nu$ an accumulation point of the sequence. According to Lemma \ref{lemmaEquiStateThinControle}, the measure~$\nu$ is of the form~$\nu=a_0\Leb_\delta+\sum_ia_i\Leb_{\gamma_i^+}$ for some non-negative numbers $a_i$. Since the homology class of a measure is continuous in the measure, we have~$\cohom{\nu}=0$. According to Lemma 
   ~\ref{lemmaMeasureCohomology}, 
    we have $\cohom{\nu}=a_0[\delta]+\sum_ia_i[\gamma_i^+]$ in $H_1(M,\RR)$. The orbit~$\delta$ is null-homologous and $\cohom{\nu}=0$, so one has~$\sum_ia_i[\gamma_i^+]=0$. The family~$([\gamma_i^+])_{1\leq i\leq n}$ is a basis of $H_1(M,\RR)$, so the numbers~$a_i$ for $i\geq 1$ are zero. Hence~$\nu$ is the unique invariant probability measure supported on~$\delta$, that is~$\frac{\Leb_{\delta}}{\length(\delta)}$.

    The sequence~$(\mu_{\qot(x_n,z_n)})_n$ in $\MM_p(\phi)$ has a unique accumulating point, which is~$\frac{\Leb_{\delta}}{\length(\delta)}$. Therefore it converges toward that measure.
\end{proof}

\begin{proof}[Proof of Theorem~\ref{theoremReebLikeGibbsMeasure}] 
    According to Lemmas~\ref{lemmaDivergingPotSequence} and~\ref{lemmaConvergingEquilibriumMeasure}, there exists a sequence of null-homologous Gibbs measure~$\mu_{\qot(x_n,z_n)}$ which accumulates toward~$\frac{\Leb_{\delta}}{\length(\delta)}$. 
    According to Lemma \ref{lemmaBStoReebLike}, $\Leb_\delta$ is Reeb-like. 
    
    Recall that the linking map is continuous (which will proven in Section \ref{subsection-LK-continuous}).  It follows from the continuity and from the compactness of the set of null-homologous invariant probability measures
    (Corollary~\ref{corollaryContinuityMinLK}), 
    that for~$n$ large enough, the measure~$\mu_{\qot(x_n,z_n)}$ has positive linking numbers with all null-homologous invariant probability measures. Namely~$\mu_{\qot(x_n,z_n)}$ has the Reeb-like property when~$n$ is large enough.
\end{proof}

\subsection{$\CoH$ differential structure}\label{sectionLocalProduct}

For $0<\alpha<1$, a map is said to be of class $\Coa$ if it is of class $\Class^1$ and its differential is $\alpha$-Hölder continuous. A map is said to be of class $\CoH$ if it is of class $\Coa$ for some $\alpha\in(0,1)$. A $\Class^1$ atlas~$\Atlas$ on a topological manifold $M$ is said to be a $\CoH$ atlas if the transition maps for the atlas are all of class $\CoH$. Alternatively, we say that $\Atlas$ induces a $\CoH$ differential structure on~$M$.

Fix a Hölder potential~$\pot\colon M\to\RR$. We construct a $\CoH$ differential structure on~$M$, and a reparametrization of~$\phi$, which is of class~$\CoH$ and Anosov for that differential structure. Additionally the Gibbs measure of~$\phi$, associated to~$\pot$, is induced by a Hölder volume form for the new differential structure. This has already be done by 
E Cawley~\cite{Cawley91} 
for Anosov diffeomorphisms and done by 
Asaoka~\cite{Asaoka07} 
for codimension one Anosov flows.

We need an additional control on the invariant volume form, so we use the local product of Gibbs measures given by 
Haydn~\cite{Haydn94}, and follow Asaoka's work.

\begin{theorem}[Haydn~{\cite[Theorem 2]{Haydn94}}]\label{theoremLocalProduct}
    Fix a Hölder potential~$\pot$ of~$\phi$. There exists positive measures~$\mu_\pot^s$ and~$\mu_\pot^u$ on the strong stable and unstable foliations (that is a measure on each leaf) satisfying the properties:
    \begin{enumerate}
        \item The measure~$\mu_\pot$ is locally equal to a positive multiple of the product measure~$\mu_\pot^s\otimes \mu_\pot^u\otimes Leb_\phi$,
        \item~$\phi_t^*\mu_\pot^s=e^{-F_t}\mu_\pot^s$ and ~$\phi_t^*\mu_\pot^u=e^{F_t}\mu_\pot^u$, where~$F_t(x)=\int_0^t(f\circ\phi_s(x)-P(f))ds$,
        \item For any points~$x,x'$ such that~$x'\in\FF^{ss}(x)$, denote the holonomy map~$h\colon\FF^{ss}_l(x)\to\FF^{ss}_l(x')$ from~$x$ to~$x'$ along the weak unstable foliation, we have~$\mu_\pot^s=e^{\omega}h^*\mu_\pot^s$, where~$\omega$ is the Hölder function given by~$\omega=\int_0^{+\infty}\big(f\circ\phi_t\circ h-f\circ\phi_t\big)dt$. A similar statement holds for the measure~$\mu_\pot^u$.
        \item The measure~$\mu_\pot^s$, restricted on a strong stable leaf, is finite on compact segments and positive on non-empty open subsets.
    \end{enumerate}
\end{theorem}

\begin{remark}
    The measure~$\mu_\pot^s$ is a measure on each leaf of the strong stable foliation. We can also interpret it as a pseudo-transverse measure to the flow, on each weak stable leaf. That is given a point~$x\in M$,~$I\subset \FF^{ss}(x)$ a strong stable segment and a continuous function~$T\colon I\to\RR$, the measure~$\mu_\pot^s$ induces the measure~$(\phi_T)^*\big(e^{F_T}\mu_\pot^s\big)$ on the curve~$\phi_T(I)$.  
\end{remark}

We equip the orbit space with a family of measures on the stable and unstable foliations.
For~$p\in\wt M$, we denote by~$\xi_p$ the point~$\pi_\OS(p)$, and by~$\theta_{p}^s\colon\FF^{ss}(\wt\pi(p))\to\LL^s(\xi_p)$ the lift/projection map, defined as follows. The maps~$\wt\pi\colon \wt M\to M$ and~$\pi_\OS\colon\wt M\to\OS$ induce two homeomorphisms~$\wt\FF^{ss}(p)\to\FF^{ss}(\wt\pi(p))$ and~$\wt\FF^{ss}(p)\to\LL^{s}(\pi_\OS(p))$. We define~$\theta_{p}^s$ as the composition of the corresponding homeomorphisms~$\FF^{ss}(\wt\pi(p))\xrightarrow[]{\wt\pi^{-1}}\wt\FF^{ss}(p)\to\LL^s(\xi_p)$.

Then define~$\mu_{p}^s=(\theta_{p}^s)^*\mu_\pot^s$, which is a measure on~$\LL^s(\xi_p)$. The measure~$\mu_{p}^s$ induces a transverse measure to the unstable foliation, on an open subset of~$\OS$ obtained as the union of unstable leaves intersecting~$\LL^s(\xi_p)$. Thanks to Theorem~\ref{theoremLocalProduct}, the transverse measure induced by~$\mu_{p}^s$ is multiplied by a Hölder continuous function when~$p$ varies. 

Fix an orientation on~$\LL^s(\xi_p)$. We denote by~$\int_a^b\mu^s_{p}$ the integral of~$\mu^s_{p}$ on the unstable segment~$[a,b]\subset\LL^s(\xi_p)$, with a negative sign if the orientation on~$\LL^s(\xi_p)$ goes from~$b$ to~$a$.

\begin{lemma}[Asaoka {\cite[Lemma 3.5]{Asaoka07}}]\label{lemmaBiHolderChart}
    The map~$\eta_{p}^s\colon \LL^s(Q(p))\to\RR$,~$\xi\mapsto\int_{\xi_p}^\xi\mu_{p}^s$ is a bi-Hölder homeomorphism. 
\end{lemma}

\begin{lemma}[Asaoka {\cite[Lemma 4.3]{Asaoka07}}]\label{lemmaCAHolderTransition}
    Given a holonomy map along the weak unstable foliation~$h_{p,q}\colon I\subset \wt\FF^{ss}_l(p)\to J\subset\wt\FF^{ss}_l(q)$, the map~$\eta_{q}^s\circ h_{p,q}\circ (\eta_{p}^s)^{-1}$ is of class~$\CoH$.
\end{lemma}

The two lemmas also hold for the measures on unstable leaves. Using these lemmas, we define a~$\CoH$ differential structure on~$\OS$ transverse to~$\LL^u$. We denote by~$\Atlas^\perp_\pot$ the atlas on~$\OS$ whose charts are given by the maps~$(\eta_{p}^s,\eta_{p}^u)$.
Thanks to the two lemmas, the transition maps are of class~$\CoH$.
So~$\Atlas^\perp_\pot$ induces a~$\CoH$ differential structure on the surface~$\OS$. From the definition follow the next lemma.

\begin{lemma}\label{lemmaInvariantStructure}
    The action~$\pi_1(M)\acts\OS$ is of class~$\CoH$ for the atlas~$\Atlas_\pot^\perp$.
\end{lemma}

The Gibbs measure $\mu_\pot$ induces a transverse measure $\mu_\pot^\perp$ on $M$. It lifts to a measure on $\OS$ which is invariant by the action~$\pi_1(M)\acts\OS$. We also denote by $\mu_\pot^\perp$ the corresponding measure on~$\OS$.

\begin{lemma}\label{lemmaHolderMeasure}
    The measure~$\mu_\pot^\perp$ on $\OS$ is induced by a Hölder area form on~$\OS$, for the differential structure~$\Atlas^\perp_\pot$.
\end{lemma}

\begin{proof}
    It follows from Theorem~\ref{theoremLocalProduct} that $\mu_\pot^\perp$ is locally a product of $\eta_{p}^s$ and $\eta_{q}^s$, multiplied by a positive function $f$ (corresponding to the functions in points 2 and 3 in the theorem). Since the foliation~$\FF^{ss}$ and~$\FF^{uu}$ are Hölder for the initial differential structure
    (see~\cite[Corollary 9.4.11]{Fisher19}), 
    the function $f$ is Hölder for the initial differential structure on $\OS$. According to Lemma~\ref{lemmaBiHolderChart}, the function $f$ is also Hölder for the atlas $\Atlas^\perp_\pot$, so $\mu_\pot^\perp$ is Hölder regular.
\end{proof}

\paragraphc{From transverse atlas to atlas.} 

Given a compact surface~$S\subset M$ (only continuously embedded), topologically transverse to the flow, the atlas~$\Atlas_\pot^\perp$ induces a~$\CoH$ differential structure on~$S$ (using a lift of $S$ to~$\wt M$ and the projection to~$\OS$). For that reason we view the atlas~$\Atlas_\pot^\perp$ on~$\OS$ as a transverse atlas to~$\phi$.

\begin{proposition}\label{lemmaTransverseStructToStruct}
    There exist a~$\CoH$ atlas~$\Atlas_\pot$ on~$M$ and a continuous reparametrization~$\psi$ of the flow~$\phi$ which satisfy the following:
    \begin{itemize}
        \item~$\psi$ is of class~$\CoH$ for~$\Atlas_\pot$,
        \item for any $\CoH$ compact surface~$S$, topologically transverse to~$\psi$, the~$\CoH$ differential structures on~$S$ induced by~$\Atlas_\pot$ and~$\Atlas_\pot^\perp$ are the same,
        \item if $\id\colon (M,\phi)\to (M,\psi)$ is seen as an orbit equivalence, the measure $\Theta_{\id}(\mu_\pot)$ is induced by a Hölder volume form for~$\Atlas_\pot$.
    \end{itemize}
\end{proposition}

We prove the proposition below.
Given a continuous and positive function~$f\colon M\to\RR^*_+$, denote by~$\psi^f$ the unique flow on~$M$, differentiable in the direction of the flow (for the initial smooth structure) and for which we have~$f\vectD{\psi^f}=\vectD{\phi}$.

We construct an atlas on~$M$ whose charts are given by flow boxes of the flow~$\psi^f$ for some function~$f$. Given a compact surface~$S\subset M$, supposed to be topologically embedded and transverse to~$\phi$, and~$T>0$, if the map~$a\colon S\times[0,T]\to M$ defined by~$a(x,s)=\psi^f_{s}(x)$ is an embedding, the image of~$a$ is called a flow box of $\psi^f$.

Fix~$T>0$ small and take finitely many compact surfaces~$(S_i)_{1\leq i\leq n}$ so that the sets~$\phi_{[0,T]}(S_i)$ are flow boxes. We also suppose that the union of the interior of the flow boxes~$\phi_{[0,T/2]}(S_i)$ gives~$M$. 

Notice that for any continuous function~$f\colon M\to [1,2]$, the set~$\psi^f_{[0,T]}(S_i)$ is a flow box containing~$\phi_{[0,T/2]}(S_i)$, and is included in~$\phi_{[0,T]}(S_i)$. In particular these flow boxes cover $M$. 
We denote by~$a^f_i\colon S_i\times[0,T]\to\psi^f_{[0,T]}(S_i)$ the parametrization of the flow box for the flow~$\psi^f$. That is~$a^f_i(x,s)=\psi^f_{s}(x)$.

The map~$a_i^f$ induces a~$\CoH$ differential structure on~$B_i$, given as the image of the product between~$\Atlas_\pot^\perp$ and the standard smooth structure on~$[0,1]$. 
Take two surfaces~$S_i,S_j$ whose flow boxes intersect. There is a maximal set~$U_{i,j}\subset S_i$ and a continuous function~$T^f_{i,j}\colon U_{i,j}\to\RR^+$, so that for any point~$(x,s)\in S_i\times [0,1]$ the point~$a_i^f(x,s)$ lies in~$S_j$ if and only if~$x$ is in~$U_{i,j}$ and~$s=T^f_{i,j}(x)$. We call the map~$T_{i,j}$ a short return time.

\begin{lemma}\label{lemmaChartTransition}
    Suppose that for all~$i,j\in I$, the map~$T^f_{i,j}$ is of class~$\CoH$ for the atlas~$\Atlas_\pot^\perp$ on~$U_{i,j}\subset S_i$. Then all transition maps of the form~$(a_i^f)^{-1}\circ a_j^f$ are of class~$\CoH$ where they are well-defined. In particular the set of charts~$(a_i)_i$ form a~$\CoH$ atlas on~$M$, compatible with the transverse atlas~$\Atlas_\pot^\perp$.
\end{lemma}

\begin{proof}
    The map~$(a_i^f)^{-1}\circ a_j^f$ coincides with the map~$(x,s)\mapsto (\psi^f_{T_{i,j}(x)}(x), s+T^f_{i,j}(x))$. The map $(x\in U_{i,j}\subset S_i)\mapsto(\psi^f_{T_{i,j}(x)}(x)\in S_j)$ is $\CoH$, by definition of the $\CoH$ differential structure on~$S_i$ and~$S_j$. So the map~$(a_i^f)^{-1}\circ a_j^f$ is as regular as~$T^f_{i,j}$.
\end{proof}

We now prove the proposition. The idea is to define locally $f$ to make the short return times of class $\CoH$. Given one short return time, one can build by hand $f$ to make that short return time regular. To extend $f$ globally, we build it piece by piece. At each step, the short return maps which are well-defined are regular. And we choose the pieces where we define $f$, so that it is afterward possible to extend it by keeping short return times regular.

Given~$i,j$ and a continuous function~$g\colon U_{i,j}\subset S_i\to [1,2]$, consider the map~$f$ defined on~$\{\phi_t(x),x\in U_{i,j}, t\in[0,T^1_{i,j}(x)] \text{ or } [T^1_{i,j}(x),0]\}$ by~$f(\phi_t(x))=T^1_{i,j}(z)g(z)$. For this function $f$, $T^f_{i,j}$ is well-defined on~$U_{i,j}$ and coincides with~$g$, so it is as regular as $g$. To define a function~$f$ which guarantees the regularity of all maps~$T^f_{i,j}$, we define~$f$ inductively using this idea.

\begin{proof}[Proof of Proposition~\ref{lemmaTransverseStructToStruct}]
    For all indexes $i<j$, up to removing a small open neighborhood of~$S_i\cap S_j$ to $S_j$, we can suppose that the surface $S_i$ and $S_j$ are disjoint. It can be done while preserving that the interior of the flow boxes $\phi_{(0,T/2)}(\Int S_k)$ cover $M$. Indeed the parts of the flow boxes we remove, are covered by either the flow boxes of $S_i$ or by the interior of the flow boxes which cover $S_i\cap S_j$. Then all surfaces $S_i$ are disjoint.

    For $i,j\leq n$, we consider pairs of the form~$(x,t)\in S_i\times\RR$ satisfying $t>0$, that $\phi_t(x)$ lies in~$S_j$, and that the orbit arc~$\phi_{(0,t)}(x)$ is disjoint from~$\cup_k\Int S_k$. We say that~$(x,t)$ is an adjacent pair for $(i,j)$. 
    For $m\geq 0$, we define $A_{i,j,m}$ as the set of adjacent pairs~$(x,t)$ for $(i,j)$, which satisfy that the curve $\phi_{(0,t)}(x)$ intersects $\cup_k\partial S_k$ in exactly $m$ points. The surfaces $S_k$ are called \emph{in-between surfaces}.
    We suppose that the surfaces $(S_k)_k$ are in general positions, so that the closures of the sets $A_{i,j,m}$ are compact manifolds of dimension $2-m$. In particular $A_{i,j,m}$ is finite for $m=2$, and empty for $m>2$.
    Notice that every point in $M\setminus\cup_k\Int S_k$ is on a unique orbit arc of the form $\phi_{(0,t)}(x)$, for an adjacent pair $(x,t)$.

    To prove that all short return times are~$\CoH$, it is enough to prove that for every adjacent pair~$(x,t)$ for $(i,j)$, and every in-between surfaces $S_k$ for $(x,t)$, the short return times $T^f_{i,j}$ and~$T^f_{i,k}$ are~$\CoH$ on a neighborhood of~$x$. Indeed all short return times are locally sums of the short return times described above.

    Now we construct~$f$ inductively. For each adjacent pair $(x,t)\in A_{i,j,2}$, choose a small flow box neighborhood $U\subset M$ of $\phi_{[0,t](x)}$. Denote by $S_{k_1},S_{k_2}$ the two in-between surfaces for the pair $(x,t)$. Choose the function $f_U\colon U\to[1,2]$, such that the short return times $T^f_{i,j}$, $T^f_{i,k_1}$ and~$T^f_{i,k_2}$ are well-defined and of class~$\CoH$ close to~$x$.
    Notice that the boundary of the 1-dimensional sub-manifold $A_{i,j,1}$ is $A_{i,j,2}$. Take a small neighborhood $V$ of the set:
    $$U\cup\bigcup_{i,j}\Big(\bigcup_{(x,t)\in A_{i,j,1}} \phi_{[0,t]}(x)\Big)$$ 
    
    Up to taking $U$ and $V$ smaller, we can extend $f_U$ to a function $f_V\colon V\to[1,2]$ satisfying the following. For every adjacent pair $(x,t)\in \cup_{i,j}A_{i,j,1}$ and in-between surfaces $S_k$ for $(x,t)$, the short return times $T^f_{i,j}$, $T^f_{i,k}$ are well-defined and~$\CoH$ close to~$x$. 
    Then we can extend $f$ on all $M$ so that it additionally satisfies that: for every adjacent pair $(x,t)\in \cup_{i,j}A_{i,j,0}$, the short return time~$T^f_{i,j}$ is well-defined and~$\CoH$ close to~$x$. 
    
    Denote by~$\psi=\psi^f$ the reparametrization of~$\phi$ corresponding to $f$. It follows from the above discussion that the flow boxes for $\psi$, and the surfaces $S_k$, induce a $\CoH$ atlas~$\Atlas_\pot$, compatible with the atlas $\Atlas_\pot^\perp$.  The atlas~$\Atlas_\pot$ is constructed using flow boxes for~$\psi$, so~$\psi$ is of class~$\CoH$ for~$\Atlas_\pot$. For a similar reason, $\Theta_{\id}(\mu_\pot)$ is Hölder regular.
\end{proof}

\begin{lemma}\label{lemmaPressurePositivity}
    There exists~$T>0$ for which we have~$\int_0^T(\pot\circ\phi_s(x)-P(\pot))ds<0$ for every~$x\in M$.
\end{lemma}

This statement is already known for symbolic dynamic. We give a proof using only the definition of the pressure.

\begin{proof}
    Take a $\phi$-invariant probability measure~$\mu\in\MM_p(\phi)$. By definition of the pressure, we have:
    $$\int_M(\pot-P(\pot))d\mu\leq -h_\mu(\phi)\leq 0$$
    
    The equality potential holds only when~$\mu$ is the Gibbs measure of~$\phi$ for the Hölder potential~$f$, in this case we have~$h_\mu(\phi)>0$ (see \cite[Theorem 6.1]{Climenhaga19}). Hence~$\int_M(\pot-P(\pot))d\mu<0$ holds in all cases. 
    Suppose that for all~$N\in\NN$, there exists~$x_N\in M$ for which we have~$\int_0^N(\pot\circ\phi_s(x_N)-P(\pot))ds\geq 0$. Denote by~$\nu_n$ the pull back of the Lebesgue measure on~$[0,T]$, by the map~$t\in [0,N]\mapsto \phi_t(x_N)\in M$, rescaled to be a probability measure. By compactness of the set of probability measures on~$M$, the sequence~$\nu_n$ accumulates to some probability measure~$\mu$. By hypothesis, we have~$\int_M(\pot-P(\pot))d\nu_N\geq 0$ for all~$N$, so at the limit~$\int_M(\pot-P(\pot))d\mu\geq 0$. Take a continuous function~$f\colon M\to\RR$ and~$t\in\RR$. For any~$N\in\NN$ we have:
    \begin{align*}
        \bigg|\int_M (f\circ\phi_t-f) d\nu_N\bigg| 
            &= \bigg|\frac{1}{N}\int_0^N(f\circ\phi_{t+s}(x_N)-f\circ\phi_t(x_N))ds\bigg| \\
            &= \bigg|\frac{1}{N}\int_0^t(f\circ\phi_{N+s}(x_N)-f\circ\phi_t(x_N))ds \bigg| \\
            &\leq \frac{2t}{N}\lVert f\rVert_\infty
    \end{align*}
    Hence the integral~$\int_M (f\circ\phi_t-f) d\nu_N$ converges toward~$\int_M(\pot\circ\phi_t-\pot) d\mu=0$. It follows that~$\mu$ is invariant by~$\phi$, which contradicts that~$\int_M(\pot-P(\pot))d\mu$ is negative. Therefore for some~$T>0$, the integral~$\int_0^T(\pot\circ\phi_s(x)-P(\pot))ds$ is negative for all~$x\in M$.
\end{proof}

\begin{lemma}\label{lemmaAnosovNewStruct}
    The flow $\psi$ is Anosov for the atlas~$\Atlas_\pot$.
\end{lemma}

Set~$Y=\vectD{\psi}$. Take~$m$ a continuous metric on~$M$, for the atlas~$\Atlas_\pot$. By construction of the product atlas~$\Atlas_\pot$, the foliation~$\FF^s$ is of class~$\Class^1$ for the atlas~$\Atlas_\pot$. Denote by~$E^{cs}$ the invariant plane field tangent to~$\FF^s$, for~$\Atlas_\pot$. 
We denote by~$|\det_m(d^{cs}\psi_t)|$ the absolute value of the determinant of the restriction of~$d\psi_t$ to~$E^{cs}$, computed with the metric~$m$.

\begin{proof}
    We first prove that~$|\det_m(d^{cs}\psi_t)|$ decreases exponentially in $t\geq 0$.
    Call stable curve any curve which remains in a single leaf of $\FF^s$, and topologically transverse to the flow.
    Denote by~$\mu_\pot^s$ the measure on the leaves of $\FF^{ss}$ coming from the decomposition of the measure $\mu_\pot$ in a local product (see Theorem \ref{theoremLocalProduct}). As discussed previously, $\mu_\pot^s$ induces a measure on every stable curve. 
    
    We define a similar measure using the metric $m$. For a point $p$ and a vector $v\in T_pM$ tangent to $\FF^s$, we denote by $\lVert v\rVert ^\perp_m$ the minimum of $\lVert v+tY(p)\rVert_m$ for $t\in\RR$. For a $\Class^1$ stable curve~$\gamma\colon[0,1]\to M$ and $A\subset [0,1]$ measurable, define the measure $\nu_m^s(A)$ by $$\nu_m^s(A)=\int_A\big\lVert \gamma'(t)\big\rVert ^\perp_mdt$$
    
    It defines a measure $\nu_m^s$ on $\Class^1$ stable curves. 
    Additionally it is of Lebesgue type. By definition of the atlas $\Atlas_\pot$, the measure $\mu_\pot^s$ is also of Lebesgue type on $\Class^1$ stable curves. Therefore there exists a continuous function $g\colon M\to\RR$, which satisfies $\mu_\pot^s=e^g\nu_m^s$. 

    Lemma~\ref{lemmaPressurePositivity} and Theorem \ref{theoremLocalProduct} imply that the measure~$\psi_{-t}^*\mu_\pot^s$ decreases exponentially in~$t>0$, with uniform constants. By compactness, $g$ is bounded. So~$\psi_{-t}^*\nu_m^s$ also decreases exponentially in~$t>0$. Therefore there exists~$A,B>0$ for which~$|\det_m(d^{cs}\psi_t)|\leq Ae^{-Bt}$ for all~$t\geq 0$. 
    
    We use a cone argument to decompose~$E^{cs}$ as a sum of two invariants line bundles. Take a bi-valued continuous section~$U$ of~$E^{cs}$, so that for every~$p\in M$,~$U(p)$ consist of two opposites vectors not parallel to~$Y(p)$. We can change the metric~$m$ so that~$(Y,U)$ gives an orthonormal basis of each plane of~$E^{cs}$. Take~$T>0$ satisfying~$Ae^{-TB}<1$. Denote by~$a\colon M\to\RR$,~$b\colon M\to\RR^+$ two continuous maps satisfying~$d\psi_T(U)=aY+bU$. Then we have:
    $$|\text{det}_m(d^{cs}\psi_T)|=\Bigg|\det\begin{pmatrix}
        1 & a \\
        0 & b 
        \end{pmatrix}\Bigg|=b$$

    It follows that~$b<1$. Take~$\lambda>\frac{\max a}{1-\max b}$ and consider the cone filed~$C\subset E^{cs}$ defined for~$p\in M$ by~$C(p)=\{xU(p)+yY(p), |y|\leq \lambda |x|\}$. The image cone~$\psi_T(C)$ is the set of points  of the form:
    $$x'b(p)U(\psi_T(p))+(y'+x'a(p))Y(\psi_T(p))$$ 
    with~$|y'|\leq\lambda |x'|$.
    We claim that~$\psi_T(C)$ contains $C$. Indeed given~$(x,y)\neq 0$ with~$|y|\leq\lambda|x|$, we have:
    \begin{align*}
        xU(\psi_T(p))+yY(\psi_T(p)) &= x'b(p)U(\psi_T(p))+(y'+x'a(p))Y(\psi_T(p)) \\
            \text{for } x' = \frac{1}{b(p)}x &\text{ and } y' = y-\frac{a(p)}{b(p)}x
    \end{align*}
    Then we have:
    \begin{align*}
        \frac{|y'|}{|x'|} &= \frac{|b(p)y-a(p)x|}{|x|} \\
            &\leq b(p)\frac{|y|}{|x|} + a(p) \\
            &\leq \lambda \max b +\max a <\lambda
    \end{align*}
    Therefore the cone field~$C$ is invariant by~$\psi_{-T}$ for large~$t>0$, and we have a uniform exponential expansion of~$\psi_{-T}$ on~$C$, so by a classical argument of invariant cone fields, there exists a continuous invariant line bundle~$E^s\subset C$, disjoint from the image of~$Y$, which is exponentially contracted by~$\psi_t$ for~$t>0$. Similarly there exists an invariant line bundle~$E^u$ which is exponentially contracted by~$\psi_t$ for~$t<0$, and with~$TM=E^s\oplus E^u\oplus \RR Y$. It follows that~$\psi$ is Anosov.
\end{proof}

\subsection{Plan of proof for Theorems~\ref{theoremSkewedAreReeb} and~\ref{theoremNullCohomVolumeFrom}}\label{sectionPlanProof}

We give a plan for the proof of Theorem~\ref{theoremSkewedAreReeb}: positively skewed $\RR$-covered Anosov flows are orbit equivalent to Reeb-Anosov flows. The plan will serve as a proper proof after the development of the technical tools in the next section. The proof of Theorem~\ref{theoremNullCohomVolumeFrom} follows the same idea.

\begin{proof}[Plan of proof for Theorem~\ref{theoremSkewedAreReeb}]
    According to Theorems~\ref{theoremReebLikeGibbsMeasure},~$\phi$ admits a Hölder potential~$\pot\colon M\to\RR$ whose Gibbs measure~$\mu_\pot$ satisfies the Reeb-like condition (null-homologous and with positive linking numbers with $\MM_p(\psi)$). In section~\ref{sectionLocalProduct}, we construct a~$\CoH$ atlas~$\Atlas_\pot$ on~$M$ and a reparametrization~$\psi$ of~$\phi$, which is of class~$\CoH$ and Anosov for the atlas~$\Atlas_\pot$. Additionally, the Gibbs measure~$\nu_\pot$ of~$\psi$, corresponding to~$\mu_\pot$, is given by a Hölder volume form~$V$ for the atlas~$\Atlas_\pot$. 

    We smooth the flow with a perturbation. Recall that having positive linking numbers with $\MM_p(\psi)$ is an open condition among null-homologous invariant signed measures. But being null-homologous is closed and not open. So we need an extra step to obtain this property after perturbation.

    We follow Asaoka's ideas~\cite{Asaoka07} to smooth the flow~$\psi$. 
    The atlas~$\Atlas_\pot$ can be extended in a class~$\Cinfty$ atlas~$\Atlas_\pot^\infty$, for which~$\psi$ is generated by a~$\CoH$ vector field~$Y$ (Lemma~\ref{lemmaExtendingSmoothStructure}). Then we approximate the Hölder volume form~$V$, corresponding to~$\nu_\pot$, by smooth volume form~$V_n$. The vector field~$Y$ can be approximated in the~$\Class^1$ topology, be a smooth vector field~$Y_n$ which preserves~$V_n$ (see Proposition~\ref{lemmaAsaokaApproximation}, where~$V$ and~$V_n$ appear as~$e^\lambda V$ and~$e^{\lambda_n}V$). Since Anosov flows are~$\Class^1$ structurally stable, the vector field~$Y_n$ generates an Anosov flow~$\psi^n$, orbit equivalent to~$\psi$ and to~$\phi$. Additionally it preserves a smooth volume form~$V_n$.

    Suppose that~$\iota_{Y_n}V_n$ is null-homologous for~$n$ large enough. Since the linking number is continuous 
    (see Corollary~\ref{corollaryContinuityMinLK}), 
    for~$n$ large, the probability measure induced by~$V_n$ has positive linking numbers with all null-homologous~$\psi^n$-invariant probability measure. Hence it satisfies the Reeb-like property. It follows from 
    Theorem~\ref{theoremReebLikeCondition} 
    that~$\psi^n$ is a reparametrization of a Reeb-Anosov flow. Hence~$\phi$ is orbit equivalent to a Reeb-Anosov flow.

    In general,~$\iota_{Y_n}V_n$ is not null-cohomologous, but its cohomology class converges toward the Poincaré dual of the homology class of~$\nu_\pot$ 
    (see Lemmas \ref{lemmaCohomologyContinuity} and \ref{lemmaCohomologyMeasureVolume}), 
    that is zero. We perturbate~$Y_n$ so that~$\iota_{Y_n}V_n$ is null-cohomologous. For that we remark that given a smooth volume form~$W$, we can find smooth vector fields~$Z$ preserving~$W$, so that~$\iota_ZW$ represent any element~$\omega\in H^2(M,\RR)$. Additionally~$Z$ can be taken with a~$\Class^1$ norm bounded by the norm$\lVert \omega\rVert $ multiplied by a function of $W$. The idea is to replace~$Y_n$ by~$Y_n+Z_n$ in the above argument, for a suitable vector field~$Z_n$ satisfying~$[\iota_{Z_n}V_n]=-[\iota_{Y_n}V_n]$, and so $[\iota_{(Y_n+Z_n)}V_n]=0$. 
    
    We would like the vector field~$Y_n+Z_n$ to converge toward~$Y$ in the~$\Class^1$ topology. But we do not control the dependence of the~$\Class^1$ norm of~$Z_n$ on~$W=V_n$. For that reason, we add a step before extending the atlas~$\Atlas_\pot$ is a smooth atlas. We fix a subset~$U\subset M$ which we call generating family of tori ($U$ later corresponds to~$\cup_iU_i$). The set~$U$ is taken small in some sense, so we can choose a atlas~$\Atlas_\pot^\infty$ which satisfies that the measure~$\nu_\pot$ is smooth inside~$U$ (see Lemma~\ref{lemmaSmoothGeneratingTori}). We approximate the volume form~$V$ by smooth volume forms~$V_n$, each of them coincides with $V$ inside~$U$. And as above, we approximate $Y$ by a sequence of vector fields $Y_n$, so that $Y_n$ preserves~$V_n$. 
    
    We choose the set $U$ so that the map~$H^2(U,\RR)\to H^2(M,\RR)$ is surjective. Then we can take a vector field~$Z_n$ as above, but supported inside~$U$ (see Lemma~\ref{lemmaTechFluxCorrection}). In particular $Z_n$ is chosen relatively to a fixed volume form $(V_n)_{|U}=V_{|U}$. Therefore its~$\Class^1$ norm is smaller than a constant (which depends on the constant $V_{|U}$ only) multiplied by the norm of $[\iota_{Y_n} V_n]$, which converges toward zero. Then~$Y_n+Z_n$ converges toward~$Y$ and $\iota_{(Y_n+Z_n)}V_n$ is null-homologous as previously. It follows that for~$n$ large enough,~$Y_n+Z_n$ generates a reparameterized smooth Reeb-Anosov flow which is orbit equivalent to $\phi$.
\end{proof}

In the proof, it is important to start with a measure which has a positive linking number with all $\MM^0_p(\phi)$ (a Reeb-like measure), and to carry that property until we can apply McDuff criterion, which itself requires a positive linking number property. 

The proof of Theorem~\ref{theoremNullCohomVolumeFrom} follows the same plan, without the consideration for the linking number.

\subsection{Smoothing the flow}\label{sectionSmoothing}

Given a class $\Class^{1+H}$ Anosov flow preserving a $\Class^{1+H}$ volume form, we smooth them with an adaptation of Asaoka's technique.

In this section~$\phi$ is supposed homologically full,~$\pot$ is a Hölder potential so that~$\mu_\pot$ is null-homologous,~$\psi$ is a reparametrization of~$\phi$, of class~$\CoH$ and Anosov for the atlas~$\Atlas_\pot$. We also denote by~$Y=\vectD{\psi}$ the vector field generator of~$\psi$, and~$\nu_\pot=cste\cdot\Theta_{\id}(\mu_\pot)$ the $\psi$-invariant probability corresponding to~$\mu_\pot$. The measure $\nu_\pot$ is null-homologous and is induced by a Hölder volume form. When~$\phi$ is additionally~$\RR$-covered and positively skewed, we additionally suppose that~$\nu_\pot$ satisfies the Reeb-like condition. 

Denote by $D^2$ the unit disc in $\RR^2$. A \emph{generating family of tori} is a finite family of disjoint compact subsets~$(U_i)_{1\leq i\leq n}$ of~$M$, such that there exists some homeomorphisms~$f_i\colon D^2\times S^1\to U_i$, satisfying that the homology classes of the curves~$f_i(0\times S^1)$ form a basis of~$H_1(M,\RR)$. When~$M$ is a~$\RR$-homology sphere, generating families of tori are empty families by convention. 
Given a generating family of tori~$(U_i)_i$, a generating sub-family of tori is a generating family of tori~$(V_i)_i$ with~$V_i\subset U_i$ for all~$i$.

The purpure of generating families of tori is to have a subset~$U\subset M$ for which~$H_1(U,\RR)\to H_1(M,\RR)$ is surjective, but which is small enough so that we can extend the~$\CoH$ structure~$\Atlas_\pot$ on~$M$, into a smooth structure for which the flow and the invariant volume form are smooth inside~$U$. 
Denote by~$\Lambda^\infty(T_M)$ the set of~$\Class^\infty$ vector fields on~$M$. Take a norm~$\lVert .\rVert $ on~$H_1(M,\RR)$. 

\def\vect{v}

\begin{lemma}\label{lemmaTechFluxCorrection}
    Let~$M$ be a smooth closed 3-manifold with a generating family of tori~$(U_i)_{1\leq i\leq n}$. Fix a volume form~$V$ on~$M$, supposed to be smooth on each~$U_i$. Then there exist~$A>0$ and a continuous map~$\vect\colon H^2(M,\RR)\to\Lambda^\infty(T_M)$, such that for each~$x\in H^2(M,\RR)$, the following holds:
    \begin{itemize}
        \item $\vect(x)\equiv 0$ outside~$\cup_i U_i$,
        \item $\vect(x)$ preserves $V$, that is $\LL_{\vect(x)}V=0$,
        \item $\lVert \vect(x)\rVert_{\Class^2}\leq A\lVert x\rVert$,
        \item $[\iota_{\vect(x)}V]=x$ in~$H^2(M,\RR)$.
    \end{itemize}
\end{lemma}

\begin{proof}
    We prove a similar result for the manifold~$N=D^2\times S^1$, and push it forward to~$M$.
    We denote by~$(x,y)$ the coordinates on~$D^2\subset\RR^2$, by~$\theta$ the coordinate on~$S^1$, and~$V_0=dx\wedge dy\wedge d\theta$. 
    Take a function~$h\colon D^2\to\RR$ such that~$h\equiv 0$ near~$\partial D^2$ and~$\int_{D^2}h(x,y)dx\wedge dy=1$. We define on~$N$ the vector field~$Y=h(x,y)\partial_\theta$. 
    
    One has~$\LL_{Y}(V_0)=d(\iota_{Y}V_0)=d(h(x,y)dx\wedge dy)=0$, so the vector field~$Y$ preserves the volume form~$V_0$. Let~$S\subset N$ be an oriented compact surface whose boundary is inside~$\partial N$. We claim that the integral~$\int_S\iota_{Y}V_0$ is equal to the algebraic intersection~$S\algcap (0\times S^1)$. To prove that, notice that~$\iota_{X}V_0$ is closed and equal to zero close to~$\partial N$, so~$\int_S\iota_{X}V_0$ depends only on the homology class of~$S$ relatively to~$\partial N$. The equality~$\int_S\iota_{X}V_0=S\algcap (0\times S^1)$ holds when~$S=D^2\times\{*\}$, so it holds for any~$S$ by linearity.

    Up to taking a generating sub-family of tori of~$(U_i)$, we can suppose that there exist some smooth parametrizations~$f_i\colon N\to U_i$. 
    Applying Moser's lemma, we can take~$f_i$ such that~$f_i^*V=\lambda_iV_0$ for some~$\lambda_i>0$. 
    Denote~$Y_i=\frac{1}{\lambda_i}df_i(Y)$, which we extend to a smooth vector field on~$M$ by setting~$Y_i\equiv 0$ outside~$U_i$. We have~$\LL_{Y_i}V=0$. Additionally for any closed surface~$S\subset M$, the integral~$\int_S\iota_{Y_i}V$ is equal to~$S\algcap \gamma_i$, where~$\gamma_i$ is the closed curve corresponding to~$f_i(0\times S^1)$. Therefore~$[\iota_{Y_i}V]\in H^2(M,\RR)$ is Poincaré-dual to~$[\gamma_i]\in H_1(M,\RR)$.

    For any~$x\in H^2(M,\RR)$, we write $x=\sum_ix_i[\gamma_i]$ and define~$\vect(x)=\sum_ix_iY_i$. Then we have~$[\iota_{\vect(x)}V]=x$ and $\lVert \vect(x)\rVert_{\Class^2} \leq C\lVert x\rVert \max_i\lVert Y_i\rVert_{\Class^2}$ for some $C>0$.
\end{proof}

We denote by~$f^*Y$ the vector field~$x\mapsto df(f^{-1}(x))(Y\circ f^{-1}(x))$, conjugated to~$Y$ by a diffeomorphism~$f$. 

\begin{lemma}\label{lemmaSmoothGeneratingTori}
    There exists a generating family of tori~$(U_i)_i$ and an embedding~$f\colon \cup_iU_i\to\RR^3$ of class~$\CoH$ for the atlas~$\Atlas_\pot$, such that
    \begin{itemize}
        \item~$f^*\nu_\pot$ and $f^*Y$ are smooth measure on each~$U_i$,
        \item~$f(V_i)$ is smooth.
    \end{itemize}
\end{lemma}

\begin{proof}
    First construct a generating family of tori by taking a family of closed~$\Class^{1+H}$ curves~$\gamma_i$ whose homology classes form a basis of~$H_1(M,\RR)$. 
    These curves can be taken embedded and disjoint. Then fix a family of disjoint compact neighborhoods~$(U_i)_i$ of the curves~$\gamma_i$.

    Inside~$U_i$, thicken $\gamma_i$ in a~$\CoH$ compact annulus~$A_i$, transverse to the flow~$\psi$. We take a~$\CoH$ embedding~$g_i\colon A_i\to\RR^2$. The measure~$\nu_\pot$ is given by a Hölder volume form~$V$ on~$M$, and~$Y$ is Hölder, so~$i(Y)V$ is a Hölder differential 2-form on~$A_i$. Since~$A_i$ is transverse to~$Y$,~$i(Y)V$ is additionally an area form on $A_i$. Hence~$(g_i^{-1})^*(i(Y)V)$ is a Hölder area form on~$g_i(A_i)$. 
    
    We can suppose that~$A_i$ has been taken such that~$g_i(A_i)$ has a smooth boundary in~$\RR^2$. 
    In~\cite[Theorem 1]{MoserDacorogna90}
    is proven a generalization of Moser's Lemma, which state that there exists a~$\CoH$ diffeomorphism~$h_i\colon g_i(A_i)\to g_i(A_i)$ such that~$(h_i\circ g_i)^*\nu_\pot$ is given by a smooth volume form on~$g_i(V_i)$. Now define the map~$f_i\colon \psi_{[0,\epsilon]}(A_i)\to \RR^3$ by~$f_i(\psi_t(x))=(h\circ g(x), t)\in\RR^2\times \RR$. Since~$\psi$ is of class~$\CoH$,~$f_i$ is also of class~$\CoH$. By construction~$df(Y)=\partial_z$ is constant. Take a closed simple curve~$\delta_i$ in~$\psi_{[0,\epsilon]}(A_i)$ homotopic to~$\gamma_i$ inside~$U_i$ such that~$f(\delta_i)$ is smooth, and take a neighborhood~$V_i$ of~$\delta_i$ inside~$\psi_{[0,\epsilon]}(A_i)$, such that~$f_i(V_i)$ is smooth. The map~$f$ defined as~$f_{|V_i}\equiv f_i+c_i$ for some constants~$c_i\in\RR^3$ satisfies the lemma.    
\end{proof}

\paragraphc{Regularity of the generator.}

The flow $\psi$ is of class $\CoH$. We improve the regularity of $\psi$ in two steps. First we find a smooth structure for which the vector field $Y$ generating $\phi$ is also of class $\CoH$. Then we perturbate $Y$ into a smooth Anosov flow. We need additional considerations when we perturbate $Y$, so that the resulting flow preserves a smooth volume with the Reeb property. We use a generating sub-family of tori to improve the regularity of the Gibbs measure on an open subset. It improves our control when doing the last perturbation.

\begin{lemma}\label{lemmaExtendingSmoothStructure}
    There exists a generating sub-family of tori~$U_i$ and a~$\Class^\infty$ smooth structure~$\Atlas_\pot^\infty$ on~$M$ which extends the~$\CoH$ structure~$\Atlas_\pot$, and such that for the atlas $\Atlas_\pot^\infty$, we have:
    \begin{itemize}
        \item~$Y$ is of class~$\CoH$ for~$\Atlas_\pot^\infty$,
        \item~$\nu_\pot$ is smooth inside each~$U_i$.
    \end{itemize} 
\end{lemma}

\begin{proof}[Proof of Lemma \ref{lemmaExtendingSmoothStructure}]
    Let~$(U_i)_i$ be a generating family of tori. Take a~$\CoH$ map~$f\colon \cup_i U_i\to\RR^3$ given by Lemma~\ref{lemmaSmoothGeneratingTori}, so that~$f^*\nu_\pot$ is smooth. Take~$n\geq 2$ and consider a~$\CoH$ map~$g\colon M\to\RR^3\times\RR^{2n}$ which extends~$f\times 0$ on~$\cup_iU_i$. If~$g$ is chosen generic outside each~$U_i$, then it is an embedding. Up to doing a $\Coa$ isotopy of $g$, we can suppose that the image of $g$ is a smooth sub-manifold of $\RR^{2n+3}$. Then the smooth structure on $\RR^{2n+3}$ pushes back to a smooth structure on $M$ which extends the $\CoH$ atlas $\Atlas_\pot^\infty$. 

    Take $\alpha\in(0,1)$ so that $\psi$ is of class $\Coa$ for the new smooth structure. Hart proved 
    \cite[Theorem A]{Hart83}
    that there exists a $\Class^1$ diffeomorphism $h\colon M\to M$ which conjugates $\psi$ with a $\Class^1$ flow whose generating vector field is also $\Class^1$. Replacing $\Class^k$ by $\Class^{k+\alpha}$ in Hart's proof yields a $\Coa$ diffeomorphism $h\colon M\to M$ which conjugates $Y$ with a $\Coa$ vector field, so that $h$ is the identity on $\cup_iU_i$. An adaptation of Hart's proof to the case $\Class^{k+\alpha}$ is provided 
    in Appendix \ref{appendixHart}.
    Then $\nu_\pot$ is smooth for the new smooth structure since $h$ coincides with the identity inside $\cup_iU_i$.
\end{proof}

\paragraphc{Smoothing the generator.}

For this subsection, we use the notation from Lemma~\ref{lemmaExtendingSmoothStructure} and consider the atlas denoted by~$\Atlas_\pot^\infty$. According to the lemma, there exists a generating family of tori $(U_i)_i$, on which $\nu_\pot$ is smooth. Take~$\alpha\in(0,1)$ such that~$Y$ is of class~$\Coa$ and such that~$\nu_\pot$ is induced by an~$\alpha$-Hölder volume form on~$M$. Let~$V$ be a smooth volume form on~$M$ and $\lambda\colon M\to\RR$ be the function which satisfies $\nu_\pot=e^\lambda V$. It is $\alpha$-Hölder continuous on $M$ and smooth on each tori $U_i$. Take $\lambda_k\colon M\to\RR$ a sequence of $\Cinfty$ function, which approximates $\lambda$ for the $\Ca$ topology, and which coincide with $\lambda$ on each $U_i$. 

\begin{proposition}[Asaoka {\cite[Proposition 2.3]{Asaoka07}}]\label{lemmaAsaokaApproximation}
    There exists a sequence of~$\Cinfty$ vector fields~$Y_k$ on~$M$ such that:
    \begin{itemize}
        \item~$(Y_k)_k$ converges toward~$Y$ in the~$\Coa$ topology,
        \item~$Y_k$ preserves~$e^{\lambda_k}V$ for all~$k$,
    \end{itemize}
\end{proposition}

Asaoka state a weaker version for this proposition. But he proves the version stated above.

Denote by~$\psi^k$ the flow generated by~$Y_k$. By $\Class^1$ structural stability of Anosov flows,~$\psi^k$ is Anosov and orbit equivalent to~$\phi$ for~$k$ large enough
(see \cite[Proposition~3]{Katok91}).
Denote by~$h_n\colon (M,\psi^k)\to (M,\psi)$ an orbit equivalence  so that~$(h_n)_n$ converges toward~$\id$ for the~$\Class^0$ topology. Also denote by~$\mu_k$ the probability measure induced by~$e^{\lambda_k}V$, invariant by~$\psi^k$.

\begin{lemma}\label{lemmaReparametrizationMeasureConvergence}
    The sequence of measures~$\Theta_{h_n}(\mu_n)\in\MM_p(\psi)$ converges toward~$\nu_\pot$. 
\end{lemma}

\begin{proof}
    We prove that the sequence of transverse measures~$h_k^*(\mu_k^\perp)$ converges toward~$\nu_\pot^\perp$, which implies the lemma. Recall that $(\lambda_k)_k$ converges toward $\lambda$, so $(\mu_k)_k$ converges toward $\nu_\pot$.

    Let~$S\subset M$ be a compact surface, transverse to~$Y$ and~$Y_k$ for all~$k$. So there exists~$\epsilon>0$ for which the maps~$\psi^k$ are all injective on the set~$[0,\epsilon]\times S\to M$. Then for any~$k$, we have~$\mu_k^\perp(S)=\frac{1}{\epsilon}\mu_k(\psi^k_{[0,\epsilon]}(S))\leq \frac{1}{\epsilon}$.

    Take a continuous function~$f\colon S\to\RR$ so that~$f\equiv 0$ on a neighborhood of~$\partial S$. We can extend~$f$ to a continuous function~$M\to\RR$, such that~$f\equiv 0$ on a neighborhood of~$\partial S$ inside~$M$. Fix three neighborhoods~$U_1,U_2,U_3\subset S$ of~$\partial S$ for which~$f_{U_i}=0$ and~$\wb U_i\subset\Int U_{i+1}$.
    
    Since the sequence~$(h_k)_k$ converges toward the identity, there exists for~$k$ large enough an embedding~$g_k\colon h_k^{-1}(S\setminus U_1) \to S$, obtained as a small isotopy along the flow~$\psi^k$ from~$h_k^{-1}(S\setminus U_1)$ to~$S$. For $k$ large enough, the image of~$g_k$ contains~$U_2$. Since~$h_k$ converges toward~$\id$ and~$Y_k$ converges toward~$Y$, the map~$g_k\circ h_k^{-1}$ converges toward the inclusion~$S\setminus U_1\hookrightarrow S$. Hence for~$k$ large we have~$S\setminus U_2\subset g_k\circ h_k^{-1}(S\setminus U_1)$. Therefore for large~$k$, the following holds:

    \begin{align*}
        \int_Sf d h_k^*(\mu_k^\perp)
            &=  \int_{S\setminus U_1}f d h_k^*(\mu_k^\perp) \\
            &=  \int_{h_k^{-1}(S\setminus U_1)}f\circ h_k d \mu_k^\perp \\
            &=  \int_{g_k\circ h_k^{-1}(S\setminus U_1)}f\circ h_k\circ g_k^{-1} d \mu_k^\perp \\
            &=  \int_{S\setminus U_2}f\circ h_k\circ g_k^{-1} d \mu_k^\perp \\
        \bigg|\int_S f d\nu_\pot^\perp - \int_Sf d h_k^*(\mu_k^\perp)\bigg|
            &\leq \bigg|\int_S f d\nu_\pot^\perp -  \int_{S\setminus U_2}f\circ h_k\circ g_k^{-1} d \mu_k^\perp\bigg|\\
            &\leq \bigg|\int_{S\setminus U_2}f d \big(\nu_\pot^\perp - \mu_k^\perp\big)\bigg| + \int_{S\setminus U_2}\big|f-f\circ h_k\circ g_k^{-1}\big| d \mu_k^\perp \\
            &\leq  \bigg|\int_{S\setminus U_2}f d \big(\nu_\pot^\perp - \mu_k^\perp\big)\bigg| + \mu_k^\perp(S) \lVert f-f\circ h_k\circ g_k^{-1}\rVert_{S\setminus U_2} \\
            &\xrightarrow[k\to+\infty]{}0
    \end{align*}
    The first term converges toward zero since $\mu_k$ converges toward $\nu_\pot$. The second term converges toward zero since we have $\mu_k^\perp(S)\leq\frac{1}{\epsilon}$, and $(f\circ h_k\circ g_k^{-1})_k$ converges toward $f$ on $S\setminus U_2$. Therefore the measure~$h_k^*(\mu_k^\perp)$ converges toward~$\nu_\pot^\perp$. 
\end{proof}

\begin{lemma}\label{lemmaFluxCorrection}
    There exists a sequence~$Z_k$ of~$\Cinfty$ vector fields on~$M$ such that:
    \begin{itemize}
        \item~$Z_k$ converges toward zero in the~$\Coa$ topology,
        \item~$Z_k$ preserves~$e^{\lambda_k}V$ for all~$k$,
        \item~$\iota_{Y_k}e^{\lambda_k}V$ and~$\iota_{Z_k}e^{\lambda_k}V$ are cohomologous.
    \end{itemize}
\end{lemma}

\begin{proof}
    Let~$\vect\colon H^2(M,\RR)\to\Lambda^\infty(TM)$ be the map constructed in Lemma~\ref{lemmaTechFluxCorrection} for the volume form~$e^\lambda V$, which is smooth inside each~$U_i$. Notice that for~$x\in H^2(M,\RR)$,~$\vect(x)$ is a vector field supported inside the union of the tori~$U_i$, where~$\lambda$ and~$\lambda_k$ coincide. Hence the vector field~$\vect(x)$ preserves the volume forms~$e^{\lambda_k}V$. Now take~$Z_k=\vect([\iota_{Y_k}e^{\lambda_k}V])$. By hypothesis, we have~$[\iota_{Z_k}e^{\lambda_k}V] = [\iota_{Y_k}e^{\lambda_k}V]$, which converges toward zero according to 
    Lemmas~\ref{lemmaReparametrizationMeasureConvergence} and~\ref{lemmaCohomologyMeasureVolume}. 
    It follows that~$Z_k$ converges toward zero for the~$\Class^2$ topology, and so also for the~$\Coa$ topology.
\end{proof}

\begin{corollary}\label{lemmaSmoothFlowFluxProperties}
    For~$k$ large enough, the smooth vector field~$Y_k-Z_k$ generates an Anosov flow, orbit equivalent to~$\phi$, which preserves the smooth volume form~$e^{\lambda_k}V$. Additionally the probability measure induced by~$e^{\lambda_k}V$ satisfies the Reeb-like property for~$k$ large enough.
\end{corollary}

    The corollary and Theorem~\ref{theoremReebLikeCondition} give a proof to Theorem~\ref{theoremSkewedAreReeb}.

\begin{proof}
    Denote by~$\psi^n$ the flow generated by the vector field~$Y_n-Z_n$.
    The vector fields~$(Y_k-Z_k)$ converges toward~$Y$ for the~$\Class^1$ topology. So by structural stability of Anosov flows
    \cite[Proposition~3]{Katok91}, 
    for~$k$ large enough,~$\psi^n$ is Anosov and is orbit equivalent to~$\phi$. Denote by~$\eta_n$ the $\psi^n$-invariant probability measure induced by~$\iota_{(Y_k-Z_k)}e^{\lambda_k}V$. It follows from the previous lemma that~$\eta_n$ is null-homologous.
    Take~$h_n$ a sequence of orbit equivalences between~$\psi^n$ and~$\phi$, so that~$h_n$ converge toward the identity map. According to 
    Lemma~\ref{lemmaReparametrizationMeasureConvergence}, 
    the measure~$\Theta_{h_n}(\eta_n)$ converges toward~$\nu_\pot$. So for~$n$ large enough, the minimum of~$\link_\phi(\Theta_{h_n}(\eta_n),\nu)$ with all~$\nu\in\MM_p^0(\phi)$, is positive
    (see Corollary~\ref{corollaryContinuityMinLK}). 
    It follows from
    Lemma~\ref{lemmaInvarianceLKbyTopEquiv} 
    that the measure~$\eta_n$ satisfies the Reeb-like property.
\end{proof}

\section{On the linking number for invariant measures} \label{sectionLK}

In this section, we prove the theorems and lemmas stated in 
Section~\ref{sectionLKintro}. 
Fix one transitive Anosov flow~$\phi$ on a closed oriented 3-dimensional manifold~$M$. Except in the proof of Theorem~\ref{theoremMeasureLinkingNumber},
we also suppose that the stable and unstable foliations of~$\phi$ are orientable. We define the linking number between two null-homologous invariant signed measures starting for two signed measures $\nu,\mu$ of the type $\Leb_\gamma$, then allowing $\nu$ to be any null-homologous invariant signed measure, and then for all $\nu,\mu$. To distinguish the several definitions, we use the notations $\link_\phi^{ff}$, $\link_\phi^f$ and respectively $\link_\phi$ ($f$ for finitely supported). We prove that they coincide on their common domains.

Let~$\Gamma_1$ and~$\Gamma_2$ be two algebraic multi-orbits, supposed to be null-homologous and with disjoint supports. We define the following linking number 
$$\link_\phi^{ff}(\Leb_{\Gamma_1},\Leb_{\Gamma_2})=S_1\algcap \Gamma_2$$
where $S_1$ is any 2-chain bounded by $\Gamma_1$. It does not depend on the choice of $S_1$. 

The proof of Theorem \ref{theoremMeasureLinkingNumber}
follows two steps. First consider a null-homologous algebraic multi-orbit~$\Gamma$ and a signed measure~$\mu\in\MM^0_s(\phi)$. We define the linking number~$\link_\phi^f(\Leb_\Gamma,\mu)$ roughly as the integral of the transverse measure~$\mu^\perp$ on a 2-chain~$S$ bounding~$\Gamma$. The linking number~$\link_\phi^f(\Leb_\Gamma,\mu)$ is well-defined for all~$\mu$, even when~$\mu$ charges the support of~$\Leb_\Gamma$. 
The linking number maps $\link_\phi^{ff}$ and $\link_\phi^f$ coincide where they are simultaneously defined.

In a second part, we use Markov partitions to give a precise description of the linking number. This part requires much more efforts and technical arguments. We summarize the strategy in Section \ref{section-LK-MarkovPartition}.

\subsection{Linking number with $\Leb_\Gamma$}

\def\MMO{\MM^0_s(\phi^\Delta)}
\def\PM{{\mathbb{P}_\phi M}}

Fix a null-homologous algebraic multi-orbit~$\Gamma$. We define the linking number between~$\Leb_\Gamma$ and any~$\mu\in\MM^0_s(\phi)$. 

Take a set~$\Delta\subset M$, obtained as a union of periodic orbits, containing the support of~$\Gamma$. We consider the \emph{blow-up}~$M_\Delta$ of~$M$ along each orbit of~$\Delta$, constructed below.

Denote by~$X=\vectD{\phi}$ the vector field generating~$\phi$. Given a periodic orbit~$\gamma$ in~$\Delta$, we denote by~$N_\gamma$ the normal vector bundle to~$\gamma$. It is obtained as the set of equivalence classes of vectors~$v$ based on a point~$p$ in~$\gamma$, non-parallel to~$X(p)$, up to the equivalence relation~$v_1\simeq v_2$ if there exist~$a>0, b\in\RR$ such that~$v_1=av_2+bX(p)$. The set~$M_\Delta$ is the union of~$M\setminus\Delta$ and of the circle bundle~$\cup_\gamma N_\gamma$, equipped with the atlas induced by~$M\setminus\Delta$ and by the angular sectors around~$\Delta$. 
We define the map~$\pi_\Delta\colon M_\Delta\to M$ by~$\pi_\Delta(x)=x$ for~$x\in M\setminus\Delta$, and~$\pi_\Delta(v)=x$ if~$v$ is a vector based on~$x\in\Delta$.

For a periodic orbit~$\gamma\subset\Delta$, we denote by~$\TT_\gamma=\pi_\Delta^{-1}(\gamma)\subset M_\Delta$ the boundary component corresponding to~$\gamma$. Since~$M$ is supposed orientable,~$\TT_\gamma$ is a torus. If the flow~$\phi$ is of class~$\Class^k$, $k\geq 1$, it lifts to a~$\Class^{k-1}$  flow~$\phi^\Delta$ on~$M_\Delta$, by setting~$\phi_t^\Delta(v)=d\phi_t(v)$ for a vector~$v$ based on a point in~$\Delta$. We denote by~$X_\Delta=\vectD{\phi^\Delta}$ the vector fields generating~$\phi^\Delta$.

Notice that on each torus~$\TT_\gamma$, the restriction flow~$\phi^\Delta_{|\TT_\gamma}$ has one or two periodic orbits corresponding to the stable direction, and one or two corresponding to the unstable direction. We refer to them as the stable and unstable orbits of~$\TT_\gamma$. Every other orbit converges in the past toward one stable orbit, and in the future toward one unstable orbit. From this remark follows the following lemma.

\begin{lemma}\label{lemmaMeasureSupportOnTorus}\label{lemmaMeasureSupportedOnSUdirection}
    For any signed measure~$\mu\in\MM_s(\phi^\Delta)$, the signed measure~$\mu_{|\TT_\delta}$ is supported on the stable and unstable orbits of~$\TT_\gamma$.

    For any signed measure~$\nu\in\MM_s(\phi)$, there exists a signed measure~$\mu\in\MM_s(\phi^\Delta)$ so that~$\pi_\PP^*\mu=\nu$ and whose total variation is the same as the total variation of~$\nu$. 
\end{lemma}

Below, all homology and cohomology modules are taken with real coefficients, except when specified. We write~$\MMO=(\pi_\Delta^*)^{-1}\big(\MM^0_s(\phi)\big)$, that is the set of signed measures~$\mu\in\MM_s(\phi^\Delta)$ satisfying $\cohom{\pi_\Delta^*\mu}=0$. Denote by
$$H^1(M_\Delta)\xrightarrow[]{P_\Delta} H_2(M_\Delta,\partial M_\Delta)\xrightarrow[]{\partial} H_1(\partial M_\Delta)$$
the Lefschetz-Poincaré duality and the boundary map.

\begin{lemma}
    For any signed measure~$\mu\in\MMO$ and~$\omega\in H^1(M_\Delta)$, the quantity~$\cohomG{\mu}\cdot\omega$ depends only on~$\mu$ and on~$(\pi_\Delta)_*(\partial P_\Delta(\omega))\in H_1(\Delta)$.
\end{lemma}

\begin{proof}
    Consider the following diagram

    \begin{center}
        \tikz[
        overlay]{\draw[draw=gray] (2.7,-0.9) rectangle (9.9,-0.3);}
        \begin{tikzcd}
            H^1(M_\Delta) \arrow[rr, "P_\Delta"] & {} & {H_2(M_\Delta,\partial M_\Delta)} \arrow[d, "(\pi_\Delta)_*"] \arrow[r, "\partial"] & H_1(\partial M_\Delta) \arrow[d, "(\pi_\Delta)_*"] \\
            H^1(M) \arrow[r, "P"] \arrow[u,, "\pi_\Delta^*"] & H_2(M) \arrow[r, "j^*"] & {H_2(M,\Delta)} \arrow[r, "\partial"] & H_1(\Delta)   
        \end{tikzcd}
    \end{center}

    Where $P$ is the Poincaré duality map and the gray rectangle is a sub-sequence of the long exact sequence in homology, coming from the inclusion $\Delta\subset M$. It is clear that the right part is commutative.
    The cohomology set $H^1(M,\ZZ)$ is naturally isomorphic to the set of homotopy class of maps of the form $f\colon M\to\RR/\ZZ$. For a smooth map $f\colon M\to\RR/\ZZ$, $P$ can be described by~$P([f^*dt])=[f^{-1}(x)]$ for any regular value~$x$ of~$f$. A similar description holds for $P_\Delta$. Using these characterizations and the fact that $H^1(M,\ZZ)$ spans $H^1(M,\RR)$, one can verify that the left part of the diagram also commutes.

    Let~$\omega\in H^1(M_\Delta)$ be a 1-form which satisfies~$(\pi_\Delta)_*(\partial P_\Delta(\omega))=0$. By exactness, there exists $S\in H_2(M)$ satisfying $j^*S=(\pi_\Delta)_*(P_\Delta(\omega))$. By commutativity, we have $\omega=\pi_\Delta^*\circ P^{-1}(S)$. It follows that $$\cohomG{\mu}\cdot\omega=\cohom{\pi_\Delta^*\mu}\cdot P^{-1}(S)=0$$ 
\end{proof}

Given~$\mu\in\MM_s(\phi^\Delta)$ satisfying~$\cohom{\pi_\Delta^*\mu}=0$, we define~$\link_\Gamma^\Delta(\mu)=\cohomG{\mu}\cdot[\alpha]$ for any smooth differential 1-form~$\alpha$ with~$(\pi_\Delta)_*(\partial P_\Delta(\omega))=[\Gamma]\in H_1(\Gamma,\RR)$. Notice that~$(\mu\in\MMO)\mapsto\link_\Gamma^\Delta(\mu)$ is continuous.

\begin{lemma}\label{lemmaTransverseMeasureOnBoundary}
    Let~$\mu\in\MM_s(\phi^\Delta)$ be an invariant signed measure,~$\gamma$ be a periodic orbit in~$\Delta$ and~$\alpha$ be a closed 1-form on~$\TT_\gamma$.
    Then the integral~$\int_{\TT_\gamma}\iota_{X_\Delta}\alpha d\mu$ depends only on~$\mu(\TT_\gamma)$ and on~$[\alpha]\in H^1(\TT_\gamma,\RR)$. 
\end{lemma}

\begin{proof}
    Thanks to
    Lemma~\ref{lemmaInvarianceImpliesClosed}, 
    the integral~$\int_{\TT_\gamma}\iota_{X_\Delta}\alpha d\mu$ does not depend on the choice of the closed 1-form~$\alpha$ in its cohomology class.
    Denote by~$dt$ the closed 1-form on~$\gamma$ satisfying~$\iota_Xdt=1$. We denote by~$\beta^\perp=\pi_\Delta^*dt$ and by~$\beta^\parallel$ a closed 1-form which is not null-cohomologous and such that the stable and unstable periodic orbits on~$\TT_\gamma$ are tangent to~$\ker(\beta^\parallel)$. Then~$\iota_{X_\Delta}\beta^\parallel$ is zero on the support of~$\mu_{|\TT_\gamma}$ (see Lemma~\ref{lemmaMeasureSupportedOnSUdirection}), so one has 
    \begin{align*}
        \int_{\TT_\gamma}\iota_{X_\Delta}\beta^\parallel d\mu 
            &= 0 \qquad \text{and} \\
        \int_{\TT_\gamma}\iota_{X_\Delta}\beta^\perp d\mu 
            &= \int_{\TT_\gamma}\pi_\Delta^*(\iota_{X}dt) d\mu \\
            &= \int_{\TT_\gamma}1d\mu \\
            &= \mu(\TT_\gamma)
    \end{align*}
    The lemma follows by linearity and from the fact that~$[\beta^\perp]$ and~$[\beta^\parallel]$ form a basis of~$H^1(\TT_\gamma,\RR)$.
\end{proof}

\begin{lemma}\label{lemmaLinkOrbitSupported}
    Given a null-homologous algebraic multi-orbit~$\Gamma$ and a signed measure~$\mu\in\MMO$, the quantity~$\link_\Gamma^\Delta(\mu)$ depends only on~$\Gamma$ and~$\pi_\Delta^*(\mu)\in\MM^0_s(\phi)$. In particular, it does not depend on~$\Delta$.
\end{lemma}

\begin{proof}
    The previous lemma implies that $\link_\Gamma^\Delta(\mu)$ depends on $\pi_\Delta^*\mu$ and not on $\mu$ itself. To prove that it does not depend on $\Delta$, take $\Delta_1,\Delta_2$ two choices of finite unions of periodic orbits containing the support of $\Gamma$. By taking $\omega_1\in H^1(M_{\Delta_1},\RR)$ as in the definition of $\link_\Gamma^{\Delta_1}(\mu)$, and pushing it back to an element in $H^1(M_{\Delta_1\cup\Delta_2},\RR)$ it follows that we have 
    $$\link_\Gamma^{\Delta_1}(\mu)=\link_\Gamma^{\Delta_1\cup\Delta_2}(\mu)=\link_\Gamma^{\Delta_2}(\mu)$$
\end{proof}

\begin{definition}\label{defLinkPhif}
    For a null-homologous algebraic multi-orbit~$\Gamma$ in~$M$, and~$\mu\in\MM^0_s(\phi)$, we define $\link_\phi^f(\Leb_\Gamma,\mu)$ as $\link_\Gamma^\Delta(\nu)$ for any union of periodic orbits~$\Delta$ containing the support of~$\Gamma$, and any signed measure~$\nu\in\MMO$ satisfying~$\pi_\Delta^*\nu=\mu$.
\end{definition}

\begin{proposition}\label{propLinkOrbitSupported}
    Let~$\Delta$ be a finite union of periodic orbits. The map~$\link_\phi^f$ is continuous on $\{\nu\in\MM^0_s(\phi), \support(\mu)\subset\Delta\}\times\MM^0_s(\phi)$.
    Additionally, we have
    $$\link_\phi^f(\Leb_{\Gamma_1},\Leb_{\Gamma_2})=\link_\phi^{ff}(\Leb_{\Gamma_1},\Leb_{\Gamma_2})$$ for any null-homologous algebraic multi-orbit~$\Gamma_1,\Gamma_2$ with disjoint support. 
\end{proposition}

\begin{remark}\label{remarkSignedMeasureConvergence}
    A consequence of Banach-Steinhaus Theorem is the following. Given a sequence of signed measures~$\mu_n$ converging toward a signed measure~$\mu_\infty$, the total variation of~$\mu_n$ is bounded by some constant independent on~$n$.
\end{remark}

\begin{proof}
    Take a sequence~$(\nu_n,\mu_n)_n$ in~$\{\nu\in\MM^0_s(\phi), \support(\mu)\subset\Delta\}\times\MM^0_s(\phi)$, converging toward~$(\nu_\infty,\mu_\infty)$. Denote by~$\Gamma_n$ and~$\Gamma_\infty$ the algebraic multi-orbit corresponding to~$\nu_n$ and~$\nu_\infty$. As explained in 
    Remark \ref{remarkSignedMeasureConvergence}, 
    the total variations of $\nu_n$ and $\mu_n$ are bounded by some $C>0$.
    
    We take a sequence~$(\mu_{\Delta,n})_n$ in~$\MMO$ so that~$\pi_\Delta^*\mu_{\Delta,n}=\mu_n$ and the total variation of~$\mu_{\Delta,n}$ is at most~$C$. We can decompose~$\mu_{\Delta,n}$ as a difference of two invariant positive measures~$\mu_{\Delta,n}^+, \mu_{\Delta,n}^-$, so that each of them as total mass at most~$C$. By compactness, the sequence~$(\mu_{\Delta,n}^+, \mu_{\Delta,n}^-)_n$ accumulates toward a pair of positive measures~$(\mu_{\Delta,\infty}^+, \mu_{\Delta,\infty}^-)$, which satisfies~$\pi_\PP^*(\mu_{\Delta,\infty}^+-\mu_{\Delta,\infty}^-)=\mu_{\infty}$. So we have~$\link_\phi^f(\nu_\infty,\mu_\infty)=\link_{\Gamma_\infty}^\Delta(\mu_{\Delta,\infty}^+)-\link_{\Gamma_\infty}^\Delta(\mu_{\Delta,\infty}^-)$.

    Since~$\Delta$ is a finite union of periodic orbits, there exists a continuous map~$\Gamma\mapsto\omega_\Gamma$, so that~$\omega_\Gamma$ is in~$H^1(M_\Delta,\RR)$ and satisfies~$(\pi_\Delta)^*(\partial P_\Delta(\omega_\Gamma))=[\Gamma]\in H_1(\Delta,\RR)$. The quantities~$\link_{\Gamma_n}^\Delta(\mu_{\Delta,n}^\pm)=\cohomG{\mu_{\Delta,n}^\pm}\cdot \omega_{\Gamma_n}$ converges toward $\link_{\Gamma_\infty}^\Delta(\mu_{\Delta,\infty}^\pm)$. It follows that~$\link_\phi^f(\nu_n,\mu_n)$ converge toward $\link_\phi^f(\nu_\infty,\mu_\infty)$.

    The last point is a consequence of Poincaré duality. Take two algebraic multi-orbits~$\Gamma_1$ and~$\Gamma_2$, an algebraic multi-orbit~$\wt\Gamma_2=\sum_im_i\gamma_i$ in~$M_\Delta$ with~$\pi_*(\wt\Gamma_2)=\Gamma_2$. Take a closed 1-form~$\alpha$ on~$M_\Delta$ and a 2-chain~$S$ satisfying~$[S]=P_\Delta[\alpha]$ and~$(\pi_\Delta)_*(\partial S)=\Gamma_1$, then we have:
    \begin{align*}
        \int_{M_\Delta}\iota_{X_\Delta}\alpha d\Leb_{\wt\Gamma_2} 
            &= \sum_i m_i \int_{\gamma_i}\alpha \\
            &= \sum_i m_i \gamma_i\algcap S \qquad \text{by duality} \\
            &= \Gamma_2\algcap (\pi_\Delta)_* (S) \\
            &= \link_\phi^{ff}(\Leb_{\Gamma_1},\Leb_{\Gamma_2})
    \end{align*}
\end{proof}

As a consequence, one obtains the self-linking number of a null-homologous algebraic multi-orbit as the linking number between a Seifert surface of the algebraic multi-orbit and the stable direction, which is taken as a definition of the self-linking number by
Ghys \cite{Ghys09}.

\subsection{Linking function for a Markov partition}\label{section-LK-MarkovPartition}

\def\conj{\pi_\PP}

Suppose that~$\phi$ is an Anosov flow whose stable and unstable foliations are orientable. 
To prove
Theorem~\ref{theoremMeasureLinkingNumber},
we transfer the problem to a question on invariant signed measures for a finite sub-shift using a Markov partition $\PP$. Markov partitions makes the linking number computable with precise formulas. Lets us describe the strategy. 

Given an invariant signed measure~$\mu$ on $\PW$, we associate to~$\mu$ a signed measure~$\normalize(\mu)$ supported on finitely many periodic orbits, so that the linking number with $\mu-\normalize(\mu)$ can be computed as an integral. To see that, start with a periodic orbit $\gamma$ of the flow which intersects exactly twice a given rectangle $\PR$. As described in Section \ref{section-FriedSection}, there exists a transverse pair of pants bounded by $\gamma$ and by two shorter periodic orbits $\delta_1,\delta_2$. The pair of pants can be seen as a 4-gon $P$ inside the rectangle $\PR$ glued to two strips tangent to the flow ($P$ is not unique and depends on some choices). The linking number between $\gamma-\delta_1-\delta_2$ and another null-homologous algebraic multi-orbit $\Gamma$ (with disjoint support) is equal to the algebraic intersection between $\Gamma$ and $P$ (assuming that $\partial P$ is disjoint from the multi-orbit). Similarly, the linking number between $\gamma-\delta_1-\delta_2$ and any invariant signed measure $\nu$ is equal to 
$$\link_\phi^f(\Leb_\gamma-\Leb_{\delta_1}-\Leb_{\delta_2},\nu)=\pm\int_P\nu^\perp$$
where the sign depends on if $P$ is positively oriented or not (it depends on the relative position of $\gamma$ and $\delta_1,\delta_2$). A small combinatorial trick can be used to replace $P$ in the integral by a sub-rectangle $\PR_\gamma$ of $\PR$ which depends on $\gamma,\delta_1,\delta_2$ and on no more choices. 

Denote by $u$ a periodic $\PP$-word whose realization is $\gamma$ and by $\eta_u$ the $\sigma$-invariant signed measure supported by the orbit of $u$. There exists a map $\LinkForm\colon\PW^2\to\{-1,0,1\}$ which satisfies that $|\LinkForm(u,v)|=1$ if and only if $v\in\PR_\gamma$. Additionally the linking number can be computed as follows
$$\link_\phi^f(\Leb_\gamma-\Leb_{\delta_1}-\Leb_{\delta_2},\nu)=\int_{\PW^2}\LinkForm(x,y)d\eta_u\oplus\nu^\perp$$

Replace $\Leb_\gamma$ by a signed measure supported on finitely many orbits. A similar equation holds (where $\delta_i$ has to be replaced by a suitable collection of periodic orbits). We describe in Section~\ref{section-ReductionFiniteOrbit} a general construction. Given a $\sigma$-measure $\mu$ (playing the role of $\Leb_\gamma$), we define another $\sigma$-invariant signed measure $\normalize(\mu)$ supported on finitely many orbits (playing the role of $\sum_i\Leb_{\delta_i}$) for which we can apply the above formula. We actually use the equation to define the linking number $\link_\phi(\mu-\normalize(\mu),\nu)$ as the integral $\int_{\PW^2}\LinkForm(x,y)d\mu\oplus\nu^\perp$. This integral is proven to be continuous in $(\mu,\nu)$.

The signed measure $\normalize(\mu)$ is supported on the set of orbits which intersects the rectangle $\PR$ at most once (in some sense), so we view $\normalize$ as a reduction map. Applying it successively for all rectangles of the Markov partition yields a reduction map whose image is the set of signed measures supported by the (finite set of) periodic orbits which intersect all rectangles at most once. The set of normalized signed measures is finite dimensional, so the linking number for these signed measures is well defined (measures supported on finitely many orbits) and is continuous. One defines the linking number $\link_\phi(\mu,\nu)$ as the sum $\link_\phi(\mu-\normalize(\mu),\nu)+\link_\phi(\normalize(\mu),\nu)$, which is continuous in $(\mu,\nu)$.

\vvline 

We fix a Markov partition~$\PP$ for the flow~$\phi$. Denote by~$(\PR_1,\hdots,\PR_p)$ the family of Markov rectangles in~$\PP$. We fix two orientations on the stable and unstable foliations, and consider the induced vertical and horizontal orders on $\PW$.

Recall that we denote by~$\Theta_\PP\colon\MM_s(\sigma)\to\MM_s(\phi)$ the map sending $\sigma$-invariant signed measures to the corresponding $\phi$-invariant signed measures. For~$\mu\in\MM_s(\sigma)$, we denote by~$\cohom{\mu}=\cohom{\Theta_{\PP}(\mu)}\in H_1(M,\RR)$ the homology class corresponding to~$\mu$. 
We denote by $\MM^0_s(\sigma)$ the set of signed measure $\mu\in\MM_s(\sigma)$ satisfying $\cohom{\mu}=0$.

\begin{definition}
    For two signed measures~$\nu,\mu\in\MM^0_s(\sigma)$ supported by finitely many periodic orbits, we define the linking number of $\nu$ and $\mu$ by  $$\link_\sigma^{ff}(\nu,\mu)=\link_\phi^{f}(\Theta_{\PP}(\nu),\Theta_{\PP}(\mu))$$
\end{definition}

\paragraphc{Measure zero sets.}
Because of the Portmanteau Lemma, it is useful to identify some sets which have zero measures for all signed measures in~$\MM_s(\sigma)$. A set~$X\subset \PW$ is said to have \emph{uniformly zero measure} if for all~$\mu\in\MM_s(\sigma)$, one has~$|\mu|(X)=0$. Notice that it is enough to satisfy the property on probability measures.

\begin{lemma}[Portmanteau Lemma for almost continuous functions]\label{lemmaPortmanteau}
    Let~$X$ be a metrizable topological space endowed with the Borel algebra. Let~$Y\subset X$ be a measurable set and~$f\colon X\to\RR$ a bounded function, which is continuous on~$X\setminus Y$. Given a sequence~$(\mu_k)_k$ of signed measures on~$X$, weakly converging toward a signed measure~$\nu$, one has  $$\uplim_{k\to+\infty}\bigg|\int_X f d\mu_k - \int_X fd\nu\bigg|\leq 2|\nu|(Y)\cdot\lVert f\rVert_\infty$$
\end{lemma}

A particular case is when~$|\nu|(Y)=0$, where we have convergence of the integrals of~$f$. This lemma is adapted from a standard Portmanteau Lemma where the function is required to be continuous everywhere (see \cite[Theorem 13.16]{Klenke2008} for the usual Portmanteau Lemma). 

\begin{proof}
    Take a distance $d$ compatible with the topology on $X$. For~$t\in\RR^+$, denote by~$h_t^+\colon X\to \RR$ the map defined by $h_t^+(x)=\sup_{y\in X}(f(y)-td(x,y))$. The map~$h_t^+$ is well-defined since~$f$ is bounded, and~$t$-Lipschitz as supremum of~$t$-Lipschitz maps. Notice that we have~$f\leq h_t\leq \sup f\leq\lVert f \rVert_\infty$ for all~$t$. It follows that    $$\uplim_k\int_X f d\mu_k \leq \uplim_k \int_X h_t^+d\mu_k = \int_X h_t^+ d\nu$$

    For any~$x$ in $X\setminus Y$,~$f$ is continuous at~$x$ so when~$t$ goes to~$+\infty$,~$h_t^+(x)$ is decreasing and converges toward~$f(x)$. Hence for all~$k$ one has 
    \begin{align*}
        \uplim_{t\to+\infty}\int_X h_t^+d\nu 
            &= \int_X\big(\uplim_t h_t^+\big)d\nu \qquad \text{according to Fatou's Lemma}\\
            &= \int_Y\big(\uplim_t h_t^+\big)d\nu + \int_{X\setminus Y}\big(\uplim_t h_t^+\big)d\nu \\
            &\leq |\nu|(Y)\cdot \lVert f \rVert_\infty + \int_{X\setminus Y}fd\nu \\
            &\leq |\nu|(Y)\cdot \lVert f \rVert_\infty + \int_{X}fd\nu -\int_Yfd\nu \\
            &\leq \int_{X}fd\nu + 2|\nu|(Y)\cdot \lVert f \rVert_\infty
    \end{align*}
    
    It follows that  $$\uplim_k\int_X f d\mu_k - \int_X f  d\nu \leq 2|\nu|(Y)\cdot \lVert f \rVert_\infty$$ 
    
    Using the map~$h_t^-(x)=\inf_y( f(y)+td(x,y))$, the same arguments prove that  $$-\lowlim_k\bigg(\int_X f d\mu_k - \int_X f  d\nu\bigg) \leq 2|\nu|(Y)\cdot \lVert f \rVert_\infty$$
\end{proof}

We now prove two lemmas stating that two families of sets have uniformly zero measures. Given a rectangle~$\PR\in\PP$, define the sets:
$$E_\PR^{n+}=\{u\in\PW,u_n=\PR, \text{ and for all } k>n, u_k\neq \PR\}$$ 
\centerline{and}
$$E_\PR^{n-}=\{u\in\PW,u_n=\PR, \text{ and for all } k<n, u_k\neq \PR\}$$ 

We also define~$E_\PR^\pm=\cup_{n\in\ZZ}(E_\PR^{n-}\cup E_\PR^{n+})$, the set of~$\PP$-word containing a~$\PR$ at some index, and admitting a half sub-word without~$\PR$.

\begin{lemma}\label{lemmaZeroMeasure}
    The set~$E_\PR^\pm$ has uniformly zero measure.
\end{lemma}

\begin{proof}Take a signed measure~$\nu\in\MM_s(\sigma)$.
    By definition, we have~$\sigma^m(E_\PR^{n+})=E_\PR^{(n+m)+}$ and~$E_\PR^{n+}\cap E_\PR^{m+}=\emptyset$ when~$n\neq m$. Since~$\nu$ is~$\sigma$-invariant, one has~$\nu(\PW)\geq \sum_{n\in\ZZ} \nu(E_\PR^{n+})\geq \sum_n \nu(E_\PR^{0+})$. The signed measure~$\nu$ is finite so~$\nu(E_\PR^{n+})=0$ for all~$n$. Similarly~$\nu(E_\PR^{n-})=0$ and~$\nu(E_\PR^\pm)=0$.
\end{proof}

A second family of sets with uniformly zero measures are the union of stable and unstable leaves, after removing the periodic points. Take a point~$v\in\PW$. Its strong stable manifold~$W^s(v)$ is the set of point~$w$ such that~$d(\sigma^n(v),\sigma^n(w))\xrightarrow[n\to +\infty]{}0$. Notice that~$W^s(v)$ is the set of~$\PP$-words which coincide with~$v$ for large enough indexes.

\begin{lemma}\label{lemmaLeafZeroMeasure}
    Let~$w^s\subset\PW$ be a strong stable leaf. If~$w^s$ contains a periodic point~$u$, then it has a unique periodic point $u$ and~$w^s\setminus \{u\}$ has uniformly zero measure. Otherwise~$w^s$ has uniformly zero measure.
\end{lemma}

\begin{proof}
    Suppose that~$u,v\in w^s$ are periodic. Denote by~$n$ a common multiple of the periods of~$u$ and~$v$. Then the distance between~$\sigma^{kn}(u)=u$ and~$\sigma^{kn}(v)=v$ goes to zero when~$k$ goes to~$+\infty$. Hence we have~$u=v$. 

    Suppose that~$w^s$ contains a periodic point~$u$, of period~$N\geq 1$. Denote by~$w^s_n=\{v\in\PW, v_{n-1}\neq u_{n-1}, \forall k\geq n, v_k=u_k\}$, so that~$w^s\setminus \{u\}=\cup_n w^s_n$. The sets~$w^s_n$ are pairwise disjoints and~$\sigma^N(w^s_n)=w^s_{n-N}$. So using the same argument as in the previous proof, the sets~$w^s_{n-N}(v)$ have uniformly zero measures. Then~$w^s\setminus \{u\}$ has uniformly zero measure.

    Suppose now that~$w^s$ contains no periodic point. If there would exist two distinct integers~$n<m\in\ZZ$ and an element $u$ in~$\sigma^n(w^s)\cap\sigma^m(w^s)$, then $u$ would be $(m-n)$-periodic after a certain index. So the limit $\lim_k\sigma^{k(m-n)+n}(u)$ would exist and be in $w^s$. It is not possible so for any~$n\neq m$, we have~$\sigma^n(w^s)\cap\sigma^m(w^s)=\emptyset$. Then the previous argument shows that~$w^s$ has uniformly zero measure.
\end{proof}

\subsection{Reduction to finitely many orbits}\label{section-ReductionFiniteOrbit}

Let $\PP$ be a Markov partition and fix a rectangle $\PR\in\PR$. We define a reduction map, which takes an orbit $\gamma$ of the flow which intersects $\PR$, cuts it along its intersection with $\PR$, and transforms (using a Shadowing Lemma) each orbit arc into a periodic orbit intersecting $\PR$ exactly once. The precise definition is written below. It is somehow the opposite operation to the one we use in Section \ref{section-FriedSection}: we build a Fried section by merging two periodic orbits into one. 

Then we push back a measure by this reduction map. Given an invariant measure $\mu$, it yields another invariant measure $\normR(\mu)$ that is supported on orbits which intersects $\PR$ at most once (in some sense). Doing this successively on all rectangles of the Markov partition yields a measure $\normalize(\mu)$ which is supported on finitely many periodic orbit. We relate in the later section the reduction map with the linking number.

\vvline

Fix a rectangle~$\PR\in\PP$.
We call \emph{$\PR$-reduction map} the map~$\normR\colon\PW\to\PW$ defined as follows. For~$n,m\in\ZZ$ with~$n\leq m$, we denote~$S_{n,m}=\{u\in\PW,u_n=\PR=u_m, \PR\not\in u_{\intint{n+1}{m-1}}\}$. Then we define:

$$\normR(u)=
\begin{cases*}
    \wb{u_{|\intint{n}{m-1}}} & \text{ if~$u$ belongs to~$S_{n,m}$ for some~$n\leq 0$ and~$m\geq 1$} \\
    u & \text{ otherwise}
\end{cases*}$$
In the first case,~$\wb{u_{|\intint{n}{m-1}}}$ is the unique periodic~$\PP$-word of period~$m-n$, which coincides with~$u$ on the set of indexes~$\intint{n}{m}$. Notice that for all~$u\in\PW$,~$\normR(u)$ coincides with~$u$ for the indexes~$0$ and~$1$.
After giving a formula for~$\normR^*\mu$ for finitely supported signed measures~$\mu$, we prove the following. The map~$\normR^*$ induces a continuous map~$\MM_s(\phi)\to\MM_s(\phi)$, which additionally preserves the homology class of signed measures. 

Let~$u$ be a cyclic~$\PP$-word. Denote by~$\eta_u\in\MM_s(\sigma)$ the measure $\sum_{k=0}^{|u|-1}\caract{\sigma^k(\wb u)}$, which is the unique invariant measure supported on the orbit of~$\wb u$, normalized by~$\eta_u(\wb u)=1$.

\begin{lemma}\label{lemmaNormalizeMeasureOneOrbit}
    Let~$u$ be a cyclic~$\PP$-word with~$u(0)=\PR$. Denote by~$v_1\hdots v_n$ the unique family of cyclic~$\PP$-word which satisfies~$v_1v_2\cdots v_n=u$ and~$v_j^{-1}(\PR)=\{0\}$ for all~$j$. Then we have  $$\normR^*\eta_u=\sum_{j=1}^n \eta_{v_j}$$
\end{lemma}

\begin{proof}
    If~$k\in\intint{0}{|u|-1}$ satisfies~$\sum_{i=1}^{j-1}|v_i|\leq k<\sum{i=1}^{j}|v_i|$ for some~$j\in\intint{1}{n}$, then we have~$\normR\circ\sigma^k(\wb u)=\sigma^{k-\sum_{i=1}^{j-1}|v_i|}(\wb v_j)$. It follows that:
    \begin{align*}
        \normR^*\caract{\sigma^i(\wb u)}
            &= \caract{\normR\circ\sigma^k(\wb u)} \\
            &= \caract{\sigma^{k-\sum_{i=1}^{j-1}|v_i|}(\wb v_j)} 
    \end{align*}
    Hence we have~$\sum_{k=0}^{|u|-1}\normR^*\caract{\sigma^k(\wb u)}=\sum_j\sum_{k=0}^{|v_j|-1}\normR^*\caract{\sigma^k(\wb v_j)}=\sum_j\normR^*\eta_{v_j}$.
\end{proof}

\begin{lemma}\label{lemmafRcontinuous}
    The map~$\normR\colon\PW\to\PW$ is continuous on~$\PW\setminus E_\PR^{\pm}$, and is measurable.
\end{lemma}

\begin{proof}
    The topology on~$\PW$ is generated by the family of open sets of the form~$C_v=\{u\in\PW,u_{|\intint{-m}{m}}=v\}$ for all~$m\in\NN$ and for $\PP$-words~$v\colon\intint{-m}{m}\to\PP$. We determine the preimages of these sets by~$\normR$.

    When the letter~$\PR$ does not appear inside~$v$, we have~$\normR^{-1}(C_v)=C_v$. Suppose that~$v$ contains the letter~$\PR$ in position~$n_1\leq 0$ and~$n_2\geq 1$. We take~$n_1$ maximal and~$n_2$ minimal with these properties. When~$v$ is the restriction to~$[-m,m]$ of the periodic~$\PP$-word of period~$v_{|\intint{n_1}{n_2-1}}$, then we have~$\normR^{-1}(C_v)=C_{v_{|\intint{n_1}{n_2}}}$. When~$v$ is not such a restriction, we have~$\normR^{-1}(C_v)=\emptyset$.

    Lastly suppose that~$v$ contains the letter~$\PR$ in position~$n\leq 0$, but it does not contain this letter on positive positions. The last case is similar to this one. We can take~$n\leq 0$ maximal with~$v_n=\PR$. For~$u\in\PW$, the image~$\normR(u)$ is in $C_v$ if and only if one of the two following cases are satisfied: 
    \begin{itemize}
        \item~$u$ contains no~$\PR$ on positive positions, and~$u_{|\intint{-m}{m}}=v$. The set of such word~$u$ is a closed subset of~$\PW$, and is contained inside~$E_\PR^{\pm}$.
        \item For some~$k>m$,~$u$ contains a~$\PR$ on position~$k$, contains no~$\PR$ on positions~$[1,k-1]$,~$u_{|\intint{-m}{m}}=v$ and~$u_{|\intint{k-n-m}{k}}$ coincides with~$v_{|\intint{-m}{n}}$. For~$k$ fixed, the set of such words~$u$ is open.
    \end{itemize}

    Hence in any case,~$\normR^{-1}(C_v)$ is a union of open sets and of one closed set contained in~$E_\PR^\pm$. So it is measurable and continuous outside~$E_\PR^\pm$.
\end{proof}

\begin{lemma}\label{lemmafR*invariant}
    For all~$\mu\in\MM_s(\sigma)$,the pushed forward signed measure~$\normR^*\mu$ is invariant by~$\sigma$.
\end{lemma}

\begin{proof}
    We prove that~$\mu((\normR\circ\sigma)^{-1}(A))=\mu((\sigma\circ\normR)^{-1}(A))$ by considering the intersection between $(\normR\circ\sigma)^{-1}(A)$ and two subspaces of~$\PW$, which are~$S_{1,1}=\{u\in\PW,u_1=\PR\}$ and its complement.
    
    For any element~$u$ in $\PW\setminus S_{1,1}$, we have~$\sigma\circ\normR(u)=\normR\circ\sigma(u)$. Indeed, when~$u$ contains no coefficient equal to~$\PR$ for all positive indexes or for all non-negative indexes, then~$\normR(u)=u$ and~$\sigma(u)$ satisfies the same condition. When there exist~$n\leq 0$ maximal and~$m\geq 2$ minimal such that~$u_n=\PR=u_m$, it follows from the definition that~$\sigma\circ\normR(u)=\normR\circ\sigma(u)$. Hence  $$\mu((\normR\circ\sigma)^{-1}(A)\setminus S_{1,1})=\mu((\sigma\circ\normR)^{-1}(A)\setminus S_{1,1})$$

    Consider a~$\PP$-words~$u\in S_{1,1}$ with~$u_1=\PR$. Denote by~$E_\PR^\pm=\cup_{n\in\ZZ}(E_\PR^{n+}\cup E_\PR^{n-})$ the set of~$\PP$-words which contain the letter~$\PR$ somewhere and admit an infinite half sub-word without~$\PR$. According to Lemma~\ref{lemmaZeroMeasure}, we have~$\mu(E_\PR^\pm)=0$. Notice that~$S_{1,1}\setminus E_\PR^\pm=(\sqcup_{n\leq 0}(S_{n,1}))\setminus E_\PR^\pm=\sqcup_{m\geq 2}(S_{1,m})\setminus E_\PR^\pm$. For~$n\leq 0$, we define the homeomorphism~$\tilde g_n\colon S_{n,1}\mapsto S_{1,2-n}$ by~$\tilde g_n(u)=\sigma^{-n-1}(u)$. Notice that~$\normR\circ\sigma(\tilde g_n(u))=\sigma\circ\normR(u)$, so~$\tilde g_n$ restricts to a homeomorphism~$g_n\colon(\sigma\circ \normR)^{-1}(A)\cap S_{n,1}\setminus E_\PR^\pm\mapsto(\normR\circ\sigma)^{-1}\cap S_{1,2-n}\setminus E_\PR^\pm$. It follows that 
    \begin{align*}
        \mu\big((\sigma\circ\normR)^{-1}(A)\cap S_{1,1}\big) 
            &= \mu\big((\sigma\circ\normR)^{-1}(A)\cap (S_{1,1}\setminus E_\PR^\pm)\big)\\
            &= \sum_{n\leq 0}\mu\big((\sigma\circ\normR)^{-1}(A)\cap (S_{n,1})\setminus E_\PR^\pm\big)\\
            &= \sum_{n\leq 0}\mu\big(\sigma^{-n-1}\big((\sigma\circ\normR)^{-1}(A)\cap (S_{n,1})\setminus E_\PR^\pm\big)\big)\\
            & \text{ \hspace{1cm} since~$\mu$ is~$\sigma$-invariant}\\
            &= \sum_{n\leq 0}\mu\big((\normR\circ\sigma)^{-1}(A)\cap (S_{1,2-n})\setminus E_\PR^\pm\big)\\
            &= \mu\big((\normR\circ\sigma)^{-1}(A)\cap S_{1,1}\setminus E_\PR^\pm\big)\\
            &= \mu\big((\normR\circ\sigma)^{-1}(A)\cap S_{1,1}\big)
    \end{align*}
    
    Hence one has~$\sigma^*(\normR^*\mu)(A)=\mu(\normR^{-1}\circ\sigma^{-1}(A))=\mu(\sigma^{-1}\circ\normR^{-1}(A))=\mu(\normR^{-1}(A))=(\normR^*\mu)(A)$, concluding.
\end{proof}

\begin{proposition}\label{lemmafR*continuous}
    The map~$\normR^*\colon\MM_s(\sigma)\to\MM_s(\sigma)$ is continuous. 
\end{proposition}

\begin{proof}
    Consider a sequence of signed measures~$\mu_k\in\MM_s(\sigma)$ weakly converging toward the signed measure~$\nu\in\MM_s(\sigma)$. Take a continuous map~$g\colon\PW\to\RR$. The map~$g\circ\normR\colon\PW\to\RR$ is bounded, measurable and continuous outside the set~$E_\PR^\pm$ which has uniformly zero measure. According to the Portmanteau Lemma~\ref{lemmaPortmanteau}, the sequence of integrals~$\int gd(\normR^*\mu_k)=\int g\circ\normR d\mu_k$ converges toward~$\int gd(\normR^*\nu)$, which proves the continuity of~$\normR^*$.
\end{proof}

\begin{proposition}\label{propCohomologyInvariantByReducing}
    For every signed measure~$\mu\in\MM_s(\sigma)$, we have~$\cohom{\normR^*\mu}=\cohom{\mu}$. Additionally we have $|\normR^*\mu|(\PW)\leq |\mu|(\PW)$.
\end{proposition}

The statement on the total variation is clear. For~$x,y\in\PP$, denote by~$A_{x,y}$ the set~$\{v\in\PW,v_0=x, v_{1}=y\}$.
Notice that for all~$u\in\PW$,~$\normR(u)$ and~$u$ have the same coefficients at the indexes~$0$ and~$1$. Hence for~$\nu=\normR^*\mu-\mu$, we have~$\nu(A_{x,y})=0$ for all~$x,y\in\PP$. The proposition follows from the next lemma.

\begin{lemma}\label{lemmaCohomologyZeroCriterium}
    For any signed measure~$\nu\in\MM_s(\sigma)$ satisfying~$\nu(A_{x,y})=0$ for all~$x,y\in\PP$, we have~$\cohom{\nu}=0$. 
\end{lemma}

\begin{proof}
    Take a smooth closed 1-form~$\alpha$ on~$M$. For each Markov rectangle~$\PR\in\PP$ and its corresponding cuboid $\PC_\PR$, we denote by~$f_\PR\colon\PC_\PR\to\RR$ a smooth function whose differential coincides with~$\alpha$. We denote by~$i_\PR\colon\PR\to\partial^+\PC_\PR$ the map obtained by pushing $\PR$ along the flow.
    
    Recall that~$\sus$ is the suspension set over~$\PW$ and~$h_\PP\colon \sus\to M$ is the corresponding semi-conjugation. Denote by~$T\colon\PW\to\RR^+$ the suspension time and by 
    $$\sus^\PR=\{(x,t)\in\PW\times\RR, 0\leq t<T(x)\}$$
    the subset of~$\sus$ corresponding to the cuboid~$\PC_\PR$ (minus the boundary $\partial^+\PC_\PR$). We have 
    \begin{align*}
        \int_M\iota_{X}\alpha d\Theta(\nu) 
            &= \int_{\sus}(\iota_{X}\alpha)\circ h_\PP d(\nu\otimes\Leb) \\
            &= \sum_{\PR\in\PP} \int_{\sus^\PR}(df_\PR(X))\circ h_\PP d(\nu\otimes\Leb) 
    \end{align*}
    Denote by $\SR$ the set of $u\in\PW$ satisfying $u_0=\PR$. Using the local product of~$\sus$ one proves
    $$\int_{\sus^\PR}(df_\PR(X))\circ h_\PP d(\nu\otimes\Leb)=\int_{\Sigma_\PR}(f_\PR\circ i_\PR-f_\PR)\circ h_\PP(x,0)d\nu(x)$$
    Denote by~$f^+,f^-\colon\PW\to\RR$ the maps whose restrictions to~$\Sigma_\PR$ are given respectively by~$f_\PR\circ i_\PR\circ\pi_\PP$ and~$f_\PR\circ\pi_\PP$. Then:
    \begin{align*}
        \int_M\iota_{X}\alpha d\Theta(\nu) 
            &= \sum_{\PR\in\PP} \int_{\PWR}(f_\PR^+-f_\PR^-)d\nu \\
            &= \int_{\PW}(f^+-f^-)d\nu \\
            &= \int_{\PW}(f^+-f^-\circ\sigma)d\nu
    \end{align*}
    Notice that in a neighborhood of a point~$u\in\PW$,~$f^+$ and~$f^-\circ\sigma$ are obtained from integrating~$\alpha$ on the Markov sub-rectangle given by~$\partial^+\PC_{u_0}\cap u_1$. It follows that~$f^+-f^-\circ\sigma$ is constant on each set~$A_{x,y}$. Therefore the integral~$\int_{A_{x,y}}(f^+-f^-\circ\sigma)d\nu$ is zero. Since the non empty sets of the form $A_{x,y}$ give a partition of $\PW$, we have~$\int_{\PW}(f^+-f^-\circ\sigma)d\nu$. It follows that~$\int_M\iota_{X}\alpha d\Theta(\nu)=0$ for all $\alpha$, and that~$\cohom{\nu}=0$.
\end{proof}

Denote by~$(\PR_1,\hdots,\PR_p)$ the Markov rectangles in the Markov partition~$\PP$. Following Fried's terminology \cite{FriedCrossSection}, we call \emph{minimal loop} any periodic orbits of~$\sigma$ which intersects every rectangle at most once. That is~$u\in\PW$ is a minimal loop if it is periodic of period some~$T\geq 1$, and for any~$\PR\in\PP$ there is at most one~$i\in\intint{1}{T}$ such that~$u_i=\PR$. Denote by $\prim$ the set of minimal loops and notice that it has a finite cardinal.

Consider the product map~$\prod_{i=1}^p \norm{i}\colon\PW\to\PW$ (taken in any order).

\begin{lemma}\label{lemmaReductionSupport}
    For any~$\mu\in\MM_s(\sigma)$, the support of~$\big(\prod_{i=1}^p \norm{i}\big)^*\mu$ is included inside~$\prim$. 
\end{lemma}

\begin{proof}
    For a signed measure~$\mu\in\MM_s(\sigma)$, denote by~$\support(\mu)\subset\PW$ the support of~$\mu$.
    First notice that for any~$i$ and~$u\in\PW$,~$\norm{i}(u)$ is either equal to~$u$ or is a minimal loop. It follows that given a signed measure~$\mu\in\MM_s(\sigma)$, we have~$\support(\norm{i}^*\mu)\subset \support(\mu)\cup\prim$. Hence
    $$\support\Bigg(\Bigg(\prod_{i=1}^k \norm{i}\Bigg)^*\mu\Bigg)\cup\prim$$
     is a decreasing sequence in~$k$.
    
    Recall that $\SR=\{u\in\PW,u_0=\PR\}$. Denote by~$A_\PR$ the set of points whose orbit intersects $\Sigma_\PR$ in a single $\PP$-word. That is either it takes the value $\PR$ at most once, or it is periodic and its fundamental domain contains at most one $\PR$. Notice that $A_\PR$ is closed. For any~$u$ in~$\PW$, the element~$\normR(u)$ satisfies one of the following properties:
    \begin{enumerate}
        \item $\normR(u)$ is in $A_\PR$, when $u$ is in $A_\PR$ or when we have~$u_n=\PR=u_m$ for some~$n\leq 0,m\geq 1$,
        \item~$\normR(u)$ takes the value $\PR$ on exactly one of the two intervals $\intoint{0}$, $\intinto{1}$, when~$u$ satisfies the same property (and then $\normR(u)=u$).
    \end{enumerate}
 
    The set of~$u$ in the second case has uniformly zero measure according to Lemma~\ref{lemmaZeroMeasure}. So the support~$\support(\normR^*\mu)$ is included in~$A_\PR$. It follows from the above discussion that the support of~$\big(\prod_{i=1}^p \norm{i}\big)^*\mu$ is included in~$\cap_\PR A_\PR\cup\prim$. 

    It remains to see $\cap_\PR A_\PR=\prim$. It is clear that~$\prim$ is included in all $A_\PR$ and so in their intersection. Reciprocally a $\PP$-word $u\in\cap_\PR A_\PR$ is periodic (it cannot contains at most once every rectangle of $\PP$, so it is in the other case for all $\PR$). It follows that its fundamental domain contains at most once every rectangle of the Markov partition $u\in\prim$.
\end{proof}

\paragraphc{Density of finitely supported measures.}
For hyperbolic dynamical systems, density of measures supported on finitely many periodic orbits is a well known fact. We prove the density result for measures in a given homology classes (here for transitive symbolic dynamical systems and in Lemma \ref{lemmaFinitBoundedMeasureDense} for transitive Anosov flows).

\begin{lemma}\label{lemmaDensityFiniteOrbitPW}
    Given~$\nu\in\MM_s(\sigma)$, there exists a sequence of signed measures~$(\nu_n)_n$ in $\MM_s(\sigma)$, converging toward~$\nu$, so that for all~$n$ we have:
    \begin{itemize}
        \item~$\nu_n$ is finitely supported,
        \item~$\cohom{\nu_n}=\cohom{\nu}$,
        \item~$|\nu_n|(\PW)\leq |\nu|(\PW)$.
    \end{itemize}
    Additionally if~$\nu$ is positive,~$\nu_n$ can be taken positive.  
\end{lemma}

\begin{proof}
    For~$n\geq 1$, we denote by~$\PP_n$ the set of~$\PP$-words of the form~$u\colon\intint{-n}{n}\to\PP$. An element~$u\in\PP_n$ corresponds to a sub-rectangle of~$u_{-n}$, consisting of points admitting a short itinerary along~$\sigma^{-n}(u)$. Then~$\PP_n$ corresponds to the Markov partition of~$\phi$ given by these rectangles.    
    
    Denote by~$\Sigma_n$ the set of bi-infinite~$\PP_n$-words, and~$\sigma_n$ the shift map on~$\Sigma_n$. We have a natural map~$\pi_n\colon \Sigma_n\to\PW$, so that for~$v\in\Sigma_n$,~$\pi_n(v)$ and $i\in\ZZ$, we have~$(\pi_n(v))_i=(v_i)_0$. Notice that~$\pi_n$ is continuous bijective, so~$\nu\in\MM_s(\sigma_n)\to\pi_n^*\nu\in\MM_s(\sigma)$ is continuous and bijective. 

    For~$n$ and~$\PR\in\PP_n$, we denote by~$\normR^n\colon\Sigma_n\to\Sigma_n$ the~$\PR$-reduction map for the Markov partition~$\PP_n$ of~$\phi$.
    Given~$n$, we fix an order on~$\PP_n$ and denote by~$\normalize_n$ the product of the maps~$\norm{i}^n$ for all~$\PR_i\in\PP_n$. Here the product is ordered using the previously fixed order. 
    
    We claim that the sequence~$\pi_n\circ\normalize_n\circ\pi_n^{-1}$ converges uniformly toward the identity.
    Notice that for any~$u\in\Sigma_n$, we have~$\normalize_n(u)_0=u_0$. It follows that for~$v\in\PW$,~$\pi_n\circ\normalize_n\circ\pi_n^{-1}(v)$ coincides with~$v$ on~$\intint{-n}{n}$. Hence~$\pi_n\circ\normalize_n\circ\pi_n^{-1}(v)$ converges uniformly toward~$v$.
    
    Take~$\nu\in\MM_s(\sigma)$ and define the invariant signed measure~$\nu_n=(\pi_n\circ\normalize_n\circ\pi_n^{-1})^*\nu$. From the claim follows that the sequence~$\nu_n$ converges toward~$\nu$. The second and third points follow from 
    Proposition~\ref{propCohomologyInvariantByReducing}.
\end{proof}

\begin{lemma}\label{lemmaDensityFiniteMeasureDisjoint}
    Given a signed measure~$\nu\in\MM_s(\sigma)$ and a finite set~$\Delta\subset\PW$, there exists a sequence of signed measures~$(\nu_n)_n$ converging toward~$\nu$, so that for all~$n$:
    \begin{itemize}
        \item~$\nu_n$ has a finite support which is disjoint from~$\Delta$,
        \item~$\cohom{\nu_n}=\cohom{\nu}$.
    \end{itemize}
\end{lemma}

\begin{proof}
    Given a cyclic~$\PP$-word~$w$, denote by~$\eta_w$ the invariant measure~$\sum_{i=0}^{|w|-1}\caract{\sigma^i(\wb w)}$. Thanks to Lemma~\ref{lemmaDensityFiniteOrbitPW}, it is enough to prove the lemma when we have~$\nu=\eta_v$ for some cyclic~$\PP$-word~$v$. We fix a cyclic~$\PP$-word~$u$ containing all rectangles of~$\PP$, so that~$u_0=v_0$. We additionally suppose that~$u$ and~$v$ are not power of a common cyclic~$\PP$-word. The lemma follows from the two claims:
    \begin{claim}
        The sequence~$(\frac{1}{k}\eta_{v^ku})_{k\geq 0}$ converges weakly toward the measure~$\eta_v$.
    \end{claim}
    \begin{claim}
        The signed measure~$\frac{1}{k}(\eta_{v^ku}-\eta_u)$ is homologous to~$\eta_v$.
    \end{claim}
    
    \begin{claimproof}
        A large proportion of the words~$\sigma^i(\wb{v^ku})$ are close to the words of the form~$\sigma^j(\wb{v})$ when $k$ is large. More precisely, the word~$\sigma^{|v|n+m}(\wb{v^ku})$ coincides with the word~$\sigma^{m}(\wb{v})$ on the indexes $\intint{-|v|n}{|v|(k-n-1)}$, for $0\leq n<k$ and $0\leq m\leq |v|$. It follows that with $m$ fixed and for any sequence for which $n$ and $k-n$ goes to infinity,~$\sigma^{|v|n+m}(\wb{v^ku})$ converges toward~$\sigma^{m}(\wb{v})$. Decompose the measure $\eta_{v^ku}$ as a sum of three terms as follows
        
        \begin{align*}
            \frac{1}{k}\eta_{v^ku} 
                &= \frac{1}{k}\sum_{i=0}^{k|v|+|u|} \caract{\sigma^i(\wb{v^ku})} \\
                &= \frac{1}{k}\sum_{i=0}^{\sqrt{k}-1} \caract{\sigma^i(\wb{v^ku})} + \frac{1}{k}\sum_{i=\sqrt{k}}^{k|v|+|u|-\sqrt{k}} \caract{\sigma^i(\wb{v^ku})} + \frac{1}{k}\sum_{i=k|v|+|u|-\sqrt{k}+1}^{k|v|+|u|} \caract{\sigma^i(\wb{v^ku})} 
        \end{align*}

        The second term converges toward $\eta_v$. The total mass of the first and third terms are $\frac{1}{\sqrt{k}}$, so these terms converge toward zero.
    \end{claimproof}
    \begin{claimproof}
        According to Lemma \ref{lemmaNormalizeMeasureOneOrbit}, we have $\normR^*(\eta_{v^ku}-\eta_u)=k\normR^*\eta_v$, and the map $\normR^*$ preserves the homology according to Proposition \ref{propCohomologyInvariantByReducing}. 
    \end{claimproof}
\end{proof}

\subsection{Linking function}\label{section-LinkingFunction}

Recall the notation $\SR=\{u\in\PW,u_0=\PR\}$. In this section, we construct a measurable function~$\LinkFormR\colon \SR\times\SR\to\RR$. Its key feature (proven in the next section) is that for all signed measures~$\nu,\mu\in\MM_s(\sigma)$ supported on finitely many orbits, one has~$\link_\sigma(\nu-\normR^*\nu,\mu)=\int_{\SR^2}\LinkFormR d\nu\otimes\mu$. We prove another important property: the continuity of $\link_\sigma(\nu-\normR^*\nu,\mu)$ in~$(\nu,\mu)$.

Denote by~$\primR$ the set of cyclic~$\PP$-word~$u$ with~$u_0=\PR$ and~$u_i\neq\PR$ for all~$i\neq 0$. For~$u\in\primR$, We denote by~$\Su$ the set of $w\in\PW$ which satisfy $w_k=u_k$ for all $k\in\intint{0}{|u|}$. Notice that for~$w\in\Su$, the element~$\sigma^{|u|}(w)$ is in~$\SR$. For~$u\in\primR$, we define the two following measurable sets:

$$\Rec_{u}^+ = \{(w,x) \in\Su\times\SR, \wb u<_h x \text{ and } \sigma^{|u|}(w)\leq_v x<_v w \}$$
$$\Rec_{u}^- = \{(w,x) \in\Su\times\SR, \wb u<_h x \text{ and } w\leq_v x<_v \sigma^{|u|}(w)\}$$

\begin{figure}
    \begin{center}
        \begin{picture}(86,41)(0,0)
        \put(0,0){\includegraphics[width=80mm]{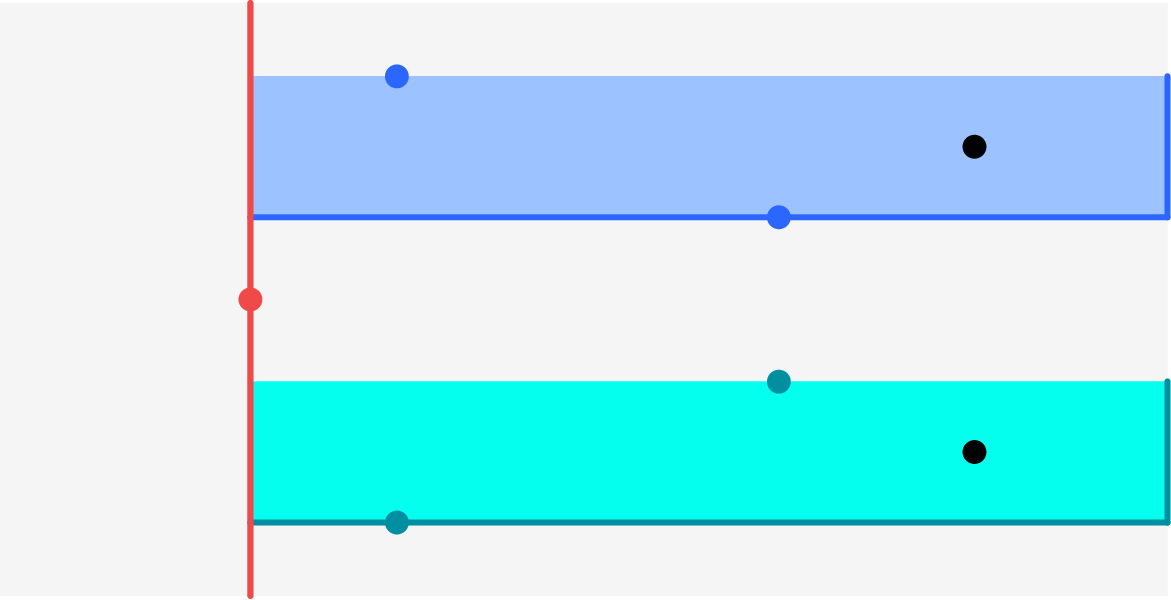}}
        \put(13.5,20){$u$}
        \put(52,23){$w$}
        \put(23,37){$\sigma^{|u|}(w)$}
        \put(52,16.5){$w'$}
        \put(23,1){$\sigma^{|u|}(w')$}
        \put(69,30){$x$}
        \put(69,9.5){$x'$}
        \put(81,30){$\Rec_u^-$}
        \put(81,9.5){$\Rec_u^+$}
        \end{picture}
    \end{center}
    \caption{Illustration of $(w,x)\in\Rec_{u}^-$ and $(w',x')\in\Rec_{u}^+$.}
    \label{figureLinkFunction}
\end{figure}

The definitions of the sets $\Rec_u^\pm$ can seem coming out of nowhere, but they are justified by a geometric picture in the next section. One may think about the Fried section in Section \ref{section-FriedSection}, obtained as the union of a transverse 4-gon and two tangent strips. The sets $\Rec^\pm_u$ play the role of the 4-gon, which we can use to compute linking number.

We call \emph{linking function} the function~$\LinkFormR\colon\PW^2\to\RR$ defined by  $$\LinkFormR=\sum_{u\in\primR}\big(\caract{\Rec_u^+}-\caract{\Rec_u^-}\big)$$
Here $\caract{X}$ is the characteristic map of the set $X$. From the definitions follows the next lemma.

\begin{lemma}\label{lemmaLinkFunctionFact}
    For any~$u\in\primR$, the following holds:
    \begin{enumerate}
        \item For any~$(w,x)\in \Rec_u^+$, we have~$w<_v\wb u$,
        \item For any~$(w,x)\in \Rec_u^-$, we have~$w>_v\wb u$,
        \item~$\Rec_u^+\cap \Rec_u^-=\emptyset$,
        \item~$\LinkFormR$ is measurable,
        \item~$ \lVert \LinkFormR \rVert_\infty=1$.
    \end{enumerate}
\end{lemma}

The next theorem and proposition are critical to define the linking number between two invariant signed measures.

\begin{theorem}\label{theoremLinkingFormProperty}
    For~$\nu,\mu\in\MM_s(\sigma)$ supported by finitely many
    periodic orbits and satisfying~$\cohom{\nu}=\cohom{\mu}=0$, we have  $$\link_\sigma^{ff}(\nu-\normR^*\nu,\mu)=\int_{\SR^2}\LinkFormR d\nu\otimes\mu$$
\end{theorem}

The proof of the theorem is postponed to next subsection.

\begin{proposition}\label{propContinuityIntLKForm}
    The map~$\MM_s(\sigma)\times \MM_s(\sigma)\to\RR$ given by~$(\nu,\mu)\mapsto\int_{\SR^2}\LinkFormR d(\nu\otimes\mu)$ is continuous.
\end{proposition}
 
Let us sketch the proof of the proposition. The key idea is to use the Portmanteau Lemma on $\LinkFormR$, so we would like to have that $\LinkFormR$ is continuous outside a set of small measure. $\LinkFormR$ is the difference between the characteristic maps of two sets, so it is continuous outside the boundary of these sets. Let us look at these boundaries.

As suggested in Figure \ref{figureLinkFunction}, for $u$ and $w\in\Sigma_u$ fixed, the projections of $\Rec_u^-$ and $\Rec_u^+$ onto the~$x\in\Sigma_\PR$ coordinate correspond to two sub-rectangles of $\Sigma_\PR$. The topology of $\Sigma_\PR$ is so that the vertical boundaries on the right of the figure are not actually boundaries inside $\PW$ (the cylinders in $\PW$ are both closed and open, so they have no boundary). The vertical boundaries on the left are included in unstable leaves (minus the periodic point) so they have mass zero. The vertical boundaries are included in stable leaves, so they have mass zero outside the potential periodic points ($w$ and $\sigma^{|u|}(w)$ when $w$ is periodic). When they contain periodic points, we cut in $\Rec^\pm_u$ a small rectangle around $w$ and glue it to $\sigma^{|u|}(w)$ (we use $\sigma^{|u|}$ for the identification). Using the argument above, the boundary of the small rectangles have mass zero, and we have removed two points with positive mass. The mass of the sets remain the same since the measures are $\sigma$-invariant.

The general argument has to take into account two more difficulties. First $w$ varies so the boundary are slightly more difficult to handle. Secondly $u$ varies too. We prove that the mass is mostly concentrated on finitely many values of $u$. That is the union of $\Rec^\pm_u$ have small masses when taking all values of $u$ excepts for finitely many ones. The ideas are precisely stated in the three following lemmas.

\begin{lemma}[Estimation Lemma 1]\label{lemmaTechLink1}
    For any~$\epsilon>0$, there exists~$A\in\NN$, such that for all signed measures~$\nu,\mu\in\MM_s(\sigma)$, we have:
    $$|\nu\otimes\mu|\Bigg(\bigcup_{u\in\primR, |u|\geq A}\Rec_u^+\Bigg)\leq \epsilon|\nu\otimes\mu|(\PW^2)$$
\end{lemma}

For~$u$ in~$\primR$ and~$k\geq 0$, we denote by~$u^k$ the concatenation of~$k$ copies of~$u$, with~$u^0=\PR$ by convention.
We define the sets  $$\Sigma_{u,k}=\Sigma_{u^k}\setminus\Sigma_{u^{k+1}}$$ \centerline{and} $$\Rec_{u,k}^+=\{(w,x)\in \Rec_u^+, w\in\Sigma_{u,k}\}$$

Recall that we denote $W^s_l(v)=\{w\in\PW,\text{ for all }k\geq 0, w_k=v_k\}$.
Notice that~$\Su\setminus W^s_l(\wb u)$ is the union of~$\Sigma_{u^k}\setminus\Sigma_{u^{k+1}}$ for all~$k\geq 1$. Given~$w$ in $W^s_l(\wb u)$, we do not have~$\sigma^{|u|}(w)<_vw$ nor $\sigma^{|u|}(w)>_vw$, so no~$(w,x)$ can be in~$\Rec_u^+\cup\Rec_u^-$. Therefore the sets~$\Rec_u^+$ and $\Rec_u^-$ are the unions of the set~$\Rec_{u,k}^+$ and $\Rec_{u,k}^-$ for all~$k\geq 1$.
    
\begin{lemma}[Estimation Lemma 2]\label{lemmaTechLink2}
    For all~$\epsilon>0$ and~$u\in\primR$, there exists~$B\in\NN$ such that for all signed measures~$\nu,\mu\in\MM_s(\sigma)$, we have  $$|\nu\otimes\mu|\Bigg(\bigcup_{k\geq B}\Rec_{u,k}^+\Bigg)\leq \epsilon|\nu\otimes\mu|(\PW^2)$$
\end{lemma}

\begin{lemma}[Estimation Lemma 3]\label{lemmaTechLink3}
    Given~$u\in\primR$ and~$k\geq 1$, the quantity~$(\nu\otimes\mu)(\Rec_{u,k}^+)$ is continuous in~$(\nu,\mu)\in\MM_s(\sigma)\times \MM_s(\sigma)$.
\end{lemma}

\begin{proof}[Proof of Proposition~\ref{propContinuityIntLKForm}]
    Assume the three estimation lemmas. 
    Let~$(\nu_n,\mu_n)_n$ be a sequence of pairs of signed measures in~$\MM_s(\sigma)$, converging toward a pair of signed measures~$(\nu_\infty,\mu_\infty)$. We should prove that~$\int_{\SR^2}\LinkFormR d(\nu_n\otimes\mu_n)$ converges toward~$\int_{\SR^2}\LinkFormR d(\nu_\infty\otimes\mu_\infty)$.
    
    According to Remark \ref{remarkSignedMeasureConvergence}, 
    the total variations of the signed measures $\nu_n\otimes\mu_n$ are bounded by some $C>0$. Take~$\epsilon>0$ and~$A\in\NN$ given by the first estimation lemma for~$\epsilon$. Denote by~$N\in\NN$ the number of elements~$u\in\primR$ with~$|u|<A$. Take~$B\in\NN$ the largest integer given by the second estimation lemma, applied for all~$u$ with~$|u|<A$ and for~$\epsilon/N$. We have 
    \begin{align*}
        \uplim_{n\to +\infty}\Bigg(\nu_n\otimes\mu_n\Bigg(\bigcup_{u\in\primR}\Rec_{u}^+\Bigg)\Bigg) 
            &=\uplim_n \Bigg(\nu_n\otimes\mu_n\Bigg(\bigcup_{\substack{|u|< A, k<B}}\Rec_{u,k}^+\Bigg)\Bigg. \\
            &\hspace{2cm} + \nu_n\otimes\mu_n\Bigg(\bigcup_{|u|\geq A}\Rec_{u}^+\Bigg) \\ 
            &\hspace{2cm} \Bigg.+ \nu_n\otimes\mu_n\Bigg(\bigcup_{\substack{|u|< A, k\geq B}}\Rec_{u,k}^+\Bigg) \Bigg)\\ 
            &\leq \uplim_n \Bigg(\nu_n\otimes\mu_n\Bigg(\bigcup_{\substack{|u|< A, k<B}}\Rec_{u,k}^+\Bigg)\Bigg) \\
            & \hspace*{2cm} + C\epsilon + N\cdot\frac{C\epsilon}{N} \\
            &\leq \nu_\infty\otimes\mu_\infty\Bigg(\bigcup_{\substack{|u|< A, k<B}}\Rec_{u,k}^+\Bigg) + 2C\epsilon \\
            & \qquad \text{ according to the estimation Lemma 3} \\
            &\leq \nu_\infty\otimes\mu_\infty\Bigg(\bigcup_u\Rec_{u}^+\Bigg) +4C\epsilon \text{ for all } \epsilon>0 \\
            &\leq \nu_\infty\otimes\mu_\infty\Bigg(\bigcup_u\Rec_{u}^+\Bigg)  \\
        \lowlim_n\Bigg(\nu_n\otimes\mu_n\Bigg(\bigcup_{u\in\primR}\Rec_{u}^+\Bigg)\Bigg) 
        &\geq \nu_\infty\otimes\mu_\infty\Bigg(\bigcup_u\Rec_{u}^+\Bigg) 
    \end{align*}
    
    Hence~$(\nu_n\otimes\mu_n)(\cup_u \Rec_{u}^+)$ converges toward~$(\nu_\infty\otimes\mu_\infty)(\cup_u \Rec_{u}^+)$. The same statement holds for~$\Rec_u^-$. Therefore the integrals~$\int_{\SR^2}\LinkFormR d(\nu_n\otimes \mu_n)$ converge toward~$\int_{\SR^2}\LinkFormR d(\nu_\infty\otimes \mu_\infty)$.
\end{proof}

Now we prove the estimation lemmas.

\begin{proof}[Proof of Lemma~\ref{lemmaTechLink1}]  
    Given $A\in\NN$, denote by~$L_A$ the union of the sets~$\Su$ for~$u\in\primR$ with~$|u|\geq A$. Notice that if~$w$ is in~$L_A$, then~$w_0=\PR$ and~$w_k\neq \PR$ for all~$k\in\intint{1}{A-1}$. Hence for two distinct indexes~$i,j\in\intint{0}{A-1}$, the set~$\sigma^i(L_A)$ and~$\sigma^j(L_A)$ are disjoint. For a signed measure~$\nu\in\MM_s(\sigma)$, we have:
    $$A|\nu|(L_A)=\sum_{i=0}^{A-1}|\nu|(\sigma^i(L_A))\leq |\nu|(\PW)$$

    so~$|\nu|(L_A)\leq\frac{|\nu|(\PW)}{A}$. We take~$A$ satisfying~$1\leq A\epsilon$. Then we have:
        $$|\nu\otimes\mu|\Bigg(\bigcup_{|u|\geq A} \Rec_u^+\Bigg) \leq |\nu\otimes\mu|\big(L_A\times\PW\big) \leq \epsilon|\nu\otimes\mu|(\PW^2)$$
\end{proof}    

\begin{proof}[Proof of Lemma~\ref{lemmaTechLink2}] 
    Take~$(w,x)$ in~$\Rec_{u,k}^+$. Then both~$w$ and~$\sigma^{|u|}(w)$ are in~$\Sigma_{u^{k-1}}\setminus\Sigma_{u^{k+1}}$. By definition, we have~$\sigma^{|u|}(w)\leq_v x<_v w$, so we also have~$x\in\Sigma_{u^{k-1}}$. According to Lemma~\ref{lemmaLinkFunctionFact}, we have~$x<_vw<_v\wb u$. So if~$x$ is in~$\Sigma_{u^{k+1}}$,~$u$ is in~$\Sigma_{u^{k+1}}$, so~$w$ is also in~$\Sigma_{u^{k+1}}$. Since~$w$ is not is this set,~$x$ lies to~$\Sigma_{u^{k-1}}\setminus\Sigma_{u^{k+1}}$.
    Take $B>1$. Notice that we have~$\sigma^{|u|}(\Sigma_{u^k})\subset\Sigma_{u^{k-1}}$, so given $i$ with~$0\leq 2i< B-1$, one has:
    $$(\id\times\sigma^{2i})(\Rec_{u,k}^+)\subset (\Sigma_{u^k}\setminus\Sigma_{u^{k+1}})\times \big(\Sigma_{u^{k-1-2i}}\setminus\Sigma_{u^{k+1-2i}}\big)$$
    
    It follows that the sets~$(\id\times\sigma^{2i})(\cup_{k\geq B}\Rec_{u,k}^+)$ for~$0\leq 2i< B-1$ are pairwise disjoint. For any~$\nu,\mu\in\MM_s(\sigma)$, we have~$(\nu\otimes\mu)\big((\id\times\sigma^{2i})(\cup_{k\geq B}\Rec_{u,k}^+)\big)=(\nu\otimes\mu)\big(\cup_{k\geq B}\Rec_{u,k}^+\big)$, so using a similar argument as above, we have:
    $$(\nu\otimes\mu)\big(\cup_{k\geq B}\Rec_{u,k}^+\big)\leq \frac{|\nu\otimes\mu|(\PW^2)}{\Big\lfloor\frac{B-1}{2}\Big\rfloor}$$
\end{proof}

Denote by~$S_k=\big\{w\in\Sigma_{u,k}, w<_v \wb u, \wb u<_h w\big\}$ for~$k\geq 1$. 

\begin{lemma}\label{lemmaSkOrder}
    We have~$S_{k}<_v S_{k+1}<_v\wb u$ for all~$k\geq 1$.
\end{lemma}

\begin{proof}
    Take~$w_1\in S_{k}$ and~$w_2\in S_{k+1}$. Suppose that~$w_2\leq_v w_1$. So~$w_1$ lies vertically between~$w_1$ and $\wb u$, which start by~$u^{k+1}$. By definition of the vertical order,~$w_1$ starts by~$u^{k+1}$, which is not possible since~$\Sigma_{u,k}=\Sigma_{u^k}\setminus\Sigma_{u^{k+1}}$. 
\end{proof}

We prove the estimation Lemma 3: $\nu\otimes\mu(\Rec_{u,k}^+)$ is continuous in $(\nu,\mu)\in\MM_s(\sigma)^2$. The idea is to use the Portmanteau Lemma. The boundary of $\Rec_{u,k}^+$ can have a non-zero mass, so we need a first step before using the Portmanteau Lemma. For that we notice that a signed measure $\nu\otimes\mu$ charges $\partial\Rec_{u,k}^+$ only on pairs $(z,z)$ and $(z,\sigma^{|u|}(z))$, where $z$ is periodic. We cut a neighborhood of $(z,z)$ in $\Rec_{u,k}^+$ and glue it on $(z,\sigma^{|u|}(z))$, to remove these points from the boundary of $\partial\Rec_{u,k}^+$. Then we use the Portmanteau Lemma.

\begin{proof}[Proof of Lemma~\ref{lemmaTechLink3}] 
    Let~$(\nu_n,\mu_n)_n$ be a sequence of pairs of signed measures in~$\MM_s(\sigma)$ converging toward a pair of signed measure~$(\nu_\infty,\mu_\infty)$. We use a Portmanteau Lemma to prove the convergence of the integral. 
    Take~$\epsilon>0$ and $k\geq 1$. We can take a subset~$X_\epsilon\subset S_k$ for which the sum of~$|\nu_\infty|(\{x\})$ for~$x\in S_k\setminus X_\epsilon$ atom of~$|\nu_\infty|$, is less that~$\epsilon$. That is  $$\sum_{x\in S_k\setminus X_\epsilon}|\nu_\infty|(\{x\})\leq \epsilon$$
    
    For each~$z\in X_\epsilon$, we take a small neighborhood~$T_z\subset \PW$ of~$z$. We can take it such that the following are satisfied:

    \begin{claim}
        There exist neighborhoods~$T_z$ of the points~$z\in X_\epsilon$ satisfying the following:
        \begin{enumerate}
            \item~$\partial T_z=\emptyset$, 
            \item the sets~$T_z$ and the sets~$\sigma^{|u|}(T_z)$ all~$z\in X_\epsilon$ are pairwise disjoint,
            \item~$\sigma^{|u|}(T_z)<_v\sigma(T_z)$,
            \item $\sigma^{|u|}(T_z)>_h\wb u$ for all $z\in X_\epsilon$.
        \end{enumerate}
    \end{claim}

    \begin{claimproof}
    We take~$T_z=\{w\in\PW, \text{ for all } |k|\leq N, w_k=z_k\}$ for some~$N>0$. It is open and close, so it has no boundary. If we prove the points 2, 3 and 4 when replacing $T_z$ by $z$, then points 2 to 4 will be satisfied when taking $N$ large enough.
    If~$z$ is in~$X_\epsilon\subset S_k$, then~$\sigma^{|u|}(z)$ is not in~$\Sigma_{u^k}$, so it is not in~$S_k$. Hence the set~$X_\epsilon$ and~$\sigma^{|u|}(X_\epsilon)$ are pairwise disjoint. 
    Since~$z$ is in~$\Sigma_u$ and~$z<_v\wb u$, we have~$\sigma^{|u|}(z)<_vz$ according to Lemma \ref{lemmaSkOrder}. For $z\in X_\epsilon$, we have $\wb u<_hz$, and $\wb u$ and $z$ coincide on $\intint{0}{|u|}$. Hence we have $\wb u<_h\sigma^{|u|}(z)$. Therefore the points 2 to 4 are satisfied for large $N$.
    \end{claimproof}

    \vline

    We define the sets  $$T=\Bigg(\bigcup_{z\in X_\epsilon}T_z\times T_z\Bigg)\bigcap \Rec_{u,k}^+$$ \centerline{and} $$\Rec'=\big(\Rec_{u,k}^+\setminus T\big)\bigcup \Big(\id\times\sigma^{|u|}\Big)(T)$$ 
    
    \begin{claim}
       ~$\Rec_{u,k}^+\setminus T$ and~$(\id\times\sigma^{|u|})(T)$ are disjoint, so~$\nu_\infty\otimes\mu_\infty(\Rec_{u,k}^+)=\nu_\infty\otimes\mu_\infty(\Rec')$ 
    \end{claim}

    \begin{claimproof}
        When~$(w,x)$ is in~$\Rec_{u,k}^+$ we have~$\sigma^{|u|}(w)\leq_v x$. When~$(w,x)$ is in~$(\id\times\sigma^{|u|})(T)$, we have~$(w,\sigma^{-|u|}(x))$ is in~$T\subset \Rec_{u,k}^+$, so~$\sigma^{-|u|}(x)<_vw$ and~$x<_v\sigma^{|u|}(w)$. 
    \end{claimproof}

    We now prove that the boundary of~$\Rec'\subset\PW^2$ has small measure for~$\nu_\infty\otimes\mu_\infty$. Notice that for~$z$ in $X_\epsilon$, the point~$(z,z)$ is in the interior of~$T$, so it is not in the adherence of~$\Rec'$. 

    \begin{claim}
        The set~$T_z\times \sigma^{|u|}(T_z)$ is included in~$\Rec'$.
    \end{claim}

    Assume the claim. Denote $x\equiv_v w$ when we have $x\leq_vw\leq_vx$, that is $x\in W^s_l(w)$. The boundary~$\partial \Rec'$ is contained in union of the three sets:
    \begin{enumerate}
        \item~$\{(w,x)\in \Sigma_{u,k}\times\SR,x\equiv_h \wb u, \wb u<_hx\}$
        \item~$\{(w,x)\in (\Sigma_{u,k}\setminus X_\epsilon)\times\SR, x\leq_v\wb u, \wb u\leq_h x, x\equiv_v w\}$
        \item~$\{(w,x)\in (\Sigma_{u,k}\setminus X_\epsilon)\times\SR, x\leq_v\wb u, \wb u\leq_h x, x\equiv_v \sigma^{|u|}(w)\}$
    \end{enumerate}

    The mass for~$|\nu_\infty\otimes\mu_\infty|$ of the first set is zero 
    (see Lemma \ref{lemmaLeafZeroMeasure}). 
    For the two other sets, we have:
    \begin{align*}
        |\nu_\infty\otimes\mu_\infty|&\big(\{(w,x)\in (\Sigma_{u,k}\setminus X_\epsilon)\times\SR, x\leq_v\wb u, \wb u\leq_h x, x\equiv_v w\}\big) \\
            &= \int_{w\in\Sigma_{u,k}\setminus X_\epsilon}\bigg(\int_{x\in W^s_l(w)}\caract{x\leq_v\wb u, \wb u\leq_h x}d|\mu_\infty|(x)\bigg)d|\nu_\infty|(w) \\
            &= \int_{\substack{w\in\Sigma_{u,k}\setminus X_\epsilon\\ w\text{ periodic}}}\bigg(\int_{x\in W^s_l(w)}\caract{x\leq_v\wb u, \wb u\leq_h x}d|\mu_\infty|(x)\bigg)d|\nu_\infty|(w) \text{ according to Lemma~\ref{lemmaLeafZeroMeasure}} \\
            &= \int_{\substack{w\in\Sigma_{u,k}\setminus X_\epsilon\\ w\text{ periodic}}}\caract{w\leq_v\wb u, \wb u \leq_h w}|\mu_\infty|(\{w\})d|\nu_\infty|(w) \\
            &= \int_{w\in S_k\setminus X_\epsilon}|\mu_\infty|(\{w\})d|\nu_\infty|(w) \\
            &\leq \epsilon|\nu_\infty|(\PW)
    \end{align*}

    A similar argument shows that the third sets has mass at most~$\epsilon|\nu_\infty|(\PW)$.
    It follows from the Portmanteau 
    Lemma~\ref{lemmaPortmanteau} 
    that we have  $$\uplim_n\big|\nu_n\otimes\mu_n(\Rec')-\nu_\infty\otimes\mu_\infty(\Rec')\big|\leq 2\epsilon|\nu_\infty|(\PW)$$ 
    This holds for all~$\epsilon>0$, and~$\nu\otimes\mu(\Rec')=\nu\otimes\mu(\Rec_{u,k}^+)$ for all~$\nu,\mu\in\MM_s(\sigma)$. Therefore~$\nu_n\otimes\mu_n(\Rec_{u,k}^+)$ converges in $n$ toward~$\nu_\infty\otimes\mu_\infty(\Rec_{u,k}^+)$.

    \begin{claimproof}
        Recall that for $x\in\sigma^{|u|}(T_z)$, we have $\wb u<_hx$.
        If~$(w,x)$ is in~$T_z\times \sigma^{|u|}(T_z)$ and satisfies~$w\leq_v \sigma^{-|u|}(x)$, then we have~$\sigma^{|u|}(w)\leq_v x<_v w$. The second equation is due to~$\sigma^{|u|}(T_z)<_v\sigma(T_z)$. In that case,~$(w,x)$ lies in~$\Rec_{u,k}^+$ and not in~$T$, so~$(w,x)\in \Rec'$. 
        
        If~$(w,x)\in T_z\times \sigma^{|u|}(T_z)$ satisfies~$\sigma^{-|u|}(x)<_v w$, we also have~$\sigma^{|u|}(w)\leq_v\sigma^{-|u|}(x)$. Indeed~$\sigma^{|u|}(w)$ belongs to~$S_{k-1}$,~$\sigma^{-|u|}(x)$ belongs to~$S_k$, and~$S_k>_vS_{k-1}$ according to Lemma~\ref{lemmaSkOrder}. Hence~$(w,\sigma^{-|u|}(x))$ is in both~$\Rec_{u,k}^+$ and~$T_z^2$, so it lies in~$T$. Hence~$(w,x)$ is in~$\Rec'$.
    \end{claimproof}
\end{proof}

\subsection{Adapted cells}

We prove Theorem~\ref{theoremLinkingFormProperty}, that is the linking number~$\link_\sigma^{ff}(\normR^*\nu - \nu,\mu)$ is equal to~$\int_{\SR^2}\LinkFormR d\nu\otimes\mu$ when $\mu$ and $\nu$ are finitely supported. For that, recall that the Fried sections $S$ are constructed as the union of a 4-gon $P$ in a Markov rectangle $\PR$ and two strips tangent to the flow. The linking number between $\partial S$ and a null-homologous algebraic multi-orbit $\Gamma$ is equal to (plus or minus) the algebraic intersection between $P$ and $\Gamma$. Here we replace $P$ with a sub-rectangle of $\PR$ which depends on no choices (contrary to the edges $P$ which depend on arbitrary choices). We define an adapted 2-cell (an nice enough cell of $\PR$) which generalize the 4-gon $P$ and the sub-rectangle. It allows us to precisely compute the linking number. 

\vvline

Denote by~$\Delta^n$ the~$n$-dimensional simplex. We call~$n$-chains, the abstract finite sums of continuous maps~$\Delta^n\to\PR$, with coefficients in~$\RR$. We also denote by~$\partial$ the boundary map. A differential embedding~$\Delta^1\to\PR$ is said adapted if it is either transverse to the stable and the unstable foliations, or is tangent to one of them. A 1-chain is said \emph{adapted} if it is the sum of adapted embeddings~$\Delta^1\to\PR$. For~$\lambda_i\in\RR$ and~$f_i\colon \Delta^2\to\PR$, the 2-chain~$c=\sum_{i=1}^n\lambda_i f_i$ is said adapted if the maps~$f_i$ are orientation preserving immersions and~$\partial f_i$ are adapted. Here~$\Delta^2$ is endowed with the anti-clockwise orientation.

Recall that $W^s_\PR$ and $W^u_\PR$ are the stable and unstable foliations on $\PR$.
We identify the bi-foliated rectangle~$(\PR,W^s_\PR,W^u_\PR)$ with a sub rectangle of~$(\RR^2,H,V)$, with a~$\Class^1$ diffeomorphism preserving the orientations of the foliations. 
Let~$C$ be an adapted 2-chain. We define the \emph{support function} of~$C$ has~$\Sup_C\colon\PR\to\RR$ defined for~$x\in\PR\subset\RR^2$ by
$$\Sup_C(x)=\lim_{t\to 0^+}C\algcap (x+(-t,t^2))$$
In other words when $C$ is an immersed polygon (positively oriented), $\Sup_C(x)$ is equal to 1 when one the following is true:
\begin{itemize}
    \item $x$ is in the interior of $C$,
    \item $x$ is on a edge $e$ of $C$ which has a non-negative slope and $C$ is adjacent to $e$ by above,
    \item $x$ is on a edge $e$ of $C$ which has a negative slope and $C$ is adjacent to $e$ by below,
    \item $x$ is on a corner of $C$, and $C$ contains an angular sector based at $x$ and of angle $(\pi-\epsilon,\pi)$ for some small $\epsilon>0$.
\end{itemize}
The value of $\Sup_C$ on the boundary of $C$ is made to be additive in $C$. Let~$c$ be an adapted 1-cycle. We define~$\Sup_c=\Sup_C$ for any 2-chain~$C$ satisfying~$\partial C=c$.

\begin{remark}\label{lemmaTrivialSupport}
    Given an adapted 2-cell $C$, if either~$\partial C=0$ or~$C$ is contained in finitely many stable leaves and one unstable leave, then~$\Sup_C=0$.
\end{remark}

\paragraphc{Reducing a cyclic~$\PP$-word}

We fix a primitive cyclic~$\PP$-word~$r$ with~$r_0=\PR$. Denote by~$\eta_r$ the measure~$\eta_r=\sum_{k=0}^{|r|-1}\caract{\sigma^k(\wb r)}$. We construct an adapted 2-cell to prove the Theorem~\ref{theoremLinkingFormProperty} for the measure~$\eta_r$.
Write~$r=s_1\circ\hdots\circ s_n$ where all~$s_i$ are cyclic words such that~$(s_i)^{-1}(\PR)=\{0\}$. According to Lemma~\ref{lemmaNormalizeMeasureOneOrbit}, we have~$\normR^*\eta_r=\sum_{j=1}^n \eta_{s_j}$. 

We suppose that~$n\geq 2$. For~$i\in\intint{0}{n-1}$, we denote by~$r_i=s_i\circ\hdots\circ s_n\circ s_1\circ\hdots\circ s_{i-1}$. Since~$n>1$,~$\wb r_i$ and~$\wb s_i$ are different. So thanks to 
Lemma~\ref{lemmaDistinctRealization},
we can find a smooth curve~$a_i\colon[0,1]\to\PR$, from~$\pi_\PP\big(\wb{r_i}\big)$ to~$\pi_\PP\big(\wb{s_i}\big)$, which is additionally transverse to the stable and unstable foliations. 
Since the set of periodic orbits is countable, we can take the curve~$a_i$ such that its interior is disjoint from all period orbits.

Recall that when~$x\in\PR$ admits an itinerary starting with~$u_0,u_1,\cdots ,u_{|u|}=u_0$, then~$T^u(x)$ is the time realizing this itinerary on~$x$, and~$\sigma^u(x)=\phi_{T^u(x)}(x)$.
Denote by~$b_i\colon[0,1]\to\PR$ the curve given by~$b_i(t)=\sigma^{s_i}\circ a_i(1-t)$, and define the adapted 2-chain~$c=\sum_i(a_i+b_i)$.

\begin{lemma}\label{lemmaComputeLink}
    For any cyclic~$\PP$-word~$r$ and any signed measure~$\mu\in\MM^0_s(\sigma)$ finitely supported, so that the support of~$\Theta_\PP(\mu)$ is disjoint from the support of~$\Theta_\PP(\eta_r)$, we have:
    $$\link_\sigma^{ff}(\eta_r-\normR^*\eta_r,\mu)=\int_{\SR}\Sup_c\circ\pi_\PP d\mu$$ 
\end{lemma}

The map~$\Sup_c$ is used to build a 2-chain bounding the algebraic multi-orbit corresponding to~$\eta_r-\normR^*\eta_r$. It comes with a choice for the curves~$a_i$ and~$b_i$. We remove later that choice by replacing~$c$ by another adapted 1-cycle. We denote by~$v_i=a_i(0)=b_{i-1}(1)$,~$w_i=a_i(1)=b_i(0)$.

\begin{figure}
    \begin{center}
        \begin{picture}(120,38)(0,0)
        \put(0,0){\includegraphics[width=120mm]{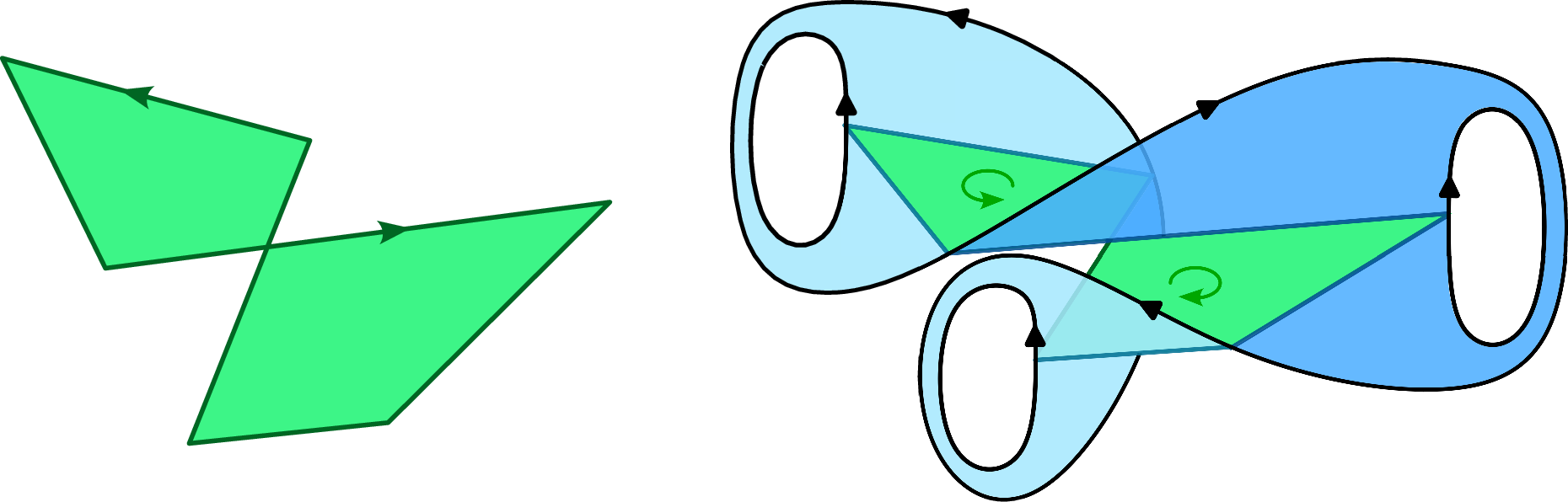}}
        \put(6,25){$\Sup_c=1$}
        \put(20,15){$\Sup_c=-1$}
        \put(12,32){$a_1$}
        \put(2,22){$b_1$}
        \put(30,22.5){$a_2$}
        \put(36.5,10){$b_2$}
        \put(22,3){$a_3$}
        \put(13,10.5){$b_3$}
        \end{picture}
    \end{center}
    \caption{An adapted cell~$c=\sum_i(a_i+b_i)$ and its support~$\Sup_c$ on the left. On the right is represented the corresponding 2-chain in~$M$, bounded the algebraic multi-orbit.}
    \label{figureAdaptedCell}
\end{figure}
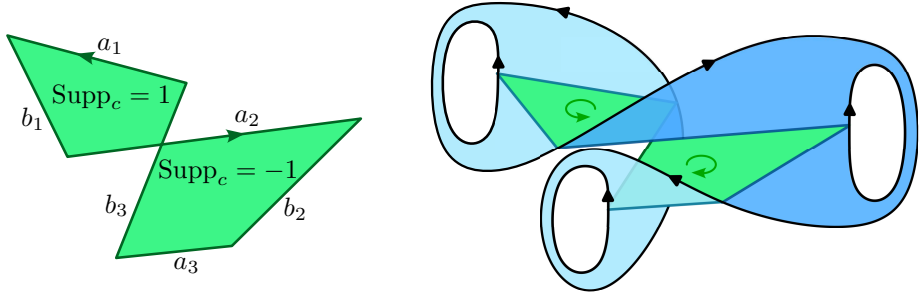

\begin{proof}
    Given a periodic~$\PP$-word~$u$, we denote by~$\gamma_u\subset M$ the algebraic multi-orbit corresponding to the periodic orbit of~$\pi_\PP(\wb u)$. According to 
    Lemma~\ref{lemmaRealizationDegree}, 
    the realization of~$u$ has degree one (here the foliations are orientable). Hence we have~$\Theta_\PP(\eta_u)=\Leb_{\gamma_u}$. We also write~$\Gamma=\gamma_r-\sum_i\gamma_{s_i}$.

    The Lemma is clear when~$r$ does not contain~$\PR$, since~$\normR^*\eta_r=\eta_r$ in that case. Otherwise write~$r=s_1\circ\hdots\circ s_n$ as above. When~$n=1$ we also have~$\normR^*\eta_r=\eta_r$. 
    
    We suppose now~$n\geq 2$. Consider the curves~$a_i$,~$b_i$ and~$c=\sum_i(a_i+b_i)$ as above.    
    Denote by~$C$ an adapted 2-chain in $\PR$, such that~$\partial C=c$. We also denote by~$e_i\subset M$ the 2-chain obtained as the image of the map~$[0,1]\times[0,1]\to M$,~$(s,t)\mapsto\phi_{tT^{s_i}\circ a_i(s)}(a_i(s))$. Its boundary is~$\partial e_i=a_i+\phi_{[0,T^{s_i}(w_i)]}(w_i)+b_i-\phi_{[0,T^{s_i}(v_i)_i]}(v_i)$. 
    Define a 2-chain by~$S=\sum_i e_i - C$, illustrated in Figure~\ref{figureAdaptedCell}.

    \begin{claim}
        We have~$\partial S=\Gamma$.
    \end{claim}

    Take an algebraic multi-orbit~$\Gamma'=\sum_\gamma m(\gamma)\gamma$, whose support is disjoint from the support of~$\Gamma$. By definition of the linking number~$\link_\phi^{ff}$, we have~$\link_\phi^{ff}(\Leb_\Gamma,\Leb_{\Gamma'})=S\algcap \Gamma'$. By construction of~$a_i$,~$e_i$ is disjoint from~$\Gamma'$, so
    \begin{align*}
        \link_\phi^{ff}(\Leb_\Gamma,\Leb_{\Gamma'})
            &= C\algcap\Gamma' \\
            &= \sum_{\gamma}m(\gamma)\sum_{x\in\gamma\cap\PR}C\algcap x \\
            &= \sum_{\gamma}m(\gamma)\sum_{x\in\gamma\cap\PR}\Sup_{c}(x) \\
            &= \int_{\PR} \Sup_c d\Leb_{\Gamma'}^\perp \\
            &= \int_{\SR} \Sup_c\circ\pi_\PP d\mu
    \end{align*}

    The last equality is due to the following. The measures $\Leb_{\Gamma'}^\perp$ and $\pi_\PP^*\mu$ coincide on the interior of $\PR$. On the boundary of $\PR$, the support $\Sup_c$ is equal to zero, maybe outside finitely many points $w_i$ and $v_i$. These points are in the support of $\eta_r$, so by hypothesis they are not in the support of $\mu$. Therefore we have $\int_{\PR} \Sup_c d\pi_\PP^*\mu=\int_{\PR} \Sup_c\Leb_{\Gamma'}^\perp$.

    \begin{claimproof}
        Notice that~$v_i=\pi_\PP(\wb{r_i})$, so we have:
        \begin{align*}
            \phi_{T^{s_i}(v_i)}(v_i)
                &= \phi_{T^{s_i}}(\pi_\PP(\wb{r_i})) \\
                &= \pi_\PP(\sigma^{|s_i|}(\wb{r_i})) \\
                &= \pi_\PP(\wb{r_{i+1}}) \\
                &= a_{i+1}(0)
        \end{align*}
        
        Similarly we have~$\phi_{T^{s_i}(w_i)}(w_i)=w_i$. By construction, we have:
        \begin{align*}
            \partial S &= \sum_i\partial e_i - \partial C\\
                &= \sum_i \big(a_i+b_i+\phi_{[0,T^{s_i}(w_i)]}(w_i)-\phi_{[0,T^{s_i}(v_i)_i]}(v_i)\big) - c \\
                &= \sum_i \phi_{[0,T^{s_i}(v_i)_i]}(v_i) - \sum_i \phi_{[0,T^{s_i}(w_i)]}(w_i) \\
                &= \phi_{[0,T^r(v_1)]}(v_1) - \sum_i \phi_{[0,T^{s_i}(w_i)]}(w_i) \\
                &= \gamma_{r}-\sum_i\gamma_{s_i}
        \end{align*}
    \end{claimproof}
\end{proof}

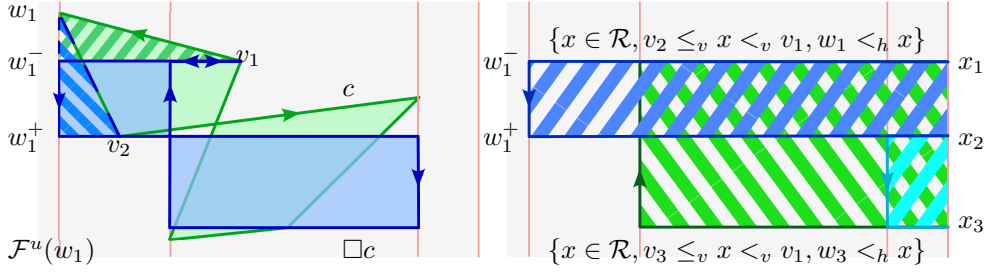
\begin{figure}
    \begin{center}
        \begin{picture}(120,35)(0,0)
        \put(0,0){\includegraphics[width=120mm]{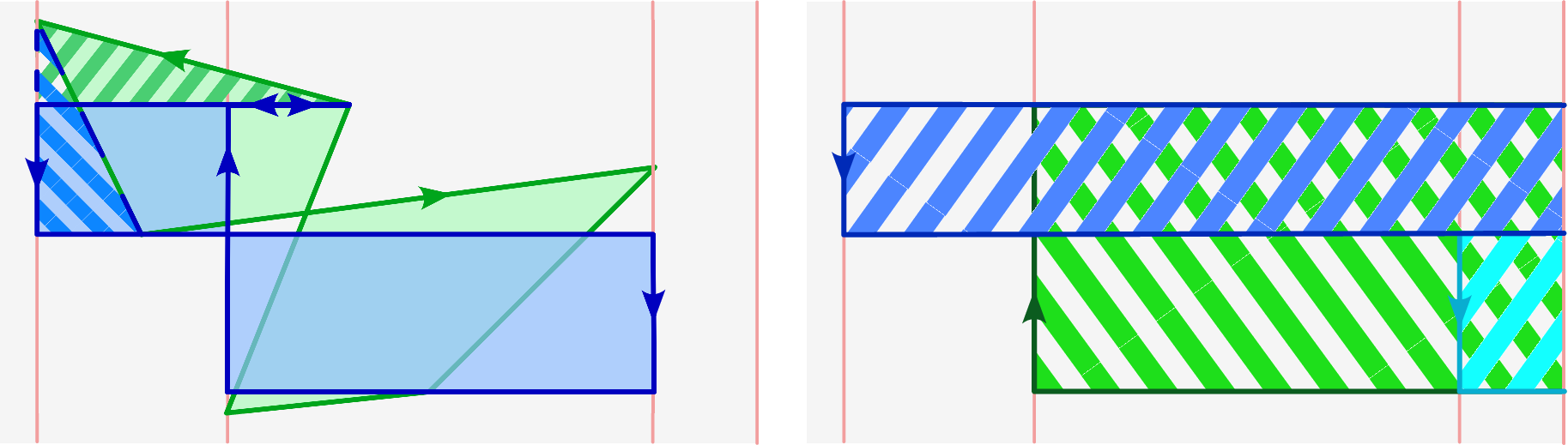}}
        \put(40,0){$\square c$}
        \put(40,21.5){$c$}
        \put(26,26){$v_1$}
        \put(9,14){$v_2$}
        \put(-4,32){$w_1$}
        \put(-4,25){$w_1^-$}
        \put(-4,15){$w_1^+$}
        \put(58.5,25){$w_1^-$}
        \put(58.5,15){$w_1^+$}
        \put(121,25){$x_1$}
        \put(121,15){$x_2$}
        \put(121,4){$x_3$}
        \put(67,28){$\{x\in\PR,v_2\leq_v x<_v v_1, w_1<_h x\}$}
        \put(67,0){$\{x\in\PR,v_3\leq_v x<_v v_1, w_3<_h x\}$}
        \put(-4,0){$\FF^u(w_1)$}
        \end{picture}
    \end{center}
    \caption{Adapted cell~$c=\sum(a_i+b_i)$ and its corresponding square model~$\square c$ on the left. In dashed filling are two sets of equal areas, used in the proof of Lemma~\ref{lemmaIndependentSurfaceModelEqual}. On the right is represented the set of points vertically between~$v_i$ and~$v_{i+1}$, and on the right of~$w_i$. These sets are used to relate the support function of~$\square c$ and the linking function~$G_\PR$, in the proof of Theorem~\ref{theoremLinkingFormProperty}.}
    \label{figureSquareModel}
\end{figure}

We denote by~$w_i^-=W^s_\PR(v_i)\cap W^u_\PR(w_i)$ and~$w_i^+=W^s_\PR(v_{i+1})\cap W^u_\PR(w_i)$. We write $[a,b]_s$ the stable segment in $\PR$ between two points $a,b$ on a same stable leaf, and similarly for $[a,b]_u$. Denote by $\square c$ the adapted 1-cycle given by $\sum_i([v_i,w_i^-]_s+[w_i^-,w_i^+]_u+[w_i^+,v_{i+1}]_s)$, which we call the square model of~$c$ (see Figure~\ref{figureSquareModel}).

\begin{lemma}\label{lemmaIndependentSurfaceModelEqual}
    For any~$\mu\in\MM_s(\sigma)$, we have~$\int_{\SR} (\Sup_c-\Sup_{\square c})\circ\pi_\PP d\mu=0$.
\end{lemma}

\begin{proof}[Proof of Lemma~\ref{lemmaIndependentSurfaceModelEqual}]
    Denote~$d_i=[v_i,w_i^-]_s+[w_i^-,w_i^+]_u+[w_i^+,v_{i+1}]_s-a_i-b_i$, which is a closed adapted 1-chain.
    We denote by~$f_i\colon\Delta_2\to{\PR}$ an adapted 2-simplex whose boundary is~$a_i + [w_i,w_i^-]_s + [w_i^-,v_i]_s$. 
    The image of~$f_i$ lies in the smallest horizontal rectangle of~${\PR}$ containing both~$w_i=\pi_\PP(\wb{s_i})$ and~$v_i=\pi_\PP(\wb{r_i})$. Both $w_i$ and $v_i$ lie in $R_{s_i}$ (the horizontal sub-rectangle of $\PR$ made of the points admitting a short itinerary along $s_i$, see Lemma~\ref{lemmaShortItinerary}).  
    Hence the image of~$f_i$ is included in~$\PR_{s_i}$, so~$\sigma^{s_i}\circ f_i$ is well-defined and~$\partial \sigma^{s_i}\circ f_i=-b_i+[w_i,w_i^+]_u+[w_i^+,v_{i+1}]_s$.

    The 1-chain~$\partial(\sigma^{s_i}\circ f_i-f_i)-d_i$ is equal to~$[w_i,w_i^+]_u-[w_i,w_i^-]_u-[w_i^-,w_i^+]_u$. According to Remark~\ref{lemmaTrivialSupport}, the support of the previous adapted 1-chain is zero. Therefore, we have~$\Sup_{d_i}=\Sup_{\sigma^{s_i}\circ f_i}-\Sup_{f_i}$.

    Let~$\mu\in\MM_s(\sigma)$ be an invariant signed measure. Since~$\sigma^{s_i}\circ f_i$ is obtained from~$f_i$ by pushing along the $\phi$, and the $\phi$ preserves the orientations of the foliations~$W^s_\PR,W^u_\PR$, we have~$\int \Sup_{\sigma^{s_i}\circ f_i}\circ\pi_\PP d\mu = \int \Sup_{f_i}\circ\pi_\PP d\mu$. It follows that~$\int \Sup_{d_i}\circ\pi_\PP d\mu=0$. Notice that we have~$\sum_i d_i=\square c-c$, so the integral~$\int_{\SR} (\Sup_c-\Sup_{\square c})\circ\pi_\PP d\mu$ is equal to zero.
\end{proof}

We now prove that the integral $\int_{\SR^2}\LinkForm d(\nu\otimes\mu)$ computes linking numbers.

\begin{proof}[Proof of Theorem~\ref{theoremLinkingFormProperty}]
    Notice that if~$\nu(\SR)=0$, then~$\normR^*\nu=\nu$ and the two terms in the equality are zero. By linearity, it is enough to prove the result when~$\nu$ is equal to the measure~$\eta_r$ for some primitive cyclic~$\PP$-word~$r$ with~$r_0=\PR$. 
    We use the notations~$s_i$,~$r_i$ used in the above discussion.
    
    Suppose that~$n=1$. Then~$\normR^*\eta_r=\eta_r$ so $\link_\sigma^{ff}(\eta_r-\normR^*\eta_r, \mu)=0$. 
    Fix~$u\in\primR$. If~$u\neq r$, then~$\eta_r\otimes \mu(\Rec_u^\pm)=0$ since no points~$w\in\Su$ are in the support of~$\eta_r$. If~$u= r$, then we also have~$\eta_r\otimes \mu(\Rec_r^\pm)=0$, since if~$w$ is in~$\Su$ and in the support of~$\eta_r$, then~$w=\wb u$ and there is no~$x\in\SR$ satisfying ~$w\leq_v x<_v \sigma^{|u|}(w)=w$ or~$\sigma^{|u|}(w)\leq_v x<_v w$. Hence we have~$\int_{\PW^2}\LinkFormR d\eta_r\otimes\mu=0=\link_\sigma^{ff}(\eta_r-\normR^*\eta_r,\mu)$.

    Now suppose that~$n\geq 2$.
    We denote by~$I^u_+$ the unstable boundary component of~$\PR$, which is on the right, and by~$x_i=W^s(v_i)\cap I^u_+$ the projection of $v_i$ on $I^u_+$, in parallel to $W^s_\PR$. Define the adapted 2-cycle~$f_i$ by~$f_i=[w_i^-,w_i^+]_u+[w_i^+,x_{i+1}]_s+[x_{i+1},x_i]_u+[x_i,w_i^-]_s$.

    \begin{claim}
        We have~$\sum_i\Sup_{f_i}=\Sup_{\square c}$ (illustrated in Figure~\ref{figureSquareModel}).
    \end{claim}

    Recall the definitions:
    $$\Rec_{u}^+ = \{(w,x) \in\Su\times\SR, \wb u<_h x \text{ and } \sigma^{|u|}(w)\leq_v x<_v w \}$$
    $$\Rec_{u}^- = \{(w,x) \in\Su\times\SR, \wb u<_h x \text{ and } w\leq_v x<_v \sigma^{|u|}(w)\}$$
    
    \begin{claim}
       For any $\mu\in\MM_s(\sigma)$, ~$\int_{\SR}\Sup_{f_i}\circ\pi_\PP d\mu=(\caract{\wb{r_i}}\otimes\mu)(\Rec_{s_i}^+)-(\caract{\wb{r_i}}\otimes\mu)(\Rec_{s_i}^-)$
    \end{claim}

    Assume the two claims for now, and take $\mu\in\MM_s(\phi)$ finitely supported, so that the support of $\Theta_\PP(\mu)$ is disjoint from the support of $\Theta_\PP(\eta_r)$. 
    Take~$u\in\primR$. 
    Notice that~$(\caract{\wb{r_i}}\otimes\mu)(\Rec_u^\pm)=0$ when~$u\neq s_i$. It follows that:
    \begin{align*}
        \int_{\SR}\LinkFormR d(\eta_r\otimes\mu) 
            &= \sum_i\int_{\SR}\LinkFormR d(\caract{\wb{r_i}}\otimes\mu) \\
            &= \sum_i\sum_{u\in\primR}\big((\caract{\wb{r_i}}\otimes\mu)(\Rec_u^+)-(\caract{\wb{r_i}}\otimes\mu)(\Rec_u^-)\big) \\
            &= \sum_i\big((\caract{\wb{r_i}}\otimes\mu)(\Rec_{s_i}^+)-(\caract{\wb{r_i}}\otimes\mu)(\Rec_{s_i}^-)\big) \\
            &= \sum_i \int_{\SR}\Sup_{f_i}\circ\pi_\PP d\mu \qquad \text{ according to Claim 2.}\\
            &= \int_{\SR}\Sup_{\square c}\circ\pi_\PP d\mu \qquad \text{ according to Claim 1.}\\
            &= \link_\sigma^{ff}(\eta_r-\normR^*\eta_r,\mu) \text{ according to Lemmas~\ref{lemmaComputeLink},~\ref{lemmaIndependentSurfaceModelEqual}}
    \end{align*}

    Which implies~$\link_\sigma^{ff}(\eta_r-\normR^*\eta_r,\mu)= \int_{\SR}\LinkFormR d(\eta_r\otimes\mu)$.

    \begin{claimproof}
        One can verify that the closed 1-chain~$\sum_i f_i-\square c$ is a chain contained in the union of the unstable leaf~$I^u_+$ and of finitely many stable leaves containing the points~$v_i$. Essentially the parts in $[w_i^-,w_i^+]_u$ cancel each other. Hence its support function is~$\Sup_{(\sum_i f_i-\square c)}=0$ (see Remark \ref{lemmaTrivialSupport}).
    \end{claimproof} 

    \begin{claimproof}
        Notice that for~$(w,x)\in \Rec_{s_i}^+$ to be in the support of~$(\caract{\wb{r_i}}\otimes\mu)$, we must have~$w=\wb{r_i}$ and~$\wb{r_{i+1}}<_v\wb{r_{i}}$. Similarly for~$(w,x)\in \Rec_{s_i}^-$ to be in the support of the same measure, we must have~$\wb{r_i}<_v\wb{r_{i+1}}$. Hence the mass of either~$\Rec_{s_i}^+$ or~$\Rec_{s_i}^-$ is zero. 
        
        Consider the case when~$\wb{r_{i+1}}<_v\wb{r_i}$.
        By construction, we have~$w_i^+<_vw_i^-$. So~$f_i$ is the boundary of a Markov rectangle~$R_i\subset\PR$ with the induced orientation. Hence the support of the support of~$f_i$ is equal to:
        $$\Sup_{f_i}(x) =
        \begin{array}{|ll}
            1 & \text{if~$x$ is in the interior of the rectangle~$R_i$}\\
            1 & \text{ if } x\in (w_i^+,x_{i+1}]_s\cup [x_{i+1},x_i)_u \\
            0 & \text{otherwise}
        \end{array}$$
        So the integral~$\int_{\SR} \Sup_{f_i}\circ\pi_\PP d\mu$ is equal to the mass for~$\mu$ of the subset of~$\SR$ of points $x$ satisfying the two following conditions
        \begin{align*}
            \pi_\PP(\wb{s_i}) = w_i &<_h \pi_\PP(x) \\
            \pi_\PP(\wb{r_{i+1}})=v_{i+1}&\leq_v \pi_\PP(x)<_vv_i=\pi_\PP(\wb{r_i})
        \end{align*}

        According to 
        Lemma~\ref{lemmaDistinctRealization} 
        this condition on~$x$ is equivalent to having 
        \begin{align*}
            \wb{s_i} &<_hx \\
            \sigma^{|s_i|}(\wb{r_{i+1}})&\leq_vx<_v\wb{r_i}            
        \end{align*} Hence we have~$\int_{\SR}\Sup_{f_i}\circ\pi_\PP d\mu=\mu(\Rec_{s_i}^+)$. In the case~$\wb{r_{i+1}}>_v\wb{r_i}$, the same argument shows that~$\int_{\SR}\Sup_{f_i}\circ\pi_\PP d\mu=-\mu(\Rec_{s_i}^-)$.  
    \end{claimproof}
\end{proof}

\subsection{Linking number~$\link_\phi$}\label{subsection-LK-continuous}

We compose the reduction maps for all Markov rectangles of~$\PP$ to reduce an invariant signed measure on one supported by a fixed finite set. We then prove Theorem~\ref{theoremMeasureLinkingNumber}. Denote by~$(\PR_1,\cdots,\PR_p)$ the Markov rectangles in~$\PP$, by~$\normB{k}=\norm{k}\circ\hdots\circ\norm{1}$ and by $\normalize=\normB{p}$. 

\begin{definition}
    For~$\nu,\mu\in\MM_s(\sigma)$ such that~$\cohom{\nu}=\cohom{\mu}=0$, we define the \emph{linking number} between~$\nu$ and~$\mu$ by  $$\link_\sigma(\nu,\mu)=\link_\phi^{f}(\Theta_\PP(\normalize^*\nu), \Theta_\PP(\mu))+\sum_{k=1}^{p-1}\int_{\PW^2}\LinkFormi{k} d(\normB{k}^*\nu)\otimes\mu$$
\end{definition}

\begin{lemma}\label{lemmaContinuityLinkSigma}
    The map~$\link_\sigma$ is continuous on~$\MM^0_s(\sigma)\times\MM^0_s(\sigma)$.
\end{lemma}

\begin{proof}
    The maps~$\nu\mapsto \normalize^*\mu$ and~$\Theta_\PP$ are continuous (see Proposition~\ref{lemmafR*continuous}). The signed measure~$\normalize^*\nu$ is supported by a subset of~$\prim$ (see Lemma~\ref{lemmaReductionSupport}), which is finite. Denote by~$\Delta\subset M$ the union of the periodic orbits containing~$\pi_\PP(\prim)$. The map~$\link_\phi^f$ is continuous on the set of pairs of signed measures~$(\nu_\phi,\mu_\phi)\in\MM^0_s(\phi)\times \MM^0_s(\phi)$ for which~$\Delta$ contains the support of~$\nu_\phi$ (see Proposition \ref{propLinkOrbitSupported}). Therefore the term~$\link_\phi^{f}(\Theta_\PP(\normalize^*\nu), \Theta_\PP(\mu))$ is continuous. The second term is continuous according to Proposition~\ref{propContinuityIntLKForm}.
\end{proof}

\begin{lemma}\label{lemmaLinkExtension}
    For any~$\nu,\mu\in\MM^0_s(\sigma)$ with finite supports, so that $\Theta_\PP(\mu)$ and $\Theta_\PP(\nu)$ have disjoint support, we have  $$\link_\sigma(\nu,\mu)=\link_\sigma^{ff}(\nu,\mu)$$
\end{lemma}

\begin{proof}
    Suppose first that the support of~$\mu$ is also disjoint from~$\prim$. Then for all~$k$, the support of~$\Theta_\PP(\normB{k}^*\nu)$ is disjoint from the support of~$\Theta_\PP(\mu)$. Hence according to Theorem~\ref{theoremLinkingFormProperty}, we have~$\int_{\PW^2}\LinkFormi{k} d(\normB{k}^*\nu)\otimes\mu =\link_\sigma^{ff}(\normB{k}^*\nu-\normB{k+1}^*\nu,\mu)$. It follows that:
    \begin{align*}
        \link_\sigma^{ff}(\nu,\mu)
            &= \link_\sigma^{ff}(\normalize^*\nu,\mu)+\sum_{k=1}^{p-1}\link_\sigma^{ff}(\normB{k}^*\nu-\normB{k+1}^*\nu,\mu)\\
            &=\link_\phi^{f}(\Theta_\PP(\normalize^*\nu), \Theta_\PP(\mu))
            + \sum_{k=1}^{p-1}\int_{\PW^2}\LinkFormi{k} d(\normB{k}^*\nu)\otimes\mu \\
            &= \link_\sigma(\nu,\mu)
    \end{align*}

    The general case follows from the density of signed measures in the first case (see Lemma~\ref{lemmaDensityFiniteMeasureDisjoint}) and from the continuity of $\link_\sigma$.
\end{proof}

\begin{lemma}\label{lemmaFactoriseLink}
    There exists a continuous bilinear map~$\link_\phi\colon\MM^0_s(\phi)\times \MM^0_s(\phi)\to\RR$ such that for all signed measures~$\nu,\mu\in\MM^0_s(\sigma)$, one has~$\link_\sigma(\nu,\mu)=\link_\phi(\Theta_{\PP}\nu,\Theta_{\PP}\mu)$. Additionally~$\link_\phi$ is continuous on~$\MM^0_s(\phi)\times \MM^0_s(\phi)$.
\end{lemma}

\begin{proof}
    It is enough to prove that for any signed measures~$\nu_1,\nu_2,\mu\in\MM^0_s(\sigma)$ with~$\Theta_\PP(\nu_1)=\Theta_\PP(\nu_2)$, we have~$\link_\sigma(\nu_1,\mu)=\link_\sigma(\nu_2,\mu)$. Once it is proven, we define the map~$\link_\phi$ by~$\link_\phi(\nu_\phi,\mu_\phi)=\link_\sigma(\nu_\sigma,\mu_\sigma)$ for any signed measures~$\nu_\sigma,\mu_\sigma\in\MM^0_s(\sigma)$ in the preimages of~$\nu_\phi,\mu_\phi$ by the surjection~$\Theta_{\PP}\colon \MM_s(\sigma)\to\MM_s(\phi)$. 

    To prove the claim, take~$\nu_1,\nu_2,\mu\in\MM^0_s(\sigma)$ with~$\Theta_\PP(\nu_1)=\Theta_\PP(\nu_2)$. Denote by~$\Gamma$ the set of periodic orbits of~$\phi$, contained in the stable/unstable boundary of~$\PP$. The map~$\Theta_\PP$ is injective outside the preimage of the set~$\cup_{\gamma\in\Gamma}(\FF^s(\gamma)\cup\FF^u(\gamma))$ for the map~$\pi_\PP$. Hence the signed measures~$\nu_i$ coincide on this set. 
    
    According to Lemma~\ref{lemmaLeafZeroMeasure}, the signed measure~$\Theta(\nu_i)$ charges the set~$\pi^{-1}(\cup_{\gamma\in\Gamma}(\FF^s(\gamma)\cup\FF^u(\gamma)))$ only on~$\Gamma$. 
    It follows that~$\nu_1-\nu_2$ has a finite support. Suppose that~$\mu$ is finitely supported, and that the support of~$\Theta_\PP(\mu)$ is disjoint from the support of~$\Theta(\nu_1-\nu_2)$. According to Lemma~\ref{lemmaLinkExtension}, we have~$\link_\sigma(\nu_1-\nu_2,\mu)=\link_\sigma^{ff}(\nu_1-\nu_2,\mu)=\link_\phi^f(\Theta_\PP(\nu_1-\nu_2),\Theta_\PP(\mu))=0$ since ~$\Theta_\PP(\nu_1-\nu_2)=0$. By continuity of $\link_\sigma$ and density of the signed measures satisfying that assumption, we have~$\link_\sigma(\nu_1,\mu)=\link_\sigma(\nu_2,\mu)$ for all signed measures $\mu\in\MM^0_s(\sigma)$. 
\end{proof}

\begin{proof}[Proof of Theorem~\ref{theoremMeasureLinkingNumber}]
    Lemma~\ref{lemmaDensityFiniteMeasureDisjoint} implies that we can approximate any pair of invariant signed measures~$(\nu,\mu)\in\MM^0_s(\phi)$ by sequence of pairs~$(\nu_n,\mu_n)\in\MM^0_s(\phi)$, where~$\nu_n$ and~$\mu_n$ are supported by disjoint and finite unions of periodic orbits. Therefore the uniqueness and symmetry of the map~$\link_\phi$ follow from the continuity. 

    When the stable and unstable foliations of~$\phi$ are orientable, the theorem follows from the above discussion (see Lemmas~\ref{lemmaLinkExtension} and~\ref{lemmaFactoriseLink}). 

    Suppose that the stable and unstable foliations are not orientable. We denote by~$\wh\pi\colon\wh M\to M$ the bundle-orientations covering (which is of degree 2), and by~$g$ the non-trivial deck transformation of~$\wh M\to\wh M$. Lift~$\phi$ to a flow~$\psi$ on~$\wh M$. 
    We denote by~$\link_{\psi}\colon\MM^0_s(\psi)^2\to\RR$ the linking number map obtained in the orientable case. 
    Take a signed measure~$\mu\in\MM^0_s(\phi)$ and denote by~$\mu\circ\wh\pi\in\MM_s(\psi)$ the signed measure which coincides with~$((\wh\pi_{|U})^{-1})^*\mu$ on each set~$U\subset\wh M$ where~$\wh\pi$ is invertible. Notice that~$\mu\circ\wh\pi$ is the unique~$\psi$ and~$g$ invariant signed measure lifting~$2\mu$.
    
    Take a closed 1-form~$\wh\alpha$ on~$\wh M$. Since~$\wh\alpha+g^*\wh\alpha$ is invariant by the action of~$g$, there exists a closed 1-form~$\alpha$ on~$M$ satisfying~$\wh\pi^*\alpha=\wh\alpha+g^*\wh\alpha$. The flow~$\psi_t$ and the signed measure~$\mu\circ\wh\pi$ are invariant by~$g$, so we have 
    \begin{align*}
        \int_{\wh M}\ivectD{\psi}\wh\alpha d\mu\circ\wh\pi 
            &= \frac{1}{2}\bigg(\int_{\wh M}\ivectD{\psi}\wh\alpha d\mu\circ\wh\pi + \int_{\wh M}\ivectD{\psi}g^*\wh\alpha d\mu\circ\wh\pi \bigg) \\
            &= \frac{1}{2}\int_{\wh M}\ivectD{\psi}\wh\pi^*\alpha d\mu\circ\wh\pi \\
            &= \int_M\ivectD{\phi}\alpha d\mu = 0  
    \end{align*}
    It follows that~$\mu\circ\wh\pi$ is null-homologous. Take~$K\subset\wh M$ a measurable subset such that~$(K,g\cdot K)$ is a partition of~$\wh M$. Then the projection~$K\to M$ is bijective. Given a continuous map~$f\colon\wh M\to\RR$, one has 
    \begin{align*}
        \int_{\wh M}fd\mu\circ\wh\pi 
            &=\int_K fd\mu\circ\wh\pi + \int_{g\cdot K}fd\mu\circ\wh\pi \\
            &=\int_M f\circ (\wh\pi_K)^{-1}d\mu + \int_M f\circ (g\cdot \wh\pi_K)^{-1}d\mu \\
            &= \int_M\Bigg(\sum_{y\in\wh\pi^{-1}(x)}f(x)\Bigg)d\mu(x)
    \end{align*}

    It follows that~$\mu\mapsto\mu\circ\wh\pi$ is continuous.
    For~$\nu,\mu\in\MM^0_s(\psi)$, we define~$\link_\phi(\nu,\mu)=\frac{1}{2}\link_{\psi}(\mu\circ\wh\pi,\nu\circ\wh\pi)$. By construction, it is continuous.
    Take~$\Gamma_1,\Gamma_2$ two algebraic multi-orbits for~$\phi$. Take~$S_1$ a 2-chain on~$M$, bounding~$\Gamma_1$. Lift~$\Gamma_i$
    to an algebraic multi-orbit~$\wh\Gamma_i$ of~$\psi$, which satisfies~$g\cdot\wh\Gamma_i=\wh\Gamma_i$. Also lift~$S_1$ to a 2-chain~$\wh S_1$ in~$\wh M$, bounding~$\wh\Gamma_1$. According to the orientable case, we have 
    \begin{align*}
        \link_\phi(\Leb_{\Gamma_1},\Leb_{\Gamma_2}) 
            &= \frac{1}{2}\link_{\psi}(\Leb_{\Gamma_1}\circ\wh\pi,\Leb_{\Gamma_2}\circ\wh\pi) \\
            &= \frac{1}{2}\link_{\psi}(\Leb_{\wh\Gamma_1},\Leb_{\wh\Gamma_2}) \\
            &= \frac{1}{2}\wh S_1\algcap\wh\Gamma_2 \\
            &= S_1\algcap\Gamma_2 \\
            &= \link_\phi^{ff}(\Leb_{\Gamma_1},\Leb_{\Gamma_2})
    \end{align*}
\end{proof}

We need to prove a last lemma. 

\begin{proof}[Proof of Lemma \ref{lemmaLinkSmoothMeasure}]
    Let $\alpha$ and $V$ be respectively smooth 1-form and 3-form, so that $V$ is invariant by $\phi$ and $d\alpha=\iota_{X}V$. Denote by~$\nu_V\in\MM^0_s(\phi)$ the signed measure induced by $V$. We prove that $\link_\phi(\nu_V,\mu)=\int_M\iota_{X}\alpha d\mu$ for any $\mu\in\MM^0_p(\phi)$.

    It is enough to verify the equality for signed measures of the form $\Leb_\Gamma$ for null-homologous algebraic multi-orbits $\Gamma$. The lemma follows from continuity of $\link_\phi$ and from the density 
    Lemma~\ref{lemmaFinitBoundedMeasureDense}.
    Denote by $\Delta$ the support of $\Gamma$ and take a rational 2-chain $S$ in $M_\Gamma$ such that $(\pi_\Delta)_*(\partial S)=\Gamma$. We can take a closed 1-form $\omega$ on $M_\Gamma$ so that $P_\Delta([\omega])=[S]$.
    From the definition $\link_\Gamma^\Delta$ follows
    \begin{align*}
        \link_\phi^f(\nu_\alpha,\Leb_\Gamma) 
                &= \int_{M_\Delta} (\iota_{X_\Delta}\omega) \pi_\Delta^*V 
                = \int_{M_\Delta} \omega\wedge (\iota_{X_\Delta} \pi_\Delta^*V) \\
                &= \int_{M_\Delta} \omega\wedge d\pi_\Delta^*\alpha \\
                &= \int_S d\pi_\Delta^*\alpha \qquad \text{by Poincaré duality of $[\omega]$ and $[S]$} \\
                &= \int_{\partial S}\pi_\Delta^*\alpha
                = \int_{\Gamma}\alpha \\
                &= \int_M \iota_{X}\alpha d\Leb_\Gamma
    \end{align*}
\end{proof}

\appendix

\section{Birkhoff section with given boundary}\label{sectionBSexistence}

We prove 
Theorem~\ref{theoremExistenceBSwithLink}, 
which states that given~$\Gamma$, there exists a Birkhoff section bounding~$N\Gamma$ for some~$N\in\NN$, if and only if $\Leb_\Gamma$ is Reeb-like. As a corollary, we obtain a second proof of Barbot 
Theorem~\ref{theoremBarbotTheorem}.

Given an oriented compact 3-manifold~$N$ and a flow~$\psi$ on~$N$, a \emph{global section} of~$\psi$ is an oriented compact surface~$S\subset N$ which is positively transverse to the flow, with~$\partial S\subset\partial N$ and which intersects all orbits of~$\psi$. It follows from compactness that~$S$ intersects every orbit arc of length~$T$ for some~$T>0$.

We start by constructing a Birkhoff section from a global section on a blown-up manifold. Take~$\Delta\subset M$ a finite union of periodic orbits. Denote by~$M_\Delta$ the blown-up manifold of~$M$ along~$\Delta$, by~$\pi_\Delta\colon M_\Delta\to M$ the blow-down projection, and by~$\phi^\Delta$ the lifted flow on~$M_\Delta$.

\begin{lemma}\label{lemmaGlobalToBirkhoffSection}
    Let~$S\subset M_\Delta$ a smooth global section of~$\phi^\Delta$. Then there exists a smooth function~$T\colon S\to\RR$ satisfying that the image~$\pi_\Delta\circ\phi_T(S)$ is a Birkhoff section of~$\phi$.
\end{lemma}

\begin{proof}
    In this proof, all isotopies are isotopies along the flow.
    There are two types of boundary components for $S$. The curves in $\partial S\cap\TT_\gamma$, for some $\gamma\subset\Delta$, which are homologous to a fiber of $\TT_\gamma\to\gamma$, and the others. We call these boundary components of type 1 and type 2.
    We need to prove that up to an isotopy on $S$, the projection $\pi_\Delta$ restricts to a nice map on~$S$. Here 'nice map' means that for each boundary component $\delta\subset\TT_\gamma$ of type 1, $\delta$ is sent to a single point in $\gamma$ and the image of $S$ is smooth close to that point. And for each boundary component $\delta$ of type 2, $\pi_\Delta$ restricts to an immersion close to $\delta$. Then the image of $S$ is a Birkhoff section.

    Take a boundary component $\delta\subset\partial S\cap\TT_\gamma$ of type 1. Then $\delta$ is isotopic to one fiber 
    (see~\cite[Theorem C]{FriedCrossSection}). 
    In fact there exists a global isotopy of $\TT_\gamma$ sending $\partial S\cap\TT_\gamma$ to a finite union of fibers of $\TT_\gamma\to\gamma$. We can extend the global isotopy on $\TT_\gamma$ to a global isotopy on $M_\Delta$. Denote by $S_1$ the image of the surface $S$ by the global isotopy along the flow. We can take the global isotopy so that $\pi_\Delta(S_1)$ is an embedded and smooth surface, close to the orbits of type 1.

    Take a boundary component $\delta\subset\partial S_1\cap\TT_\gamma$ of type 2. Then similarly $\delta$ is isotopic to a smooth curve $\delta'$ in $\TT_\gamma$, transverse to the flow, and transverse to the foliation on $\TT_\gamma$ by fibers of $\TT_\gamma\to\gamma$. Similarly there exists a global isotopy on $M_\Delta$, whose image of $S_1$ is a surface $S_2$, and which satisfies that $\pi_\Delta(S_2)$ is immersed close to the orbits of type 2. Then $\pi_\Delta(S_2)$ is a Birkhoff section.
\end{proof}

\begin{proof}[Proof of Theorem~\ref{theoremExistenceBSwithLink}]
    The implication $1\Rightarrow 2$ is already proven in Section \ref{sectionBS} (see Lemma \ref{lemmaBStoReebLike}). We prove $2\Rightarrow 1$.

    Let~$\Gamma$ be an algebraic multi-orbit with integer coefficients, so that~$\Leb_\Gamma$ is Reeb-like. Denote by~$\Delta$ the support of~$\Gamma$. 
    Fix an element~$\omega_1\in H^1(M_\Delta,\RR)$ which satisfies~$(\pi_\Delta)_*(\partial P_\Delta(\omega_1))=[\Gamma]$. Since $\Gamma$ has integer coefficients, we can choose~$\omega_1\in H^1(M_\Delta,\QQ)$ with rational coefficients.
    Define the set 
    $$F=\{(\cohom{\pi_\Delta^*\mu},\cohomG{\mu}\cdot\omega_1)\in H_1(M,\RR)\times\RR,\mu\in\MM_p(\phi^\Delta)\}$$
    By assumption, we have~$\cohomG{\mu}\cdot\omega_1=\link_\phi(\Leb_\Gamma,\pi_\Delta^*\mu)>0$ for all probability measure~$\mu\in\MM^0_p(\phi^\Delta)$. Therefore the intersection~$F\cap(0\times(-\infty,0])$ is empty. According to 
    Lemma~\ref{lemmaNullHomologousIntersectionIsEnough},
    there exists~$\eta>0$ and a linear map~$f\colon H_1(M,\RR)\to\RR$ satisfying that for all~$\mu\in\MM_p(\phi^\Delta)$, we have~$\cohomG{\mu}\cdot\omega_1\geq f(\cohom{\pi_\Delta^*\mu})+\eta$. 

    By Poincaré duality,~$f$ corresponds to an element~$\omega_2\in H^1(M,\RR)$. Since~$\MM_p(\phi^\Delta)$ is compact, we can perturbate~$f$  so that~$\omega_2$ has rational coefficients. Take an integer~$N\in\NN_{>0}$ satisfying that~$N(\omega_1-\pi_\Delta^*\omega_2)$ have integer coefficients. Then for all probability measures~$\mu\in\MM_p(\phi^\Delta)$, we have  $$\cohomG{\mu}\cdot N(\omega_1-\pi_\Delta^*\omega_2)\geq N\eta>0$$
    
    It follows from Schwartzman-Fried-Sullivan Theory 
    (see~\cite[Theorem D]{FriedCrossSection}), 
    that~$\phi^\Delta$ admits a global section~$S$ whose relative homology class~$[S]\in H_2(M_\Delta,\partial M_\Delta,\RR)$ is dual to~$N(\omega_1-\pi_\Delta^*\omega_2)$. According to
    Lemma~\ref{lemmaGlobalToBirkhoffSection}, 
    the global section can be taken so that~$\pi_\Delta(S)$ is a Birkhoff section of~$\phi$. Additionally the boundary of~$\pi_\Delta(S)$ is homologous to~$(\pi_\Delta)_*(\partial P_\Delta(N(\omega_1-\pi_\Delta^*\omega_2)))= N[\Gamma]$ in~$H_1(\Delta,\RR)$. Therefore~$N\Gamma$ bounds a Birkhoff section. 
\end{proof}

We should discuss a particular case. When there is no null-homologous invariant probability measure, 
Sullivan~\cite{Sullivan76} 
proved that the flow admits a global section. Since the flow is Anosov, the first-return map for that global section is conjugated to a linear Anosov diffeomorphism on the torus. In that case, every null-homologous algebraic multi-orbit, with integer coefficients, bounds a Birkhoff section 
(see~\cite[Section 4.2]{Marty20}). 
Hence Theorem~\ref{theoremExistenceBSwithLink} 
is a generalization of the suspension case, where we have less control over the multiplicity at the boundary.

\begin{proof}[Proof of the remaining equivalences in the Table~\ref{mainCorollaryEquivalence}]
    Recall the assertions in the table we prove to be equivalent: there exists a positive Birkhoff section (cell~2) if and only if there exists a Reeb-like invariant probability measure (cell~3) if and only if the flow is orbit equivalent to a Reeb-Anosov flow (cell~4). 
    Here we prove the implications in the first column of the table: cell~4~$\implies$~cell~3~$\iff$~cell~2. The implication cell~4~$\implies$~cell~2 together with the implication cell~2~$\implies$~cell~1 (the Anosov flow is positively skewed) stated in Theorem \ref{theoremABM} yield a second proof to Barbot's Theorem~\ref{theoremBarbotTheorem}: Reeb-Anosov are skewed.

    Suppose that there exists $\psi$ orbit equivalent to $\phi$, which preserves a contact form $\alpha$ with $\alpha\wedge d\alpha>0$ (cell 4). Then according to 
    Theorem~\ref{theoremReebLikeCondition}, 
    the $\psi$-invariant probability measure induced by $\alpha\wedge d\alpha$ has the Reeb-like property. The Reeb-like property is invariant under orbit equivalence, so there exists a $\phi$-invariant probability measure with the Reeb-like property. Namely it is null-homologous and has only positive linking number with null-homologous $\phi$-invariant probability measures (cell 3).
    
    The implication cell 2 $\implies$ cell 3 follows from the above discussion 
    (Lemma~\ref{lemmaBStoReebLike}). 
    
    Suppose that there exists a measure $\mu\in\MM^0_p(\phi)$ with the Reeb-like property (cell 3). According to 
    Lemma~\ref{lemmaFinitBoundedMeasureDense}, 
    there exists a sequence of null-homologous invariant 1-cycles $\Gamma_n$, with rational coefficients, so that $\Leb_{\Gamma_n}$ converges toward $\mu$. By continuity of the linking number, for large $n$, $\Leb_{\Gamma_n}$ has positive linking numbers with all measures in $\MM^0_p(\phi)$. It follows from 
    Theorem~\ref{theoremExistenceBSwithLink}
    that there exists a positive Birkhoff section (cell 2).
\end{proof}

\section{Fried-desingularisation} \label{appendixFS}

We prove 
Proposition \ref{propFriedDesingularization}: 
given a Birkhoff section $S_1$ and a partial section $S_2$, there exists a Birkhoff section homologous to $[S_1]+[S_2]$ in $H_2(M,\partial S_1\cup\partial S_2,\RR)$.
Take~$\Delta=\partial S_1\cup\partial S_2$, $M_\Delta$ the blown-up manifold and fix a Riemannian metric~$m$ on~$M_\Delta$.
We denote by~$D$ the diameter of~$M_\Delta$ for~$m$, and by~$\lVert \cdot \rVert$ one norm on~$H_2(M_\Delta,\partial M_\Delta,\RR)$. Define the normalization map $$N\colon H_2(M_\Delta,\partial M_\Delta,\RR)\to H_2(M_\Delta,\partial M_\Delta,\RR)$$
given by~$N(0)=0$ and~$N(x)=x/\lVert x\rVert $ when~$x\neq 0$.

For~$x\in M_\Delta$ and~$T\geq 0$, denote by~$\gamma_x^T\colon [0,T+1]\to M_\Delta$ one map given by~$\gamma_x^T(t)=\phi^\Delta_t(x)$, for~$t\in [0,T]$, and~$(\gamma_x^T)_{|[T,T+1]}$ is smooth path from~$\phi_T(x)$ to~$x$, whose length is at most~$D$. The choice of the curve~$(\gamma_x^T)_{|[T,T+1]}$ has no influence on the following notion. The accumulating points of elements of the form~$N([\gamma_x^T])$, for~$T\xrightarrow[]{}+\infty$, are called \emph{asymptotic directions}. When the ambient manifold is closed,
Schwartzman and Fried~\cite[Theorem D]{FriedCrossSection} 
proved that a compact surface is relatively homologous to a global section of flow if and only if it intersects algebraically positively every asymptotic direction. Fried's proof works with minor modification for compact manifolds with boundary (tangent to the flow lines). 

Hryniewicz \cite{Hryniewicz2020} proved the existence of a Birkhoff section in a given homology class (relatively to a given set of boundary periodic orbits) under the same condition (intersecting positively all asymptotic directions in the blown-up manifold). Note that Hryniewicz's definition of Birkhoff sections differs from ours in the general case, but they are equivalent for Anosov flows. He requires that a Birkhoff section intersects every orbit ray (in the future and in the past) instead of intersecting every long enough orbit arc. But in the hyperbolic case, he proves that it is equivalence (see \cite[Equation 30]{Hryniewicz2020}).

\begin{proof}[Proof of Proposition \ref{propFriedDesingularization}]
We lift the surface~$S_i$ to a surface~$\wh{S_i}$ inside~$M_\Delta$. We write~$\wh{S}=\wh{S_1}\cup \wh{S_2}$. Since~$S_1$ is a Birkhoff section, there exists~$T_0>0$ for which every orbit arc of length~$T_0$ of~$\phi^\Delta$ intersects~$\wh{S_1}$. Recall that~$S_i$ is positively transverse to the~$\phi$.
Denote by~$C_1>0$ the supremum of the algebraic intersection of~$\wh S$ and a smooth curve of length~$D$, which is finite. Take~$C_2>0$ for which every closed~$\Class^1$ curve~$\gamma\subset M_\Delta$ satisfies~$\lVert [\gamma]\rVert \leq C_2 length_m(\gamma)$. Also take~$C_3>0$ the maximum of~$\Big\lVert\frac{\partial \phi^\Delta_t}{\partial t}\Big\rVert_m$, so that the length for $m$ of an orbit arc $\phi_{[0,T]}(x)$ is at most $C_3T$.

Take~$x\in M_\Delta$ and~$T>0$. When~$x$ is not in the boundary of $M_\Delta$, the algebraic intersection between~$\gamma_x^T$ and~$\wh{S}$ is at least 
$$\phi^\Delta_{[0,T]}(x)\algcap ([\wh{S_1}]+[\wh{S_2}])-C_1\geq \frac{T}{T_0}-1-C_1$$
By continuity, it also holds for~$x\in\partial M_\Delta$. So we have 
\begin{align*}
    N\big([\gamma_x^T]\big) \algcap \wh{S}
            &\geq \frac{\gamma_x^T\algcap \wh{S}}{\big\lVert [\gamma_x^T]\big\rVert} \\
            &\geq \frac{\frac{T}{T_0}-1-C_1}{C_2 (C_3T+D)} \\
    \liminf_{T\xrightarrow[]{}+\infty} N\big([\gamma_x^T]\big) \algcap \wh{S} 
            &\geq \frac{1}{T_0C_2C_3}
\end{align*}

Hence~$\wh{S}$ intersects positively every asymptotic direction. Fried proved in similar settings that there exists a global section relatively homologous to~$\wh{S}$ in~$H_2(M_\Delta,\partial M_\Delta,\RR)$. More precisely when the  manifold playing the role of $M_\Delta$ is closed (Theorem D in Fried's article) or when the flow is transverse to the boundary of $M_\Delta$: in Fried's Theorem E, global section to the flow restricted to an invariant set are characterized in a similar manner. One can check that Fried's proof holds for compact manifold with boundary tangent to the flow. Alternatively one can reprove it using a doubling manifold (of $M_\Delta$) and applying the closed case. Or one can add a thick torus to all boundary components of $M_\Delta$, extend the flow to be transverse to the boundary, and apply Theorem E to the complement of the thick tori.

The existence of a Birkhoff section homologous to $[S_1]+[S_2]$ follows from 
Lemma \ref{lemmaGlobalToBirkhoffSection}.
\end{proof}

\section{Hart lemma for Hölder regularity} \label{appendixHart}

We give a proof of the following adaptation of Hart's lemma. For a flow $\Psi$ and a vector field $Z$, we denote by~$f^*\Psi_t=f\circ\Psi_t\circ f^{-1}$ the flow conjugated to~$\Psi$ by~$f$ and by $f^*Z=df(Z\circ f^{-1})$. In particular we have~$\frac{\partial f^*\Psi_t}{\partial t}=f^*\frac{\partial \Psi_t}{\partial t}$.

\begin{lemma}[Extended Hart's lemma~\cite{Hart83}]\label{lemmaHartManifold}
    Let $k\geq 1$ and~$\Psi$ be a~$\Cxa{k}$ flow on a smooth compact manifold~$M$, and suppose that $\vectD{\Psi}$ is $\Cxa{k}$ on a neighborhood of a compact set $K\subset M$. Then there exists a~$\Cxa{k}$ diffeomorphism~$f\colon M\to M$ such that~$f_{|K}=\id_K$ and $\vectD{f^*\Psi}$ is of class $\Cxa{k}$.
\end{lemma}

The key argument is, in a local chart, to integrate along the flow to regularize its generating vector field.
Take the space~$E=\RR^m$ for some~$m>0$ and~$\Psi$ a flow of class~$\Cxa{k}$ on~$E$. Let $\Psi$ be a $\Cxa k$ flow on $E$ and denote by~$Z$ the generating vector field of~$\Psi$.

Fix three open subsets~$U,V,W\subset E$ such that~$\wb W\subset V$,~$\wb V\subset U$,~$\wb V\cap \partial U=\emptyset$ and~$\wb U$ is compact. Let~$T>0$ be such that~$\Phi_T(V)\subset U$, and~$c\colon U\to[0,T]$ be a~$\Class^{2+\alpha}$ function such that~$c$ is constant on~$W$,~$c>0$ on~$V$ and~$c\equiv 0$ outside~$V$. Define the map~$F_c\colon U\to\RR^3$ by $$F_c(x)=\int_0^1\Psi_{c(x)s}(x)ds$$
\begin{lemma}
    \label{lemmaHartEuclidien}
    When~$\lVert c \rVert_{\Cxa{k}}$ is small enough,~$F_c$ is well-defined on all~$U$ and satisfies:
    \begin{itemize}
        \item~$F_c$ is a~$\Cxa{k}$ diffeomorphism on~$U$,
        \item~$F_c\equiv \id$ on~$F_c(U\setminus V)$,
        \item~$F_c^*Z$ is of class~$\Cxa{k}$ on~$F_c(W)$, 
        \item if~$Z$ is of class~$\Cxa{k}$ on a neighborhood of~$\phi_{[0,T]}(x)$, then~$F_c^*Z$ is of class~$\Cxa{k}$ on a neighborhood of~$F_c(x)$.
    \end{itemize}
    Additionally~$F_c$ convergences toward the identity for the~$\Cxa{k}$ topology, when~$\lVert c \rVert_{\Cxa{k}}$ goes to zero.
\end{lemma}

\begin{proof}
    The maps $x\mapsto\Psi_{c(x)s}(x)$ are $\Cxa{k}$-close to the identity map when~$\lVert c \rVert_{\Cxa{k}}$ is small. So~$\lVert F_c-\id \rVert_{\Cxa{k}}$ is converges to zero when~$\lVert c \rVert_{\Cxa{k}}$ goes to zero. Since the set of~$\Class^1$ embedding is open in the set of~$\Class^1$ map,~$F_c$ is a~$\Cxa{k}$ diffeomorphism on its image when~$T$ is small.

    Notice that~$F_c$ coincides with the identity map on~$U\setminus\wb V\supset\partial U$. So for~$T$ small, we have~$F_c(U)=U$ and~$F_c\colon U\to U$ is a~$\Cxa{k}$ diffeomorphism. For~$x$ in~$U$, we have:

    \begin{align*}
        dF_c(x) 
            &= \int_0^1d\Psi_{c(x)s}(x)ds + dc(x)\int_0^1sZ\circ\Psi_{c(x)s}(x)ds \\
        dF_c(x)(Z(x))  
            &= \int_0^1Z\circ \Psi_{c(x)s}(x)ds + (dc(x)(Z(x)))\int_0^1sZ\circ \Psi_{c(x)s}(x)ds \\
            &= \frac{1}{c(x)}(\Psi_{c(x)}(x)-x) + (dc(x)(Z(x)))\int_0^1sZ\circ \Psi_{c(x)s}(x)ds 
             \text{ when } c(x)\neq 0
    \end{align*}
    Since we have~$c(x)\neq 0$ and~$dc\equiv 0$ on~$W$,~$dF_c(Z)$ is of class~$\Cxa{k}$ on~$W$. 
    When~$Z$ is of class~$\Cxa{k}$ on $\Psi_{[0,T]}(x)$,the second member in $dF_c(x)(Z(x))$ is also~$\Cxa{k}$ on a neighborhood of~$x$. So~$dF_c(Z)$ is of class~$\Cxa{k}$ on~$x$.
\end{proof}

\begin{proof}[Proof of Lemma \ref{lemmaHartManifold}]
    Let~$(A,B_1,\cdots ,B_n)$ be an open covering of~$M$ such that $K\subset A$, $K\cap\wb{B_i}=\emptyset$, ~$Z$ is of class~$\Cxa{k}$ on a neighborhood of~$\wb A$, and such that there exists~$\Cinfty$ embedding from each~$\wb B_i$ to~$\RR^3$. We prove by induction on~$n$ that there exists a~$\Cxa{k}$ diffeomorphisms~$f\colon M\to M$ such that~$f^*Z$ is of class~$\Cxa{k}$ on~$M$. 
    
    The property is clear for $n=0$. Assume~that the property is satisfied for $n>0$. Take~$C$ an open neighborhood of~$\wb B_1$, disjoint from $K$ and~$g\colon C\to\RR^3$ be a~$\Cinfty$ embedding with pre-compact image. We take~$U=\im(g)$, and~$V,W\subset U$ such that~$\wb{g(B_1)}\subset W$,~$\wb W\subset V$ and~$\wb V\subset U$. By assumption, there exists a small $T>0$ so that $Z$ is of class $\Cxa k$ on the set $\Psi_{[0,T](\wb A)}$. Take a some $\Cxa{k}$-small function~$c\colon U\to[0,T]$ and denote by~$F_c\colon U\to U$ the~$\Cxa{k}$ diffeomorphism as in Lemma~\ref{lemmaHartEuclidien}. According to the same lemma, $(F_c\circ g)^*Z$ is of class $\Cxa k$ on $W\cup g(A\cap C)$.

    Define a map~$h\colon M\to M$ by~$h\equiv\id$ on~$M\setminus C$ and~$h\equiv g^{-1}\circ F_c\circ g$ on~$C$.
    The map~$F_c$ is equal to the identity on a neighborhood of~$\partial U$, so~$h$ is of class~$\Cxa{k}$ on all~$M$. We can take~$F_c$ arbitrarily~$\Class^1$ close to the identity so that~$F_c(W)$ contains~$\wb{g(B_1)}$. Then~$h^*Z$ is of class~$\Cxa{k}$ on a neighborhood of~$\wb{B_1}$. Similarly we can take~$F_c$ such that $h^*Z$ is~$\Cxa{k}$ on a neighborhood of~$g(\wb A\cap C)$. So~$h^*Z$ is~$\Cxa{k}$ on a neighborhood of~$\wb{A}\cap B_1$, and also on~$\wb{A}\cap (M\setminus B_1)$ by hypothesis. 

    The map~$h$ is a~$\Cxa{k}$ diffeomorphism such that~$h^*Z$ is~$\Cxa{k}$ on a neighborhood of~$\wb{A\cup B_1}$. According to the induction hypothesis, there exists a~$\Cxa{k}$ diffeomorphism~$i\colon M\to M$ such that~$i^*(h^*Z)$ is~$\Cxa{k}$ on~$M$, and with $i_{|K}=\id_K$. So~$f=i\circ h$ satisfies the property at the rank~$n$, which concludes the induction.
\end{proof}


\addcontentsline{toc}{section}{References}
\bibliographystyle{alpha}

\end{document}